\documentclass{amsart}
\usepackage{amsxtra}
\theoremstyle{plain}
\newcommand{\cleqn}{\setcounter{equation}{0}}
\newcommand{\clth}{\setcounter{theorem}{0}}
\newcommand {\sectionnew}[1]{\section{#1}\cleqn\clth}
\newtheorem{theorem}{Theorem}[section]
\newtheorem{lemma}[theorem]{Lemma}
\newtheorem{definition-theorem}[theorem]{Definition-Theorem}
\newtheorem{proposition}[theorem]{Proposition}
\newtheorem{corollary}[theorem]{Corollary}
\newtheorem{definition}[theorem]{Definition}
\newtheorem{example}[theorem]{Example}
\newtheorem{remark}[theorem]{Remark}
\newtheorem{notation}[theorem]{Notation}
\newtheorem{assumption}[theorem]{Assumption}
\newtheorem{lemma-definition}[theorem]{Lemma-Definition}
\newtheorem{lemma-notation}[theorem]{Lemma-Notation}
\newtheorem{question}[theorem]{Question}
\newtheorem{remark-definition}[theorem]{Remark-Definition}
\newtheorem{notation-remark}[theorem]{Notation-Remark}
\newtheorem{definition-remark}[theorem]{Definition-Remark}
\newcommand \bth[1] { \begin{theorem}\label{t#1} }
\newcommand \ble[1] { \begin{lemma}\label{l#1} }

\newcommand \bpr[1] { \begin{proposition}\label{p#1} }
\newcommand \bco[1] { \begin{corollary}\label{c#1} }
\newcommand \bde[1] { \begin{definition}\label{d#1}\rm }
\newcommand \bex[1] { \begin{example}\label{e#1}\rm }
\newcommand \bre[1] { \begin{remark}\label{r#1}\rm }

\newcommand \bnota[1] {\begin{notation}\label{n#1}\rm }
\newcommand \bas[1] { \begin{assumption}\label{a#1}\rm }

\newcommand \bln[1] { \begin{lemma-notation}\label{ln#1} }
\newcommand \bqu[1] { \begin{question}\label{q#1}\rm }

\newcommand {\eth} { \end{theorem} }
\newcommand {\ele} { \end{lemma} }

\newcommand {\epr} { \end{proposition} }
\newcommand {\eco} { \end{corollary} }
\newcommand {\ede} { \end{definition} }
\newcommand {\eex} { \end{example} }
\newcommand {\ere} { \end{remark} }
\newcommand {\enota} { \end{notation} }
\newcommand {\eas} {\end{assumption}}

\newcommand {\eln}{ \end{lemma-notation} }

\newcommand {\equ} {\end{question}}
\newcommand \thref[1]{Theorem \ref{t#1}}
\newcommand \leref[1]{Lemma \ref{l#1}}
\newcommand \prref[1]{Proposition \ref{p#1}}

\newcommand \deref[1]{Definition \ref{d#1}}
\newcommand \exref[1]{Example \ref{e#1}}
\newcommand \reref[1]{Remark \ref{r#1}}
\newcommand \lb[1]{\label{#1}}

\newcommand \asref[1]{Assumption \ref{a#1}}

\def \d {{\partial}}   %differentials and partials

\def \Rset {{\mathbb R}}         %mathsets
\def \Cset {{\mathbb C}}

\def \C  {{\mathcal{C}}}

\def \O  {{\mathcal{O}}}

\def \la {\lambda}

\def \ra  {\rightarrow}           %maps

\def \la {\langle}
\def \ra {\rangle}
\def \Ad { {\mathrm{Ad}} }

\def \g  {\mathfrak{g}}   % Lie algebra letters
\def \h  {\mathfrak{h}}
\def \f  {\mathfrak{f}}
\def \n  {\mathfrak{n}}
\def \m  {\mathfrak{m}}
\def \b  {\mathfrak{b}}
\def \p  {\mathfrak{p}}

\def \a  {\mathfrak{a}}

\def \u  {\mathfrak{u}}

\def \v  {\mathfrak{v}}
\def \q  {\mathfrak{q}}
\def \d  {\mathfrak{d}}

\def \t  {\mathfrak{t}}
\def \l  {\mathfrak{l}}
\def \c  {\mathfrak{c}}

\def \k {\mathfrak{k}}
\DeclareMathOperator \ad { {\mathrm{ad}} }

\newcommand{\beqa}{\begin{eqnarray*}}                     %added by Lu
\newcommand{\eeqa}{\end{eqnarray*}}
\def \hs {\hspace{.2in}}
\def \lara {\la \, , \, \ra}

\def \gog {\g \oplus \g}

\def \lara {\la \, , \, \ra}

\def \gdia {\g_{\diag}}

\def \lam {\lambda}
\def \V {\mathcal V}

\def \bfu {{\bf u}}

\def \bfv {{\bf v}}

\def \pist {\pi_{\rm st}}

\def \sG {{\scriptscriptstyle G}}
\def \sX {{\scriptscriptstyle X}}
\def \sY {{\scriptscriptstyle Y}}
\def \sZ {{\scriptscriptstyle Z}}
\def \sC {{\scriptscriptstyle C}}

\def \sM {{\scriptscriptstyle M}}

\def \piG {{\pi_{\scriptscriptstyle G}}}

\def \Cset {{\mathbb C}}
\def \lrw {\longrightarrow}
\def \Pist {\Pi_{\rm st}}

\def \Gdia {G_\diag}

\def \sL {{\scriptscriptstyle L}}

\def \sF {{\scriptscriptstyle F}}

\def \T {\mathbb T}

\def \delOO {\delta_{\scriptscriptstyle \O_+,\O_-}}

\def \diag {{\rm diag}}

\def \ot {\otimes}

\def \gotg {\g \ot \g}

\def \tilr {\tilde{r}}

\def \tilf {\tilde{\f}}
\def \OO {{\O_+ \cap \O_-}}
\def \sO {{\scriptscriptstyle \O}}
\def \rfp {r^\flat_+}
\def \rfm {r^\flat_-}
\def \tOO {\t_{\scriptscriptstyle{\O_+, \O_-}}}
\def \rBD {r_{\scriptscriptstyle{BD}}}

\def \Fbb {\mathbb F}

\def \sFbb {\scriptscriptstyle{\mathbb F}}

\def \sGQ {{\scriptscriptstyle{G/Q}}}
\def \sGP {{\scriptscriptstyle{G/P}}}

\def \sC {\scriptscriptstyle{C}}

\def \swF {\scriptscriptstyle{\widetilde{F}}}
\def \wF {\widetilde{F}}
\def \swFbb {\scriptscriptstyle{\widetilde{{\mathbb{F}}}}}
\def \wFbb {\widetilde{{\mathbb{F}}}}
\def \wPi {\widetilde{\Pi}}
\begin{document}

\setlength{\baselineskip}{1.2\baselineskip}
%%%%%%%%%%%%%%%%%%%%%%%%%%%%%%%%%%%%%%%%%%%%%%%%%%%%%%%%%%%%%%%%%%%%%%%%%%%
%%%%%%%%%%%%%%%%%%%%%%    Title    %%%%%%%%%%%%%%%%%%%%%%%%%%%%%%%%%%%%%%%%
\title[On the $T$-leaves of some Poisson structures related to products of flag varieties]
{On the $T$-leaves of some Poisson structures related to products of flag varieties}
\author{Jiang-Hua Lu}
\address{
Department of Mathematics   \\
The University of Hong Kong \\
Pokfulam Road               \\
Hong Kong}
\email{jhlu@maths.hku.hk}
\author{Victor Mouquin}
\address{
Department of Mathematics   \\
University of Toronto \\
Toronto, Canada}               
\email{mouquinv@math.toronto.edu}
\date{}
\begin{abstract} For a connected abelian Lie group $\T$ acting on a Poisson manifold $(Y, \pi)$
by Poisson isomorphisms, the $\T$-leaves of $\pi$ in $Y$ are, by definition, the orbits of the symplectic leaves of $\pi$ under $\T$, and the leaf stabilizer of a $\T$-leaf  is the subspace of the Lie algebra of $\T$ that is everywhere tangent to all the symplectic leaves in the $\T$-leaf. In this paper, we first develop a general theory on $\T$-leaves and leaf stabilizers for a class of Poisson structures defined by Lie bialgebra actions and quasitriangular $r$-matrices. We then apply the general theory
to four series of holomorphic Poisson structures on products of flag varieties and related spaces 
of a complex semi-simple Lie group $G$. 
We describe their $T$-leaf decompositions,
where $T$ is a maximal torus of $G$, in terms of {\it (open) extended Richardson varieties} and {\it extended double
Bruhat cells associated to conjugacy classes} of $G$, and we compute their leaf stabilizers and
the dimension of the symplectic leaves in each $T$-leaf.
\end{abstract}
\maketitle
%\tableofcontents
%%%%%%%%%%%%%%%%%%%%   Introduction   %%%%%%%%%%%%%%%%%%%%%%%%%%%%

\sectionnew{Introduction and statements of results}\lb{intro}

\subsection{Introduction}\lb{subsec-intro} 
A holomorphic Poisson structure on a complex manifold $Y$ is a holomorphic bi-vector field $\pi$ on $Y$ such that the bracket $\{\phi, \psi\} = \pi(d\phi \wedge d\psi)$ on the sheaf of holomorphic functions satisfies the Jacobi identity. Holomorphic Poisson structures form an important class of Poisson structures, and they have recently also  been studied in the context of generalized complex geometry and deformation theory (see \cite{Hitchin-ins, Kim:Poisson} and references therein).

A triple $(Y, \pi, \lambda)$, where $(Y, \pi)$ is a complex Poisson manifold and $\lambda$ a 
holomorphic action on $Y$ by a connected abelian complex Lie group $\T$ preserving $\pi$, is called a
{\it complex $\T$-Poisson manifold}. A complex $\T$-Poisson manifold
$(Y, \pi, \lambda)$ gives rise to a decomposition of $Y$ into the $\T$-orbits of symplectic leaves of $\pi$, also
called {\it $\T$-leaves}, which are of the form $\displaystyle \cup_{t \in \T} t\Sigma$, where
$\Sigma$ is a symplectic leaf of $\pi$ in $Y$  (see $\S$\ref{subsec-T-leaves} for the precise definition). While a complex manifold can not support non-symplectic holomorphic Poisson structures with finitely many symplectic leaves, as the degeneracy locus of such a Poisson structure, being a non-empty divisor, can not be the union of finitely many symplectic leaves, it is easy to construct examples of $\T$-Poisson manifolds with finitely many $\T$-leaves: for a complex torus $\T$, any smooth toric $\T$-variety with the zero Poisson structure is such an example. 

If $L$ is a $\T$-leaf in a $\T$-Poisson manifold $(Y, \pi, \lambda)$, the subspace $\t_{\scriptscriptstyle{L}}$ of the Lie algebra $\t$ of $\T$
which is tangent to every symplectic leaf in $L$ under the action $\lambda$ is called the {\it leaf stabilizer} of ($\lambda$ in)
$L$ (see $\S$\ref{subsec-T-leaves}). Every $\T$-leaf $L$ in $(Y, \pi, \lambda)$ 
admits a nowhere vanishing anti-canonical section, called the {\it Poisson $\T$-Pfaffian}, constructed as an exterior product of the Poisson bi-vector field $\pi$ and vector fields on $Y$ coming from any complement of $\t_{\scriptscriptstyle{L}}$ in $\t$
(see \reref{re-T-Pfaffian}), and the co-rank of $\pi$ in $L$ is equal to the codimension of $\t_{\scriptscriptstyle{L}}$ in $\t$. 
For a given $\T$-Poisson manifold $(Y, \pi, \lambda)$, it is a natural and interesting problem to determine the $\T$-leaf decomposition
of $\pi$ in $Y$ and the leaf stabilizers for the $\T$-leaves. 

Let $G$ be a connected complex semi-simple Lie group with Lie algebra $\g$, and fix a pair $(B, B_-)$ of opposite Borel subgroups of $G$ and a symmetric non-degenerate invariant bilinear form $\lara_\g$ on $\g$. The choice of 
$(B, B_-, \lara_\g)$ gives rise to a {\it standard} multiplicative holomorphic Poisson structure $\pist$ on $G$,
and the pair $(G, \pist)$ is known as a {\it standard complex semi-simple Poisson Lie group} 
(see $\S$\ref{subsec-rst} for detail).  
The Poisson structure $\pist$ is invariant under the action by the maximal torus $T = B \cap B_-$ by left translation.
It is well-known 
\cite{Hodges-Levasseur, HKKR} that the $T$-leaves of $\pist$ in $G$ are the double Bruhat cells
$G^{u, v} = BuB \cap B_-vB_-$, where $u, v \in W$, the Weyl group of $(G, T)$. Double Bruhat cells have been
studied intensively and have served as motivating examples of the 
theories of total positivity and cluster algebras (see \cite{BFZ:clusterIII, Fomin-Zelevinsky:double, Goodearl-Yakimov:PNAS}
and references therein).

The Poisson structure $\pist$ on $G$ projects to a well-defined Poisson structure, denoted by $\pi_{{\scriptscriptstyle{G/B}}}$, on the flag variety $G/B$ of $G$. The $T$-leaves of $\pi_{{\scriptscriptstyle{G/B}}}$ in $G/B$ have been shown in  \cite{Goodearl-Yakimov:GP} to be
precisely the open Richardson varieties, i.e, non-empty intersections $(BuB/B) \cap (B_-wB/B)$, where
$u, w \in W$ and $w \leq u$ in the Bruhat order on $W$.  Open Richardson varieties and their Zariski closures in $G/B$, called Richardson varieties,
have also been studied intensively 
from the points of view of geometric representation theory, combinatorics, and cluster algebras (see, for example, \cite{Billey-Coskun:singularities, B-Lak:SMT, KLS, K-Lak:Grass, Lak-Litt, Lenagan-Yakimov}). 
There are many other natural examples of $T$-Poisson manifolds associated to the Poisson Lie group
$(G, \pist)$, including the generalized flag varieties $G/P$ 
\cite{Goodearl-Yakimov:GP}, where $P$ is a parabolic subgroup of $G$,
twisted conjugacy classes in $G$ \cite{Evens-Lu:Grothendieck, Lu:dim-formula, Lu:twisted}, symmetric spaces of $G$
\cite{E-L:cplx}, the wonderful compactification of $G$ 
when $G$ is of adjoint type \cite{E-L:cplx}, and the variety of Lagrangian subalgebras \cite{E-L:cplx, Lu-Yakimov:DQ}. In these examples, the
$T$-leaves, and leaf stabilizers in some cases, have been determined by somewhat ad-hoc methods
(but see \cite{Yakimov:splitting} for the method of {\it weak splittings} in the study of $T$-leaves and symplectic leaves
for a class of Poisson structures including $\pi_{{\scriptscriptstyle{G/B}}}$ on $G/B$).
 .

In this paper,  we describe the $T$-leaves and the leaf stabilizers for four series of $T$-Poisson manifolds
associated to a standard complex semi-simple Poisson Lie group $(G, \pist)$, 
respectively denoted as 
\begin{equation}\lb{eq-series}
(F_n, \, \pi_n), \hs (\Fbb_n, \, \Pi_n), \hs (\wF_n, \, \tilde{\pi}_n), \hs (\wFbb_n, \, \wPi_n), \hs \hs n \geq 1.
\end{equation}
When $n = 1$, we have
\begin{align*}
&(F_1, \pi_1) = (G/B, \;\pi_{{\scriptscriptstyle{G/B}}}), \hs (\Fbb_1,\; \Pi_1) = 
((G \times G)/(B \times B_-), \;  \; \Pi_{{\scriptscriptstyle{(G \times G)/(B \times B_-)}}}), \\
&(\wF_1, \; \tilde{\pi}_1) = (G, \; \pist), \hs 
(\wFbb_1, \; \wPi_1) = (G \times G, \; \Pist),
\end{align*}
where $(G \times G,  \Pist)$ is the Drinfeld double Poisson Lie group of $(G, \pist)$ (see $\S$\ref{subsec-rst}), 
and $\Pi_{{\scriptscriptstyle{(G \times G)/(B \times B_-)}}}$ is the projection of $\Pist$ to $(G \times G)/(B \times B_-)$.
For $n \geq 1$, both $F_n$ and $\wF_n$ are quotient manifolds of $G^n$, and the Poisson structures $\pi_n$ and $\tilde{\pi}_n$ are
projections of the $n$-fold product Poisson structure $\pi_{\rm st}^n$ on $G^n$. Similarly, $\Fbb_n$ and $\wFbb_n$
are quotient manifolds of $(G \times G)^n$, with $\Pi_n$ and $\wPi_n$ projections of the product Poisson 
structure $\Pi_{\rm st}^n$ on $(G \times G)^n$. 
Precise definitions of the Poisson manifolds in \eqref{eq-series} are given in $\S$\ref{subsec-GBGB}.

This paper, a sequel to \cite{Lu-Mou:mixed}, is the second of a series of papers devoted to a detailed study of the 
four series of Poisson manifolds in \eqref{eq-series}. 
In \cite{Lu-Mou:mixed}, we have identified the Poisson structures in  \eqref{eq-series} as 
{\it mixed product Poisson structures 
defined by quasitriangular $r$-matrices}. In the present paper, 
we develop a general theory on $\T$-leaves and $\T$-leaf stabilizers for a class of Poisson structures 
defined by quasitriangular $r$-matrices and apply the theory to the Poisson manifolds in \eqref{eq-series}.
The general theory also provides a unified approach to many other Poisson structures such as those 
mentioned earlier from
\cite{E-L:cplx, Evens-Lu:Grothendieck, Goodearl-Yakimov:GP, Lu:dim-formula, Lu:twisted, Lu-Yakimov:DQ}
(see $\S$\ref{subsec-other}).

Our descriptions of $T$-leaves for the Poisson manifolds in \eqref{eq-series} naturally lead to what we call
{\it extended Bruhat cells},  {\it extended Richardson varieties}, and 
{\it extended double Bruhat cells associated to conjugacy classes} (see $\S$\ref{subsec-g-cell-BS} and $\S$\ref{subsec-g-C}).
In \cite{Balazs-Lu:BS}, the third in the series of papers on the Poisson manifolds in \eqref{eq-series},
we express explicitly in the so-called Bott-Samelson coordinates the 
Poisson structures $\pi_n$ on  
extended Bruhat cells in terms of the root strings and the structure constants
of the Lie algebra $\g$ of $G$. In particular, we show in \cite{Balazs-Lu:BS}
that each extended Bruhat cell of dimension $m$ gives rise to
a polynomial Poisson algebra $\Cset[z_1, \ldots, z_m]$ which is a {\it symmetric nilpotent semi-quadratic Poisson-Ore extension of ${\mathbb C}$} in  the sense of \cite[Definition 4]{Goodearl-Yakimov:PNAS}. Moreover, when the bilinear form $\lara_\g$ on $\g$ that comes into the 
definition of $\pist$ is suitably chosen, the Poisson bracket $\{z_i, z_j\}$ between any two coordinate functions is in fact a polynomial
with {\it integer} coefficients. In separate papers, we will further study 
extended Bruhat cells and extended double Bruhat cells in the
context of symplectic groupoids.

Let $\T$ be an algebraic torus. The $\T$-leaves in a $\T$-Poisson manifold $(Y, \pi, \lambda)$ are the semi-classical analogs 
of $\T$-prime ideals of a quantum algebra $A$ with rational $\T$-actions by automorphisms \cite{Brown-Goodearl:book}, and the $\T$-leaf decomposition of $(Y, \pi, \lambda)$  is the semi-classical analog of the 
Goodearl-Letzler partition of the spectrum ${\rm Spec}(A)$ of $A$ 
into tori indexed by $\T$-invariant prime ideals \cite{Goodearl-Letzler}. 
In particular, if a $\T$-invariant prime ideal $I$ in $A$ corresponds to a $\T$-leaf $L$ of
$(Y, \pi, \lam)$, the torus in the Goodearl-Letzler partition of ${\rm Spec}(A)$ indexed 
by $I$ should correspond to the quotient torus $\T/\T_L$, where 
$\T_L$ is the sub-torus of $\T$ preserving the symplectic leaves in $L$.
In the case of the Bruhat cell $BuB/B \subset G/B$ with the Poisson structure $\pi_1 = \pi_{{\scriptscriptstyle{G/B}}}$, 
where $u \in W$, one has the quantum algebra ${\mathcal{U}}_-^u$ constructed by
De Concini, Kac, and Procesi \cite{DKP:Uw} as a quantization of the algebra of regular functions on 
$BuB/B$ (see \cite{Yakimov:prime-2}), and 
the explicit correspondence between the Goodearl-Letzler partition of ${\rm Spec} ({\mathcal{U}}_-^u)$ and the
$T$-leaves of $\pi_{{\scriptscriptstyle{G/B}}}$ in   
$BuB/B$, namely the open Richardson varieties 
$(BuB/B) \cap (B_-wB/B)$, $w \leq u$, have been studied in detail in 
\cite{M-C, Yakimov:prime-1, Yakimov:prime-2, Yakimov:prime-3}. 
Similar studies for $(\wF_1, \tilde{\pi}_1) = (G, \pist)$ can be found in \cite{H-L-Toro, Joseph}.
It would thus be very interesting to 
study the quantizations of the four series of Poisson manifolds in \eqref{eq-series} (or of their Poisson submanifolds) and establish 
explicit correspondences between the Goodearl-Letzler partitions of the spectra of the quantizations
and the $T$-leaf decompositions and the leaf stabilizers described in the current paper.

As the Poisson manifolds treated in this paper are special classes of Poisson homogeneous spaces, the general
theory established in the paper can also be regarded a further development of Drinfeld's theory on Poisson homogeneous spaces \cite{dr:homog, Lu-Yakimov:DQ}. In particular, a generalization of Drinfeld's Lagrangian subalgebras associated to points in 
a Poisson homogeneous space \cite{dr:homog} plays an important role in our general theory (see \leref{le-l-yy} and 
and formulas \eqref{eq-delOO-3} and \eqref{eq-ty-0} for detail).

We now give an outline and the main results of the paper.
 
\subsection{Holomorphic Poisson structures related to flag varieties}\lb{subsec-GBGB}
If $G$ is a group and $n \geq 1$ an integer, let 
the product group $G^n$ act on itself from the right by 
\begin{equation}\lb{eq-Gn-Gn-1}
(g_1, g_2, \ldots, g_n) \cdot (h_1, h_2, \ldots, h_n) = (g_1h_1, \, h_1^{-1}g_2 h_2, \; \ldots, h_{n-1}^{-1}g_nh_n),
\hs g_j, h_j \in G. 
\end{equation}
For subgroups $Q_1, \ldots, Q_n$ of $G$ and subsets $S_1, \ldots, S_n$ of $G$ such that $S_1$ is right $Q_1$-invariant
and $S_j$ is left $Q_{j-1}$ and right $Q_j$-invariant for $ j = 2, \ldots, n$, let
$S_1 \times_{Q_1}  \cdots \times_{Q_{n-1}} S_n/Q_n$ denote the quotient of $S_1 \times \cdots \times S_n$ by the action of $Q_1 \times \cdots \times Q_n$ as a subgroup of $G^n$. 

If $(G, \piG)$ is a Poisson Lie group and if $Q_1, \ldots, Q_n$ are closed Poisson Lie subgroups of $(G, \piG)$, then  the product Poisson structure $\pi_\sG^n$ on $G^n$ projects to a well-defined Poisson structure on the quotient space 
$G \times_{Q_1} \cdots \times_{Q_{n-1}} G/Q_n$ (see \cite[$\S$7]{Lu-Mou:mixed}). Throughout the paper, if  
$G \times_{Q_1} \cdots \times_{Q_{n-1}} G/Q_n$ is denoted by $Z_n$, we will denote 
by $\pi_{\sZ_n}$ the projection of $\pi_\sG^n$ to $Z_n$ and also refer to $\pi_{\sZ_n}$ as a {\it quotient Poisson structure}.
Denote the image
$(g_1, \ldots, g_n) \in G^n$ in $Z_n$ by  
$[g_1, \ldots, g_n]_{\sZ_n}$, and define  
\[
\mu_{\sZ_n}: \;\; Z_n \lrw G/Q_n,  \;\;\; [g_1, g_2, \ldots, g_n]_{\sZ_n} \longmapsto g_1g_2 \cdots g_n Q_n \in G/Q_n.
\]
 Then 
the map $\mu_{\sZ_n}:  (Z_n, \, \pi_{\sZ_n}) \to (G/Q_n, \, \pi_{\scriptscriptstyle{G/Q_n}})$
is Poisson, and the action
\begin{equation}\lb{eq-GZZ}
G \times Z_n \lrw Z_n, \;\;\; (g, \,[g_1, g_2, \ldots, g_n]_{\sZ_n}) \longmapsto [gg_1, g_2, \ldots, g_n]_{\sZ_n}, 
\hs g, \,g_j \in G,
\end{equation}
is a Poisson action of the Poisson Lie group $(G, \piG)$ on the Poisson manifold $(Z_n, \pi_{\sZ_n})$.
This class of quotient Poisson structures was introduced in \cite{Lu-Mou:mixed}.

Let $(G, \pist)$ be a standard complex semi-simple Poisson Lie group, determined by 
the choice of a pair $(B, B_-)$ of opposite Borel subgroups of $G$ and a symmetric non-degenerate  invariant bilinear form
$\lara_\g$ on $\g$. Let $(G \times G, \, \Pist)$ be its Drinfeld double (see
$\S$\ref{subsec-rst}).
Both $B$ and $B_-$ are Poisson Lie subgroups of $(G, \pist)$, while $B \times B_-$  is a Poisson Lie subgroup of $(G \times G, \, \Pist)$. For an integer $n \geq 1$, let
\begin{align}\lb{eq-Fn}
F_n &= G \times_B \cdots \times_B G/B,\hs \;
{\mathbb F}_n  = (G \times G) \times_{(B \times B_-)} \cdots \times_{(B \times B_-)} (G \times G)/(B \times B_-),\\
\lb{eq-DXn}
\wF_n & = G \times_B \cdots \times_B G,\hs \;
\wFbb_n  = (G \times G) \times_{(B \times B_-)} \cdots \times_{(B \times B_-)} (G \times G),
\end{align}
and let $\pi_n$ and $\tilde{\pi}_n$ be the projections of $\pi_{\rm st}^n$ from $G^n$ to $F_n$ and $\wF_n$ respectively, and
let $\Pi_n$ and $\wPi_n$ be the projections of  $\Pi_{\rm st}^n$ from $(G \times G)^n$ to 
$\Fbb_n$ and $\wFbb_n$ respectively. 
The maximal torus $T = B \cap B_-$ of $G$ acts on $(F_n, \pi_n)$,  $(\Fbb_n, \Pi_n)$, $(\wF_n, \tilde{\pi}_n)$, and
$(\wFbb_n, \wPi_n)$ by  Poisson diffeomorphisms via
\begin{align}\lb{eq-T-Fn}
&t \cdot [g_1, \, g_2, \, \ldots,\, g_n]_{\sF_n} = [tg_1, \,g_2, \,\ldots, \,g_n]_{\sF_n},\\
\lb{eq-T-Fbbn}&  t \cdot [g_1, \,k_1, \,g_2, \,k_2,  \,\ldots, \, g_n, \, k_n]_{\sFbb_n} = 
[tg_1, \,tk_1, \,g_2, \,k_2,\,  \ldots, \,g_n, \,k_n]_{\sFbb_n},\\
\lb{eq-T-Dn}
&t \cdot [g_1, \, g_2, \, \ldots,\, g_n]_{\swF_n} = [tg_1, \,g_2, \,\ldots, \,g_n]_{\swF_n},\\
\lb{eq-T-Dbbn}&  t \cdot [g_1, \,k_1, \,g_2, \,k_2,  \,\ldots, \, g_n, \, k_n]_{\swFbb_n} = 
[tg_1, \,tk_1, \,g_2, \,k_2,\,  \ldots, \,g_n, \,k_n]_{\swFbb_n},
\end{align}
where $t \in T$ and $g_j, k_j \in G$ for $1 \leq j \leq n$. Let $\h$ be the Lie algebra of $T$. For $Z_n \in \{F_n, \Fbb_n, \wF_n,
\wFbb_n\}$, let $\lambda_{\sZ_n}: \h \to \V^1(Z_n)$ be the Lie algebra action of $\h$ on $Z_n$ generated by the action of $T$ on $Z_n$, and for $z \in Z_n$, define the {\it leaf stabilizer} of $\lambda_{\sZ_n}$ at $z$ to be 
\[
\t_z = \{x \in \h: \; \; \lambda_{\sZ_n}(x)(z) \in T_z \Sigma_z\},
\]
where $\Sigma_z$ is the symplectic leaf of $\pi_{\sZ_n}$ in $Z_n$ through $z$.

Let $W = N_G(T)/T$ be the Weyl group of $(G, T)$, where $N_G(T)$ is the normalizer of $T$ in $G$. Let $\leq$ be the Bruhat order on $W$, and recall the monoidal product $\ast$ on $W$ defined by 
\[
\overline{BuBvB} = \overline{B(u \ast v) B}, \hs\hs u, \, v \in W,
\]
where for a subset $X$ of $G$, $\overline{X}$ denotes the Zariski closure of $X$ in $G$. For $\bfu = (u_1,  \ldots, u_n) \in W^n$, let $l(\bfu) = l(u_1) + \cdots + l(u_n)$, where $l: W \to {\mathbb N}$ is the length function on $W$, and let 
$$
B \bfu B = (Bu_1B) \times_B \ldots \times_B (Bu_nB) \subset \wF_n.
$$
For another sequence $\bfv = (v_1,  \ldots, v_n) \in W^n$, let
\[
(B \!\times\! B_-\!)(\bfu, \!\bfv\!)(B \!\times\! B_-\!)= 
(Bu_1B \!\times\!B_- v_1 B_-\!)\times_{(\!B \!\times\!B_-\!)} \cdots \times_{(\!B \!\times\! B_-)} (Bu_nB \!\times\! B_- v_n B_-)\subset \wFbb_n.
\]
The images of $B \bfu B \subset \wF_n$ and $(B \times B_-)(\bfu, \bfv)(B \times B_-)\subset \wFbb_n$ under the projections
\begin{align*}
&\wF_n \longrightarrow F_n, \;\;\; [g_1, g_2, \ldots, g_n]_{\swF_n} \longmapsto [g_1, g_2, \ldots, g_n]_{\sF_n},\hs g_j \in G, \\
&\wFbb_n \longrightarrow \Fbb_n, \;\;\; [g_1, k_1, g_2, k_2, \ldots, g_n, k_n]_{\swFbb_n} \longmapsto 
[g_1, k_1, g_2, k_2, \ldots, g_n, k_n]_{\sFbb_n},\hs g_j, k_j \in G,
\end{align*}
will be respectively denoted by $B \bfu B/B \subset F_n$ and $(B \times B_-)(\bfu, \bfv)(B \times B_-)/(B \times B_-)\subset \Fbb_n$.

We now state our results on the $T$-leaves and leaf stabilizers for each one of the four series in \eqref{eq-series}.
%\[
%(Z_n, \pi_{\sZ_n}) \in \{(F_n, \pi_n), \; (\Fbb_n, \Pi_n), \; (\wF_n, \tilde{\pi}_n), \;(\wFbb_n, \wPi_n)\}.
%\]

The following \thref{th-Fn} on $(F_n, \pi_n)$, $n \geq 1$, will be proved in $\S$\ref{subsec-proof-AB}. 
For $n = 1$, Parts 1) and 2) of \thref{th-Fn} have been proved in \cite[Theorem 0.4]{Goodearl-Yakimov:GP} and 
\cite[Theorem 3.1]{Yakimov:prime-3}.

\bth{th-Fn} For $\bfu =(u_1, \ldots, u_n)\in W^n$ and  $w \in W$, let 
$$
R^\bfu_w = (B \bfu B/B)  \cap \mu_{\sF_n}^{-1} (B_-wB/B) \subset F_n \hs \mbox{and} \hs
\h^{\bfu}_w = \{x + u_1 \cdots u_n w^{-1}(x): x \in \h \} \subset \h.
$$

1) $R^\bfu_w\neq \emptyset$ if and only if $w\leq u_1 \ast \cdots \ast u_n$, and in this case,  $R^\bfu_w$ is a connected smooth submanifold of $F_n$ with
$\dim R^\bfu_w = l(\bfu) - l(w)$;

2) The decomposition of $F_n$ into $T$-leaves of the Poisson structure $\pi_n$ is
\[
F_n = \bigsqcup_{\bfu \in W^n, w \in W} R^\bfu_w,
\]
and all the symplectic leaves of $\pi_n$ in $R^\bfu_w$ have dimensions equal to 
$$
l(\bfu) - l(w) - \dim \ker(1 + u_1u_2 \cdots u_nw^{-1});
$$

3) The leaf stabilizer of $\lambda_{\sF_n}$ at $z \in R^\bfu_w$ is $\t_z =  \h^{\bfu}_w$. 
\eth

The following \thref{th-Fbbn} on $(\Fbb_n, \Pi_n)$, $n \geq 1$, will be proved in $\S$\ref{subsec-proof-AB}. 
For $n = 1$, \thref{th-Fbbn}  has been proved in \cite[$\S$4]{E-L:cplx}.

\bth{th-Fbbn} For $\bfu = (u_1, \ldots, u_n), \bfv = (v_1, \ldots, v_n) \in W^n$, and $w \in W$, let
\begin{align*}
&G(w) = G_{\rm diag}(w, e)(B\times B_-)/(B \times B_-) \subset (G \times G)/(B \times B_-),\\
&R^{\bfu, \bfv}_w =((B \times B_-) (\bfu, \bfv) (B \times B_-)/ (B \times B_-)) \cap \mu_{\sFbb_n}^{-1} (G(w))  \subset \Fbb_n,\\
&\h^{\bfu, \bfv}_w = \{ x + u_1 \cdots u_n w^{-1} (v_1 \cdots v_n)^{-1}(x) : x \in \h \} \subset \h,
\end{align*}
where $G_{\rm diag} = \{(g, g): g \in G\}$. 

1) $R^{\bfu, \bfv}_w\neq \emptyset$ if and only if $w \leq (v_1 \ast \cdots \ast v_n)^{-1} \ast u_1 \ast \cdots \ast u_n$, and in this case, 
$R^{\bfu, \bfv}_w$ is a connected smooth submanifold of $\Fbb_n$ with $\dim R^{\bfu, \bfv}_w = l(\bfu) + l(\bfv) - l(w)$;

2) The decomposition of $\Fbb_n$ into $T$-leaves of the Poisson structure $\Pi_n$ is
\[
\Fbb_n = \bigsqcup_{\bfu, \bfv \in W^n, w \in W} R^{\bfu, \bfv}_w,
\]
and all the symplectic leaves of $\Pi_n$ in $R^{\bfu, \bfv}_w$ have dimensions equal to 
$$
l(\bfu) + l(\bfv) - l(w) -\dim \ker (1 + u_1 \cdots u_n w^{-1}(v_1 \cdots v_n)^{-1});
$$

3) The leaf stabilizer of $\lambda_{\sFbb_n}$ at $z \in R^{\bfu, \bfv}_w$ is $\t_z =  \h^{\bfu, \bfv}_w$. 
\eth

The following \thref{th-Dn-intro} on $(\wF_n, \tilde{\pi}_n)$, $n \geq 1$, will be proved in $\S$\ref{subsec-proof-CD}. 
For $n = 1$, Parts 1) and 2) of \thref{th-Dn-intro} have been proved in \cite{Kogan-Zelevinsky}.

\bth{th-Dn-intro} 
1) For any $\bfu = (u_1, \ldots, u_n) \in W$ and $v \in W$, the intersection 
$(B \bfu B) \cap \mu_{\swF_n}^{-1}(B_- v B_-)$
in $\wF_n$ is a non-empty smooth connected submanifold of $\wF_n$ of
dimension $l(\bfu) + l(v) + \dim T$;

2) The decomposition of $\wF_n$ into $T$-leaves of the Poisson structure $\tilde{\pi}_n$ is 
\[
\wF_n = \bigsqcup_{\bfu \in W^n, \, v \in W} (B \bfu B) \cap \mu_{\swF_n}^{-1}(B_- v B_-),
\]
and all the symplectic leaves of $\tilde{\pi}_n$ in $(B \bfu B) \cap \mu_{\swF_n}^{-1}(B_- v B_-)$ have dimensions equal to 
\[
l(\bfu) + l(v) + \dim {\rm Im}(1- u_1 \cdots u_nv^{-1});
\]

3) The leaf stabilizer of $\lambda_{\swF_n}$ at $z \in (B \bfu B) \cap \mu_{\swF_n}^{-1}(B_- v B_-)$ is 
$\t_z =
\{ x - u_1 \cdots u_nv^{-1}(x): x \in \t\}$.
\eth

Let $\C$ be the set of all conjugacy classes in $G$. The following \thref{th-Dbbn} on $(\wFbb_n, \wPi_n)$, $n \geq 1$, will be proved in $\S$\ref{subsec-proof-CD}. 
For $n = 1$, \thref{th-Dbbn} has been proved in \cite[$\S$7.3]{Lu:twisted}.

\bth{th-Dbbn} For $C \in {\mathcal{C}}$, let 
$\Omega_{\sC} = \Gdia \cdot(C \times \{e\}) \cdot \Gdia = \{(g, h) \in G \times G: \; gh^{-1} \in C\}$, 
and for 
$\bfu = (u_1, \ldots, u_n), \,\bfv = (v_1, \ldots, v_n) \in W^n$, let
\begin{align*}
R_{\sC}^{\bfu, \bfv} & = ((B \!\times\! B_-\!)(\bfu, \!\bfv\!)(B \!\times\! B_-\!)) \cap \mu_{\swFbb_n}^{-1}(\Omega_{\sC})  \subset \wFbb_n,\\
\h^{\bfu, \bfv} &= \{ x - u_1 \cdots u_n(v_1 \cdots v_n)^{-1}(x): x \in \h \} \subset \h.
\end{align*}

1) For any $\bfu, \bfv \in W^n$ and $C \in \C$,
 $R_{\sC}^{\bfu, \bfv}$ is a connected smooth submanifold of
$\wFbb_n$ of dimension $l(\bfu) + l(\bfv) + \dim C + \dim T$; 

2) The decomposition of $\wFbb_n$ into $T$-leaves of the Poisson structure $\wPi_n$ is
\begin{equation}\lb{eq-T-leaves-Dbbn}
\wFbb_n = \bigsqcup_{\bfu, \bfv \in W^n, C \in \C} R^{\bfu, \bfv}_{\sC},
\end{equation}
and all symplectic leaves of $\wPi_n$ in $R^{\bfu, \bfv}_{\sC}$ have
dimensions equal to 
\[
l(\bfu) + l(\bfv) + \dim C + \dim {\rm Im}(1- u_1 \cdots u_n(v_1 \cdots v_n)^{-1});
\]

3) The leaf stabilizer of $\lambda_{\swFbb_n}$ at $z \in R_{\sC}^{\bfu, \bfv}$ is 
$\t_z = \h^{\bfu, \bfv}$. 
\eth

\subsection{Extended Bruhat cells, Bott-Samelson varieties, and extended Richardson varieties}\lb{subsec-g-cell-BS}
Let $n \geq 1$ and consider the disjoint union
\[
F_n = \bigsqcup_{{\bf u} \in W^n} B \bfu B/B.
\]
For $\bfu = (u_1, \ldots, u_n) \in W^n$, we will call
\[
B\bfu B/B = (Bu_1B \times_B \cdots \times_B Bu_n B)/B \subset F_n
\]
an {\it extended Bruhat cell}. By \thref{th-Fn}, extended Bruhat cells in $F_n$ are Poisson submanifolds
with respect to the Poisson structure $\pi_n$. 
If $\bfu = (s_1, \ldots, s_n)$, where each $s_j \in W$ is a simple reflection, we say that the 
extended Bruhat cell $B \bfu B/B$ in $F_n$ is {\it of Bott-Samelson type}. 

Every extended Bruhat cell $B \bfu B/B \subset F_n$  with the Poisson structure $\pi_n$ is Poisson isomorphic to an
 extended Bruhat cell of Bott-Samelson type in $F_{l(\bfu)}$ with the Poisson structure $\pi_{l(\bfu)}$. Indeed, 
consider first the case of $n = 1$: if $u \in W$ and $u = s_1 \cdots s_{l(u)}$ is a reduced decomposition of $u$, the Poisson morphism
$\mu_{\sF_{l(u)}}: (F_{l(u)}, \pi_{l(u)}) \to (G/B, \pi_1)$ then restricts to a  Poisson isomorphism
\[
(B {\bf s}(u) B/B, \; \pi_{l(u)})  
\lrw (B u B/B, \; \pi_1),
\]
where ${\bf s}(u) = (s_1, \ldots, s_{l(u)})$.
For any arbitrary $\bfu = (u_1, \ldots, u_n) \in W^n$, the choice of a reduced word 
${\bf s}(u_j) \in W^{l(u_j)}$ for
each $u_j$ gives rise to a Poisson isomorphism from $(B \bfu B/B, \pi_n)$ to the Bott-Samelson type extended Bruhat cell 
$B ({\bf s}(u_1), \ldots, {\bf s}(u_n)) B/B$ in  $F_{l(\bfu)}$.

On the other hand, recall that for any sequence $(s_1, \ldots, s_n)$ of simple reflections in $W^n$,
the Bott-Samelson variety $Z_{(s_1, \ldots, s_n)}$ can be defined as the Zariski closure of  
$B (s_1, \ldots, s_n) B/B$ in $F_n$ and thus carries the Poisson structure $\pi_n$. 
Bott-Samelson varieties, with the Poisson structures thus defined, 
are therefore the building blocks for the
Poisson manifolds $(F_n, \pi_n)$, $n \geq 1$, and are important even for the study of $(F_1, \pi_1)$. In \cite{Balazs-Lu:BS},
a sequel to the current paper, the Poisson structure $\pi_n$ on any Bott-Samelson variety $Z_{(s_1, \ldots, s_n)}$ 
is explicitly computed in each of the $2^n$ natural affine coordinate charts. In 
particular, it is shown in \cite{Balazs-Lu:BS} that in each of these affine coordinate charts, 
$\pi_n$ gives rise to a polynomial Poisson algebra $\Cset[z_1, \ldots, z_n]$ which is a Poisson-Ore extension of $\Cset$ compatible with 
the $T$-action. When the bilinear form $\lara_\g$ is such that $\frac{1}{2}\la \alpha, \alpha\ra_\g \in 
{\mathbb Z}$ for every root of $\g$, it is shown in \cite{Balazs-Lu:BS} that the Poisson structure 
on $\Cset[z_1, \ldots, z_n]$ corresponding to each of the $2^n$ affine coordinate charts has the property that
the Poisson brackets $\{z_i, z_k\}$
between the coordinate functions are polynomials with coefficients in ${\mathbb Z}$, thus giving rise to 
a Poisson-Ore extension of any field ${\bf k}$ of arbitrary characteristic.

For $\bfu =(u_1, \ldots, u_n) \in W^n$ and $w \in W$ such that $w \leq u_1 \ast \cdots \ast u_n$, it is natural to 
refer to 
\[
R^{\bfu}_w = (B \bfu B/B) \cap \mu_{\sF_n}^{-1}(B_-wB/B) \subset F_n
\]
as an (open) {\it extended Richardson variety} and its closure $\overline{R^{\bfu}_w}$ in $F_n$ an 
{\it extended Richardson variety} in $F_n$. By \thref{th-Fn}, extended Richardson varieties in $F_n$
are precisely closures of $T$-leaves of the Poisson structure $\pi_n$ in $F_n$. 
It would be very interesting to extend the work of 
Lenagan and Yakimov in \cite{Lenagan-Yakimov} on cluster structures on Richardson varieties in $F_1=G/B$ to extended Richardson varieties.

Analogous to taking Zariski closures in $F_n$ of extended  Bruhat cells of Bott-Samelson type,
 one can consider the Zariski closures 
in $\Fbb_n$ of $(B \times B_-)(\bfu, \bfv) (B \times B_-)/(B \times B_-)$, where $\bfu, \bfv \in W^n$ are
sequences of simple reflections. Such closures  are examples of 
{\it double Bott-Samelson varieties} introduced in \cite{Victor:thesis}, and carry the Poisson structure $\Pi_n$.
Some combinatorial aspects of double Bott-Samelson varieties and calculations of the Poisson structure $\Pi_n$ 
in special coordinate charts 
have been given in \cite{Victor:IMRN,Victor:thesis}. 

\subsection{Extended double Bruhat cells associated to conjugacy classes}\lb{subsec-g-C} 
Analogous to the Poisson manifold $(\wF_n = G \times_B \cdots \times_B G, \, \tilde{\pi}_n)$, $n \geq 1$, one also has the quotient manifold
\[
\wF_{-n} = G \times_{B_-} \cdots \times_{B_-} G
\]
of $G^n$ with the well-defined Poisson structure $\tilde{\pi}_{-n}$, the projection of the Poisson structure 
$\pi_{\rm st}^n$ from $G^n$ to $\wF_{-n}$. 
Consider, on the other hand, the diffeomorphism
$S_{\swFbb_n}: \wFbb_n \to \wF_n \times \wF_{-n}$ given by
\[
[g_1, k_1, \ldots, g_n, k_n]_{\swFbb_n} \longmapsto ([g_1, \ldots, g_n]_{\swF_n}, \, [k_1, \ldots, k_n]_{\swF_{-n}}),
\]
and set $\tilde{\pi}_{n, n} = S_{\swFbb_n}(\wPi_n)$. By a result in \cite[$\S$8]{Lu-Mou:mixed}, 
$\tilde{\pi}_{n, n}$  is 
a {\it two-fold mixed product Poisson
structure} on the product manifold $\wF_n \times \wF_{-n}$, i.e., it is the sum of the product Poisson structure $\tilde{\pi}_n \times \tilde{\pi}_{-n}$ and a certain mixed term
defined by the action of $B$ on $\wF_n$ and $B_-$ on $\wF_{-n}$ by left translation.  
For $(\bfu, \bfv) =(u_1, \ldots, u_n, v_1, \ldots, v_n) \in W^{2n}$ and any conjugacy class $C$ in $G$, define
\begin{align*}
G_{\sC}^{\bfu, \bfv} &= \{([g_1, \ldots, g_n]_{\swF_n}, \, [k_1, \ldots, k_n]_{\swF_{-n}})
\in (B \bfu B) \times (B_- \bfv B_-): \;
g_1g_2\cdots g_n (k_1k_2 \cdots k_n)^{-1} \in C\}\\
&  \subset \wF_n \times \wF_{-n},
\end{align*}
where $(B_- \bfv B_-) = (B_-v_1 B_-) \times_{B_-} \cdots \times_{B_-} (B_-v_n B_-) \subset \wF_{-n}$.
By \thref{th-Dbbn},  each $G^{\bfu, \bfv}_{\sC}$ is a $T$-leaf of $\tilde{\pi}_{n, n}$ in $\wF_n \times \wF_{-n}$ for the diagonal action of $T$.
We call $G^{\bfu, \bfv}_{\sC}$ an {\it extended double Bruhat cell associated to the conjugacy class $C$}. 
The case of $n = 1$ has been considered in \cite[$\S$7.3]{Lu:twisted}, where 
$G^{u, v}_{\sC} \subset G$, for $u, v \in W$ and a conjugacy class $C$ in $G$, is called a {\it double Bruhat cell
associated to the conjugacy class $C$}.
Note that for $\bfu, \bfv \in W^n$ and  
$C = \{e\}$, we have
\[
G^{\bfu, \bfv} :=G^{\bfu, \bfv}_{\{e\}} = 
\{([g_1, \ldots, g_n]_{\swF_n}, \, [k_1, \ldots, k_n]_{\swF_{-n}})\in (B \bfu B) \times (B_- \bfv B_-): \;
g_1g_2\cdots g_n =k_1k_2 \cdots k_n\},
\]
a direct generalization of double Bruhat cells in $G$. In a separate paper, we will study
extended Bruhat cells and extended double Bruhat cells via Poisson groupoids and symplectic groupoids.

\subsection{General theory}\lb{subsec-summary-general}
Let $r \in \gotg$ be a
quasitriangular $r$-matrix on a Lie algebra $\g$ (definition recalled in $\S$\ref{subsec-r}), and let $\lambda: \g \to \V^1(Y)$
be a Lie algebra action of $\g$ on a manifold $Y$. A simple observation made in \cite{Lu-Mou:mixed} (see also 
\cite[$\S$2.1]{David-Severa:quasi-Hamiltonian-groupoids}) is that if the $2$-tensor field $\lambda(r)$ on $Y$ is skew-symmetric, then it is Poisson, and in this case we say that the Poisson structure $-\lambda(r)$ (or $\lambda(r)$)
on $Y$ is {\it defined} by the Lie algebra action $\lam$ and the quasitriangular $r$-matrix $r$.
We also refer to Poisson structures obtained this way simply as {\it Poisson structures defined by quasitriangular $r$-matrices}.

Consider again the quotient manifold $Z_n = G \times_{Q_1} \cdots \times_{Q_{n-1}} G/Q_n$ with the quotient Poisson
structure $\pi_{\sZ_n}$ 
for an arbitrary Poisson Lie group $(G, \piG)$ and
closed Poisson Lie subgroups $Q_1, \ldots, Q_n$ described in $\S$\ref{subsec-GBGB}. Note that the manifold $Z_n$ is diffeomorphic to the product manifold
$G/Q_1 \times \cdots \times G/Q_n$ via the diffeomorphism
$J_{\sZ_n}: Z_n \to G/Q_1 \times \cdots \times G/Q_n$ given by
\begin{equation}\lb{eq-JZ-intro}
J_{\sZ_n}([g_1, g_2, \ldots, g_n]_{\sZ_n}) = (g_1Q_1, \; g_1g_2Q_2, \; \ldots, \;g_1g_2\cdots g_nQ_n),  \hs 
g_1, \ldots, g_n  \in G.
\end{equation}
A key fact established in \cite[$\S$7]{Lu-Mou:mixed} is that the Poisson structure $J_{\sZ_n}(\pi_{\sZ_n})$ on 
$G/Q_1 \times \cdots \times G/Q_n$ is defined by a quasitriangular $r$-matrix
(see $\S$\ref{subsec-two-main}).

After a review on the notion of $\T$-leaves in $\S$\ref{sec-T-leaves} and on Poisson structures defined by
quasitriangular $r$-matrices in $\S$\ref{sec-r}, we devote 
$\S$\ref{sec-orbits} - $\S$\ref{sec-mixed} of the paper  to a general theory on
$\T$-leaves and $\T$-leaf stabilizers for a class of $\T$-invariant Poisson structures defined by
quasitriangular $r$-matrices. \thref{th-Fn} - \thref{th-Dbbn} are  proved  
in $\S$\ref{sec-flags} as immediate 
examples of the general theory. 

To give an outline of the general theory, let $G$ be a connected Lie group with Lie algebra $\g$ and $r \in \g \otimes \g$ a factorizable 
quasitriangular
$r$-matrix on $\g$. The symmetric part of $r$ determines 
a symmetric non-degenerate invariant bilinear form $\lara_\g$ on $\g$ (see $\S$\ref{subsec-factorizable}).  If $Y$ is a manifold with a Lie group action $\lambda$ by $G$, we
say that the quadruple $(G, r, Y, \lam)$ is {\it strongly admissible} if the stabilizer subgroup
$Q_y$ of $G$ at every $y \in Y$ is connected and its Lie algebra $\q_y$ satisfies
\[
[\q_y, \q_y]\subset \q_y^\perp \subset \q_y,
\]
where 
$\v^\perp = \{x \in \g: \la x, \v\ra_\g=0\}$ for a subspace $\v$ of $\g$.
If $(G, r, Y, \lam)$ is strongly admissible, the $2$-tensor field $\lam(r)$ on $Y$ is necessarily skew-symmetric
(see $\S$\ref{subsec-sigma-r}), so it is
a Poisson structure on $Y$. 
Associated to the factorizable quasitriangular $r$-matrix $r$ one also has  
two distinguished 
Lie subalgebras  $\f_+ = {\rm Im} (\rfp)$ and $\f_- = {\rm Im}(\rfm)$  of $\g$, where if $r = \sum_i x_i \otimes x_i^\prime \in \gotg$, then 
$\rfp,  \rfm: \g \to \g$ are respectively given by 
\[
\rfp(x) = \sum_i \la x, \, x_i\ra_\g x_i^\prime \hs \mbox{and} \hs
\rfm(x) = -\sum_i \la x, \, x_i^\prime\ra_\g x_i, \hs x \in \g.
\]
 A pair $(M_+, M_-)$ of connected Lie subgroups of $G$ 
is said to be {\it $r$-admissible} if their respective Lie algebras $\m_+$ and $\m_-$
satisfy 
\[
\f_+ \subset \m_+, \hs \f_- \subset \m_-, \hs [\m_+, \, \m_+] \subset \m_+^\perp, \hs [\m_-, \, \m_-] \subset \m_-^\perp.
\]

In $\S$\ref{sec-orbits} we consider a six-tuple
$(G, r, Y, \lam, M_+, M_-)$, where $(G, r, Y, \lam)$ is a strongly admissible quadruple, and 
$(M_+, M_-)$ is a pair of $r$-admissible Lie subgroups of $G$.  
Let $\T$ be the connected component of $M_+ \cap M_-$ containing the identity element. Then $\T$ is necessarily abelian
and acts on $(Y, \lam(r))$ through $\lam$ by Poisson isomorphisms. On the other hand, one has the disjoint union
\begin{equation}\lb{eq-YOO}
Y = \bigsqcup_{\O_+, \O_-} \OO,
\end{equation}
where $\O_+$ and $\O_-$ are respectively $M_+$-orbits and $M_-$-orbits in $Y$.
The conditions $\f_+ \subset \m_+$ and $\f_- \subset \m_-$ imply that each non-empty intersection $\OO$ in 
\eqref{eq-YOO} is a $\T$-invariant Poisson submanifold of the $\T$-Poisson manifold $(Y, \lam(r), \lam)$.
We say that the six-tuple $(G, r, Y, \lam, M_+, M_-)$ is {\it admissible} if
the $\T$-leaves of the Poisson structure $\lam(r)$ in $Y$ are precisely the connected components of all the
non-empty intersections $\OO$ in \eqref{eq-YOO}.

For each pair $(\O_+, \O_-)$ of $M_+$- and $M_-$-orbits contained in the same $G$-orbit in $Y$, 
we introduce an integer $\delOO$, given in \eqref{eq-delOO-3}, and a 
subspace $\tOO$ of the Lie algebra $\t = \m_+ \cap \m_-$ of $\T$, given in \eqref{eq-ty-0}. Both $\delOO$ and $\tOO$
are defined using arbitrary points $y_+ \in \O_+$ and $y_- \in \O_-$ but are independent of the choices.
The main results of $\S$\ref{sec-orbits} are summarized in the following \thref{th-main}.

\bth{th-main}
Suppose that $\delOO =0$ for every pair $(\O_+, \O_-)$ of $M_+$- and $M_-$-orbits. Then the
six-tuple $(G, r, Y, \lam, M_+, M_-)$ is admissible;  the leaf stabilizer of each $\T$-leaf in $\OO$ is $\tOO$, and the co-rank of the Poisson structure $\lambda(r)$ in $\OO$ is equal to the co-dimension of $\tOO$ in $\t$.
\eth

The second part of the general theory, presented in $\S$\ref{sec-mixed}, is an application of the ``test" described in
\thref{th-main}: assume that $Q$ is a closed and connected Lie subgroup of $G$ whose Lie algebra satisfies 
$\f_+ \subset \q$ and $[\q, \q] \subset \q^\perp$. Equip $G$ with the Poisson structure $\piG = r^L - r^R$, where
$r^L$ (resp. $r^R$) is the
left (resp. right) invariant tensor field on $G$ with value $r$ at the identity element of $G$. Then
\cite[$\S$7]{Lu-Mou:mixed}  $Q$ is a Poisson Lie subgroup of
the Poisson Lie group $(G, \piG)$, and one thus has the quotient Poisson structure $\pi_{\sZ_n}$ on the quotient manifold
$Z_n = G\times_Q \cdots \times_Q G/Q$, as explained in $\S$\ref{subsec-GBGB}. Let $(M_+, M_-)$ be again an $r$-admissible pair of Lie subgroups of $G$ with
respective Lie algebras $\m_+$ and $\m_-$. The connected component $\T$ of $M_+ \cap M_-$ containing the identity element
then acts on $(Z_n, \pi_{\sZ_n})$ by Poisson isomorphisms via \eqref{eq-GZZ}.
We study the $\T$-leaves of $\pi_{\sZ_n}$ in $Z_n$ via the $\T$-leaves of the Poisson structure $J_{\sZ_n}(\pi_{\sZ_n})$ 
on the $n$-fold product manifold $(G/Q)^n$, where $J_{\sZ_n}: Z_n \to (G/Q)^n$ is the diffeomorphism
given in \eqref{eq-JZ-intro}, and $\T$ acts on $(G/Q)^n$ diagonally. 

More precisely, let $\lam$ be the Lie algebra action of the direct product Lie algebra $\g^n$ on $(G/Q)^n$ 
induced from the left action of $G^n$ on $(G/Q)^n$ by left translation.
By a result in \cite[$\S$8]{Lu-Mou:mixed}
(see $\S$\ref{subsec-two-main} for more detail), the Poisson structure $J_{\sZ_n}(\pi_{\sZ_n})$ 
on $(G/Q)^n$ is defined by the Lie algebra action $\lambda$ of $\g^n$ and 
a certain quasitriangular $r$-matrix $r^{(n)}$ on  $\g^n$. One can thus study the $\T$-leaves of $J_{\sZ_n}(\pi_{\sZ_n})$ in 
$(G/Q)^n$ via the six-tuple
\begin{equation}\lb{eq-six-intro}
\left(G^n, \; r^{(n)}, \; (G/Q)^n, \; \lambda,\;  M_+^{(n)}, \; M_-^{(n)}\right),
\end{equation}
where $(M_+^{(n)}, M_-^{(n)})$ is an $r^{(n)}$-admissible pair of Lie subgroups of $G^n$ determined by $(M_+, M_-)$.
Applying \thref{th-main}, we obtain in \prref{pr-D} sufficient conditions for the
six-tuple in \eqref{eq-six-intro} to be admissible. Translating to $(Z_n, \pi_{\sZ_n})$ using 
the Poisson diffeomorphism 
\[
J_{\sZ_n}: \;\; (Z_n, \,\pi_{\sZ_n}) \longrightarrow ((G/Q)^n, \,J_{\sZ_n}(\pi_{\sZ_n})),
\]
we obtain in \thref{th-Zn} a description of $\T$-leaves and leaf stabilizers for $(Z_n, \pi_{\sZ_n})$ under some sufficient conditions on the triple $(M_+, M_-, Q)$.   
\thref{th-Fn} and \thref{th-Fbbn} are then proved in $\S$\ref{subsec-proof-AB} as special cases of
\thref{th-Zn}.
\thref{th-Dn-intro} and \thref{th-Dbbn} are similarly 
proved in $\S$\ref{subsec-proof-CD} as special cases of \thref{th-Dn}, an analog of \thref{th-Zn}.

\subsection{Notation}\lb{subsec-nota-intro}
Throughout the paper, the pairing between a finite dimensional vector space $V$ (over $\Rset$ or $\Cset$)
and its dual space $V^*$ will always be denoted by $\lara$. 
The annihilator of a vector subspace $U$ of $V$ is, by definition, the subspace of $V^*$ given by
$U^0 =\{\xi \in V^*: \la \xi, U\ra =0\}$. 
For each integer $k \geq 1$, $\wedge^k V$ is identified with the subspace of skew-symmetric elements in $V^{\otimes k}$, and for $v_1, \ldots,  v_k \in V$,
\begin{equation}\lb{eq-wedge}
v_1 \wedge v_2 \wedge \cdots \wedge v_k = \sum_{\lambda \in S_k} {\rm sign}(\lambda) v_{\lambda(1)} \otimes v_{\lambda(2)} \otimes \cdots \otimes v_{\lambda(k)} \in \wedge^k V \subset V^{\otimes k}.
\end{equation}
For an element $r = \sum_i u_i \otimes v_i \in V \otimes V$, define $r^{21} = \sum_i v_i \otimes u_i \in V \otimes V$ and
\begin{equation}\lb{eq-r-sharp}
r^\#: \;\;\; V^* \lrw V, \;\;\; r^\#(\xi) = \sum_i \la \xi, u_i\ra v_i, \hs \xi \in V^*.
\end{equation}
Then $(r^\#)^* = (r^{21})^\#: V^* \to V$.

If $\g$ is a Lie algebra over $\Rset$ (resp. $\Cset$), by a {\it left} Lie algebra action of $\g$ on a manifold 
(resp. complex manifold) $Y$  we mean a Lie algebra anti-homomorphism
$\lambda: \g \to \V^1(Y)$, where $\V^1(Y)$ is the space of smooth (resp. holomorphic) vector fields on $Y$. 
For $k \geq 1$ and $X = \sum x_{i_1} \otimes \cdots \otimes x_{i_k} \in \g^{\otimes k}$, let
$\lam(X)$ be the $k$-tensor field on $Y$ given by
\[
\lambda(X) = \sum \lam(x_{i_1}) \otimes \cdots \otimes \lam(x_{i_k}).
\]
When
$G$ is a connected Lie group with Lie algebra $\g$ and $\lambda: G \times Y \to Y, (g, y) \mapsto gy,$ is a left Lie group action, we use
$\lambda$ to also denote the induced left action of $\g$ on $Y$, i.e.,
\begin{equation}\lb{eq-sigma-induced}
\lambda: \;\; \g \lrw \V^1(Y), \;\;\; \lambda(x)(y) = \frac{d}{dt}|_{t=0} \exp (tx)y, 
\hs x \in \g, \; y \in Y.
\end{equation}

\subsection{Acknowledgments}
We would like to thank Allen Knutson to whom goes the credit of \leref{le-leaf-stabilizer}. Work in this paper has been partially
supported by the Research Grants Council of the Hong Kong SAR, China (GRF HKU 704310 and 703712).

\sectionnew{Poisson actions and $\T$-leaves}\lb{sec-T-leaves}
\subsection{Regular and full Poisson actions by Lie bialgebras}\lb{subsec-full-poi}
Recall that a {\it Lie bialgebra} over the field ${\bf k} = \Cset$ or $\Rset$ is a pair $(\g, \delta)$, where $\g$ is a Lie algebra over
${\bf k}$ and $\delta: \g \to \wedge^2 \g$, called the {\it co-bracket}, is a linear map satisfying the co-Jacobi identity and the cocycle condition
\begin{equation}\lb{eq-cocycle}
\delta [x,\, y] = \ad_x(\delta(y)) -\ad_y (\delta(x)), \hs x, y \in \g.
\end{equation}
Given a Lie bialgebra $(\g, \delta_\g)$, the dual map of $\delta$ is then a Lie bracket on the dual space $\g^*$ of $\g$.

A {\it left Poisson action} of the Lie bialgebra $(\g, \delta)$ on a Poisson manifold $(Y, \pi)$ is a left Lie algebra
action $\lambda: \g \to \V^1(Y)$ such that
\[
L_{\lambda(x)} \pi = \lambda(\delta(x)), \hs x \in \g,
\]
where for $x \in \g$, $L_{\lambda(x)} \pi$ is the Lie derivative of $\pi$ in the direction of the vector field $\lambda(x)$.
Here when $\g$ is a Lie bialgebra over $\Cset$, we understand  $(Y, \pi)$ to be a complex Poisson manifold, i.e., 
$Y$ is a complex manifold and $\pi$ a holomorphic Poisson structure on $Y$. 

Lie bialgebras and Poisson actions of Lie bialgebras are the infinitesimal counterparts of Poisson Lie groups and Poisson
actions of Poisson Lie groups, and we refer to \cite{Chari-Pressley, dr:Poi-Lie, dr:quantum, Etingof-Schiffmann, STS2, STS:RIMS} for the basics of the theory. In this paper we follow the 
notation and sign conventions in the review section $\S$2 of \cite{Lu-Mou:mixed}.  

Let $\lambda: \g \to \V^1(Y)$ be a left Poisson action of a Lie bialgebra $(\g, \delta)$ on a Poisson manifold $(Y, \pi)$.
For $y \in Y$, set $\pi_y = \pi(y) \in \wedge^2 T_yY$, and 
\[
\lambda_y: \;\;\; \g \lrw T_yY: \;\;\; \lambda_y(x) = \lambda(x)(y), \hs y \in Y, \; x \in \g.
\]
Recall that ${\rm Im}(\pi_y^\#)$, where $\pi_y^\#: T^*_yY \to T_yY$ is defined using
\eqref{eq-r-sharp}, is the tangent space to the symplectic leaf of $\pi$ through $y$, and the rank of $\pi$ 
at $y$ is, by definition, $\dim ({\rm Im}(\pi_y^\#))$. 
Set
%and we call $\g_y \subset \g$ the {\it leaf stablizer} of the Poisson action $\lambda$ at $y \in Y$.
\begin{equation}\lb{eq-Sigma-y}
\Phi_y = \lambda_y(\g) + {\rm Im} (\pi_y^\#) \subset T_yY.
\end{equation}

\bde{de-full}
The Poisson action $\lambda$ of $(\g,\delta)$ on $(Y, \pi)$ is said to be {\it regular} (resp. {\it full}) if
$\dim \Phi_y$ is independent of $y \in Y$ (resp. $\Phi_y = T_yY$ for all $y \in Y$).
\ede

Given a left Poisson action $\lambda$ of a Lie algebra $(\g, \delta)$ on a Poisson manifold $(Y, \pi)$, it is shown in
\cite{Lu:Duke-Poi} that the vector bundle $(Y \times \g) \oplus T^*Y$, the direct product of
the trivial vector bundle over $Y$ with fiber $\g$ and the cotangent bundle $T^*Y$ of $Y$, has the 
structure of a Lie algebroid over $Y$ with $-\lambda + \pi^\#$ as the anchor map. We denote this Lie algebroid by
\begin{equation}\lb{eq-A-sigma}
A_\lambda = (Y \times \g) \bowtie T^*Y
\end{equation}
and refer to \cite{Lu:Duke-Poi} for details.
Thus the Poisson action $\lambda$ is regular (resp. full) if and only if the Lie algebroid $A_\lambda$ is regular
(resp. transitive), i.e., the anchor map $-\lambda + \pi^\#$ of $A_\lam$ has constant rank (resp. everywhere surjective).

For a Poisson action $\lambda$ of a Lie bialgebra $(\g, \delta)$ on $(Y, \pi)$, and for $y \in Y$, we also define
\begin{equation}\lb{eq-gy}
\g_y = \{x \in \g: \; \lambda_y(x) \in {\rm Im}(\pi_y^\#)\} \subset \g.
\end{equation}
The linear map $\lambda_y: \g \to T_yY$ then induces a well-defined vector space isomorphism
\begin{equation}\lb{eq-g-gy}
\g/\g_y \lrw \Phi_y/{\rm Im} (\pi_y^\#): \;\;\; x + \g_y \longmapsto \lambda_y(x) + {\rm Im}(\pi_y^\#), 
\hs x \in \g.
\end{equation}

\subsection{$\T$-leaves  and leaf stabilizers}\lb{subsec-T-leaves}
In $\S$\ref{subsec-T-leaves}, let $(Y, \pi)$ be a real (resp. complex) Poisson manifold and
$\T$  a connected real (resp. complex) abelian
Lie group
acting on $Y$ by Poisson isomorphisms. Denote the action by $\lambda: \T \times Y \to Y$ and the Lie algebra of $\T$ by $\t$.
Then $\lambda$ is a Poisson action of the Lie bialgebra $(\t, 0)$ on $(Y, \pi)$, where $0$ denotes the zero map
$\t \to \wedge^2 \t$. We will refer to the triple $(Y, \pi, \lambda)$ as a $\T$-Poisson manifold.
Let $A_\lambda = (Y \times \t) \bowtie T^*Y$  be the Lie algebroid of over $Y$ given in \eqref{eq-A-sigma}.

Recall
(see, for example, \cite{Rui:classes, Kirill:book}) that
the {\it orbits} of the Lie algebroid $A_\lambda$ are the integral submanifolds of the distribution on $Y$ 
(not necessarily of constant rank) defined as the image of 
the anchor map of $A_\lambda$, i.e., of the distribution 
$\displaystyle \bigsqcup_{y \in Y} (\lambda_y(\t) + {\rm Im}(\pi_y^\#))$ on $Y$.

\bde{de-T-leaves}
Orbits of the Lie algebroid $A_\lambda$ are called the $\T$-orbits of symplectic leaves, or simply the $\T$-leaves, of 
$\pi$ in $Y$.
\ede

The following \leref{le-T-leaves-0} justifies the terminology in \deref{de-T-leaves}.

\ble{le-T-leaves-0}
Let $L$ be any $\T$-leaf of $\pi$ in $Y$ and let $\Sigma$ be any symplectic leaf of $\pi$ such that $L \cap \Sigma \neq \emptyset$. Then $\Sigma \subset L$ and the action map
$\lambda: \T \times \Sigma \to L$
is a surjective submersion.
\ele

\begin{proof} As $\displaystyle \bigsqcup_{y \in Y} {\rm Im}(\pi_y^\#)$ is a sub-distribution of 
$\displaystyle \bigsqcup_{y \in Y} (\lambda_y(\t) + {\rm Im}(\pi_y^\#))$, one has $\Sigma \subset L$.
Let $t \in \T$ and $y \in \Sigma$. The differential of $\lambda$ at $(t, y)$ is the linear map
\[
\lambda_*: \; T_{(t, y)} (\T \times \Sigma) \lrw T_{ty} L, \; (r_t(x), v_y) \longmapsto \lambda_{ty}(x) + (\lambda_t)_* v_y, \;\;
x \in \t, \, v_y \in T_y \Sigma,
\]
where for $x \in \t$, $r_t (x) \in T_t \T$ is the right translate of $x$ by $t$, and $(\lambda_t)_*: T_y Y \to T_{ty}Y$ is the differential at $y$ of the map $\lambda_t: Y \to Y, y_1 \to ty_1, y_1 \in Y$. Since $\lambda_t$ preserves $\pi$, $(\lambda_t)_*(T_y\Sigma) = {\rm Im} (\pi_{ty}^\#)$. Thus  the action  map $\lambda$ is 
a submersion, and $\displaystyle \T\Sigma:=\bigsqcup_{t \in \T} t\Sigma$ is open in $L$. If
$\Sigma'$ is an other symplectic leaf of $\pi$ contained in $L$, then either $\T \Sigma = \T \Sigma'$ or $\T \Sigma \cap \T \Sigma' = \emptyset$, and since $L$ is connected, one must have $\T \Sigma = \T \Sigma'$. Thus $\lambda$ is surjective.
\end{proof}

\bde{de-leaf-stabilizer}
For $y \in Y$, the subspace 
\[
\t_y = \{x \in \t: \; \lambda_y(x) \in {\rm Im}(\pi_y^\#)\}
\]
of $\t$ is called the $\T$-leaf stabilizer (or simply the leaf stabilizer) of $\lambda$ at $y$.
\ede

The following \leref{le-leaf-stabilizer} is due to A. Knutson through private communication.

\ble{le-leaf-stabilizer}
If $y$ and $y'$ are in the same $\T$-leaf of $(Y, \pi)$, then $\t_y = \t_{y'}$.
\ele

\begin{proof}
Let $L$ be any $\T$-leaf of  $\pi$ in $Y$. Let 
 $2r$ be the rank of $\pi$ in $L$, and let $\pi^r$ be the $r$'th exterior product of $\pi$ with itself. For $x \in \t$, let $V_x = \lambda(x) \wedge \pi^r$, a $(2r +1)$-vector field on $L$. As $\T$ is abelian and acts on $L$ by Poisson isomorphisms, $V_x$ is $\T$-invariant. For any local (smooth or holomorphic) function $f$ on $L$, if $X_f = \pi^\#(df)$ is the Hamiltonian vector field of $f$, then 
$$
L_{X_f}V_x = (L_{X_f} \lambda(x)) \wedge \pi^r = X_{\lambda(x)(f)} \wedge \pi^r = 0.
$$
Thus $V_x$ is also invariant under local Hamiltonian flows. It follows that the vanishing locus $Z(V_x)$ of $V_x$ is open in $L$. As $Z(V_x)$ is also closed and as $L$ is connected, $V_x$ vanishes everywhere on $L$ if it vanishes at one point. But for any $y \in L$, $V_x(y) = 0$ if and only if $x \in \t_{y}$. It follows that the subspace $\t_{y}$ of $\t$ is independent of $y \in L$. 
\end{proof}

\bde{de-leaf-stab-L}
For a $\T$-leaf $L$ in $Y$, the subspace $\t_\sL :=\t_y$ of $\t$, where $y \in L$ is arbitrary, will be called the
{\it leaf stabilizer of $\lambda$ in $L$}. 
\ede

\bre{re-T-Pfaffian}
(1) Note that by the vector space isomorphism in 
\eqref{eq-g-gy}, the co-rank of $\pi$ in $L$ is the same as the co-dimension of
$\t_\sL$ in $\t$;

(2) Let $L \subset Y$ be a $\T$-leaf in $Y$ with $\dim L = l$, and let 
$2r$ be the rank of $\pi$ in $L$. Let  $\t_\sL^\prime \subset \t$ be any vector space complement of $\t_\sL$ in $\t$. Then 
$\dim \t_\sL^\prime = l -2r$. For any
$\xi \in \wedge^{l-2r} \t_\sL^\prime$, $\xi \neq 0$, the $l$-vector field
\begin{equation}\lb{eq-eta-L}
\eta_\sL = \lambda(\xi) \wedge \underbrace{\pi \wedge \cdots \wedge \pi}_r
\end{equation}
on $Y$ then restricts to a nowhere vanishing anti-canonical section of $L$. It is also clear that 
up to non-zero scalar multiples, the restriction of $\eta_\sL$ to $L$ is independent of the
choices of 
the complement $\t_\sL^\prime$ of $\t_\sL$ in $\t$ and of the non-zero element
$\xi \in \wedge^{l-2r} \t_\sL^\prime$.  
The restriction to $L$ of $\eta_\sL$ in \eqref{eq-eta-L} is called a Poisson {\it $\T$-Pfaffian} of the Poisson structure
$\pi$ on $L$, a term suggested by M. Yakimov and A. Knutson. \hfill$\diamond$
\ere

\sectionnew{Poisson structures defined by quasitriangular $r$-matrices}\lb{sec-r}

\subsection{Quasitriangular $r$-matrices}\lb{subsec-r} We recall some basic facts on quasitriangular $r$-matrices and Lie bialgebras. Our references are
\cite{Chari-Pressley, dr:Poi-Lie, dr:quantum, Etingof-Schiffmann, Lu-Mou:mixed}.

Let $\g$ be a Lie algebra. Recall that the Classical Yang-Baxter operator
 for $\g$ is the map ${\rm CYB}: \gotg \to \g \ot \g \ot \g$ given by 
\[
{\rm CYB}(r) = \sum_{i, j} \left([x_i, x_j] \ot y_i \ot y_j + x_i \ot [y_i, x_j] \ot y_j + x_i \ot x_j \ot [y_i, y_j]\right),
\hs \mbox{if} \; r = \sum_i x_i \ot y_i.
\]
A {\it quasitriangular $r$-matrix} on $\g$ is
an element $r \in \gotg$ such that $r+ r^{21} \in (S^2\g)^\g$ and that ${\rm CYB}(r) = 0$, where
$(S^2 \g)^\g$ is the space of $\g$-invariant elements in $S^2\g$ with respect to the adjoint action.
A quasitriangular $r$-matrix $r$ on $\g$ is said to be {\it factorizable} if $r+ r^{21} \in (S^2\g)^\g$ is
non-degenerate.

Let $r \in \gotg$ be a quasitriangular $r$-matrix on $\g$. Define
\begin{equation}\lb{eq-delta-r}
\delta_r: \;\; \g \lrw \g \ot \g, \;\;\; \delta_r(x) = \ad_x(r), \hs x \in \g.
\end{equation}
As $r+ r^{21} \in (S^2\g)^\g$, $\delta_r$ takes values in $\wedge^2 \g$, and 
$(\g, \delta_r)$ is a Lie bialgebra. Define $r_{\pm}: \g^* \to \g$ by
\begin{equation}\lb{eq-r-pm}
r_+ = r^\#, \hs r_- = -(r^{21})^\# = -r_+^*
\end{equation}
(see \eqref{eq-r-sharp}). The Lie bracket on $\g^*$, defined as the dual map of $\delta_r$, is then given by
\begin{equation}\lb{eq-gs-bra}
[\xi, \; \eta] = \ad_{r_+(\xi)}^* \eta - \ad_{r_-(\eta)}^* \xi = 
\ad_{r_-(\xi)}^* \eta - \ad_{r_+(\eta)}^* \xi, \hs \xi, \, \eta \in \g^*,
\end{equation}
where for $x \in \g$ and $\zeta \in \g^*$, the element $\ad_x^* \zeta \in \g^*$ is given by
$\la \ad_x^* \zeta,  y\ra = \la \zeta,  [y, x]\ra$ for $y \in \g$.
It is well-known (\cite[Lecture 4]{Etingof-Schiffmann}, \cite{RT:factorizable}) that both $r_+$ and $r_-$ are Lie algebra homomorphisms.
Set
\begin{equation}\lb{eq-de-l-pm}
\f_+ = {\rm Im} (r_+), \hs \f_- = {\rm Im} (r_-).
\end{equation}
Then both $\f_+$ and $\f_-$ are Lie subalgebras of $\g$, and $\delta_r(\f_\pm) \subset \wedge^2 \f_\pm$. In other words,
both $\f_+$ and $\f_-$ are sub-Lie bialgebras of the Lie bialgebra $(\g, \delta_r)$. A different proof of 
the following \leref{le-m-sub} can be found in \cite[$\S$7.1]{Lu-Mou:mixed}.

\ble{le-m-sub} 
Any Lie subalgebra $\m$ of $\g$ 
containing $\f_+$ or $\f_-$ is a sub-Lie bialgebra of $(\g, \delta_r)$.
\ele

\begin{proof}
Assume that  $\m \supset \f_+$. One needs to show that $\m^0$, the annihilator of $\m$ in 
$\g^*$, is a Lie ideal with respect to the Lie bracket on $\g^*$ given in \eqref{eq-gs-bra}. Let $\xi \in \m^0$ and
$\eta \in \g^*$. Since $\m \supset \f_+$, one has $\m^0 \subset \f_+^0 = \ker r_-$, and it follows from \eqref{eq-gs-bra} that
$[\xi, \, \eta] = -\ad_{r_+(\eta)}^* \xi$. As $r_+(\eta) \in \f_+ \subset \m$, one has $\ad_{r_+(\eta)}^* \xi \in \m^0$.
Thus $\m^0$ is a Lie ideal of $\g^*$. The case when $\m \supset \f_-$ is proved similarly. 
\end{proof}

\subsection{Factorizable quasitriangular $r$-matrices}\lb{subsec-factorizable}
Recall that a quadratic Lie algebra is a pair $(\g, \lara_\g\!)$, where $\g$ is a Lie algebra and $\lara_\g$ is a 
symmetric non-degenerate invariant bilinear form on $\g$.
Given a quadratic Lie algebra $(\g, \lara_\g)$, for any vector subspace $\v$ of $\g$, let
\begin{equation}\lb{eq-v-perp}
\v^\perp = \{x \in \g: \la x, \, \v\ra_\g = 0\}.
\end{equation}
By a {\it Lagrangian subalgebra} of $(\g, \lara_\g)$ we mean a 
Lie subalgebra $\l$ of $\g$ that is also Lagrangian with respect to $\lara_\g$, i.e., $\l^\perp = \l$. By a {\it Lagrangian splitting} 
of $(\g, \lara_\g)$ we mean 
a decomposition $\g = \u + \u'$ where both $\u$ and $\u'$ are Lagrangian subalgebras of $\g$. 
The notion of quadratic Lie algebras with Lagrangian splittings is then equivalent to that of Manin triples
\cite[Lecture 4]{Etingof-Schiffmann}.

Given a Lie bialgebra $(\u, \delta_\u)$, recall that the {\it Drinfeld double Lie algebra} of $(\u, \delta_\u)$ is
the quadratic Lie algebra $(\g, \lara_\g)$, where $\g = \u \oplus \u^*$ as a vector space, $\lara_\g$ is the symmetric bilinear 
form on $\g$ given by 
\[
\la x + \xi, \; y + \eta\ra_\g = \la x, \eta\ra + \la \xi, \, y\ra, \hs x, y \in \u, \; \xi, \eta \in \u^*,
\]
and the Lie bracket on $\g$ is the unique one with respect to which the bilinear form $\lara_\g$ is invariant
and both $\u \cong \u \oplus 0$ and $\u^*\cong 0 \oplus \u^*$ are Lie subalgebras. The decomposition $\g = \u + \u^*$ is thus a Lagrangian splitting of
$(\g, \lara_\g)$. One also refers to $((\g, \lara_\g), \u, \u^*)$ as the {\it Manin triple of the Lie bialgebra
$(\u, \delta_\u)$}.
Conversely,  a Lagrangian splitting
$\g = \u + \u'$ of a quadratic Lie algebra $(\g, \lara_\g)$ 
gives rise to a Lie bialgebra $(\u, \delta_\u)$, where $\delta_\u: \u \to \wedge^2 \u$
is the map dual to the Lie bracket on $\u'$, the latter being identified with $\u^*$
via the pairing between $\u$ and $\u'$ defined by $\lara_\g$, and the Manin triple $((\g, \lara_\g), \u, \u')$
is isomorphic to the Manin triple of $(\u, \delta_\u)$.

Assume now that $r$ is a factorizable quasitriangular $r$-matrix on a Lie algebra $\g$, i.e., the linear map
\[
r_+-r_- = (r+r^{21})^\#: \;\;\; \g^* \lrw \g
\]
is invertible. The symmetric bilinear form
$\lara_{\g}$ on $\g$ given  by 
\begin{equation}\lb{eq-lara-g}
\la x_1, \; x_2 \ra_\g = \la (r_+-r_-)^{-1} x_1, \; x_2\ra = \la x_1, \; (r_+-r_-)^{-1} x_2 \ra, \hs \;\;\; x_1, \, x_2 \in \g.
\end{equation}
is then non-degenerate and invariant, making $(\g, \lara_{\g})$ into a quadratic Lie algebra. We will refer to $\lara_{\g}$ as the {\it symmetric bilinear form on $\g$
associated to $r$}. Set
\begin{equation}\lb{eq-r-flats}
r^\flat_\pm = r_\pm \circ (r_+-r_-)^{-1}: \;\;\; \g \lrw \g.
\end{equation}
One thus has $r^\flat_+-r^\flat_- = {\rm Id}_\g$, and
\[
\la r^\flat_+(x_1), \; x_2\ra_\g + \la x_1, \; r^\flat_-(x_2)\ra_\g = 0, \hs x_1, \, x_2 \in \g.
\]
It also follows from the definitions that $\f_\pm = {\rm Im} (r^\flat_\pm)$ and that
\begin{equation}\lb{eq-ker-rr}
\ker (r_+) = \f_-^0, \hs \ker (r_-) = \f_+^0, \hs \ker (\rfp) = \f_-^\perp, \hs \ker (\rfm) = \f_+^\perp.
\end{equation}
In particular, if $x \in \f_+^\perp$, then $x = \rfp(x)-\rfm(x) = \rfp(x) \in \f_+$, so 
$\f_+^\perp \subset \f_+$. Similarly, $\f_-^\perp \subset \f_-$.

Still assuming that $r$ is factorizable, 
consider now the direct product Lie algebra $\gog$ and the  bilinear form $\lara_{\gog}$ on $\gog$ given by
\begin{equation}\lb{eq-bra-gog}
\la (x_1, \; x_2), \; (x_1^\prime, \; x_2^\prime)\ra_{\gog} = \la x_1,  \; x_1^\prime\ra_\g - 
\la x_2, \; x_2^\prime\ra_\g, \hs x_1, x_2, x_1^\prime, x_2^\prime \in \g.
\end{equation}
One then has the Lagrangian splitting
\begin{equation}\lb{eq-splitting-lr}
\gog = \gdia + \l_r,
\end{equation}
of the quadratic Lie algebra $(\gog, \lara_{\gog})$, 
where 
$\gdia = \{(x, x): x \in \g\}\subset \gog$, and
\begin{equation}\lb{eq-lr}
\l_r = \{(r_+(\xi), \; r_-(\xi)): \; \xi \in \g^*\} =\{(\rfp(x), \; \rfm(x)): \; x \in \g\} \subset \gog.
\end{equation}
One checks that $\delta_r: \g \to \wedge^2 \g$ coincides with the
co-bracket on $\g \cong \gdia$ induced by the Lagrangian splitting in \eqref{eq-splitting-lr}.
The assignment $\gotg \ni r \mapsto \l_r \subset \gog$
gives a one to one correspondence between the set of all  
factorizable quasitriangular $r$-matrices on $\g$  that have $\lara_\g$ as the associated 
symmetric bilinear form and the set of Lagrangian subalgebras $\l$ of $(\gog, \lara_{\gog})$ such that $\gog =\gdia + \l$.

\bex{ex-double-0} Let $(\g, \lara_\g)$ be any quadratic Lie algebra and let $\g = \u + \u'$ be a Lagrangian splitting
of $(\g, \lara_\g)$. 
The quadratic Lie algebra $(\gog, \lara_{\gog})$ (see \eqref{eq-bra-gog})  has the Lagrangian splitting
$\gog = \gdia + \l$, where $\l = \{(\xi, x): \xi \in \u', x \in \u\} \subset \gog$. The factorizable 
quasitriangular $r$-matrix $r_{(\u, \u')}$ on $\g$ such that $\l_{r_{(\u, \u')}} = \l$ is given by
\begin{equation}\lb{eq-r-Manin-0}
r_{(\u, \u')}=\sum_{i=1}^m x_i \otimes \xi_i \in \gotg,
\end{equation}
where $\{x_i\}_{i=1}^m$ is any basis of
$\u$ and  $\{\xi_i\}_{i=1}^m$ the basis of $\u'$ such that $\la x_i, \xi_j\ra_\g = \delta_{ij}$
for $1 \leq i, j \leq m$. We will call $r_{(\u, \u')}$ the {\it $r$-matrix on $\g$ defined by
the Lagrangian splitting $\g = \u + \u'$} of $(\g, \lara_\g)$. It is easy to see that 
the Lie subalgebras $\f_-$
and $\f_+$ of $\g$ associated to $r_{(\u, \u')}$ (see definitions in \eqref{eq-de-l-pm}) are respectively given by
$\f_- = \u$ and $\f_+ = \u'$. In particular, $\f_- \cap \f_+ = 0$. Conversely, let $r$ be a factorizable quasitriangular
$r$-matrix on $\g$ with $\lara_\g$ as the associated symmetric form. If $\f_- \cap \f_+ =0$, then $\g = \f_- + \f_+$ is
a Lagrangian splitting and $r = r_{(\f_-, \f_+)}$.
\hfill$\diamond$
\eex

\bex{ex-BD-0}
Let $\g$ be a complex simple Lie algebra, let $\lara_\g$ be any non-zero scalar multiple of
the Killing form of $\g$, and let 
$\lara_{\gog}$ be the bilinear form on 
$\gog$ given in \eqref{eq-bra-gog}. In \cite{BD1}, Belavin and Drinfeld classified all 
Lagrangian splittings of $(\gog, \lara_{\gog})$ of the form
$\gog = \gdia + \l$ and wrote down explicitly the  corresponding factorizable quasitriangular $r$-matrices on $\g$.
Details on the so-called {\it standard $r$-matrix} $r_{\rm st}$ on $\g$ will be recalled in $\S$\ref{subsec-rst}. 
\hfill$\diamond$
\eex

\subsection{Poisson structures defined by quasitriangular $r$-matrices}\lb{subsec-sigma-r}
Let $Y$ be a manifold with a left Lie algebra action $\lambda: \g \to \V^1(Y)$.
For $r \in \gotg$, let $\lambda(r)$ be the 
tensor field on $Y$ given by
\begin{equation}\lb{eq-pi-Y}
\lambda(r) = \-\sum_i \lambda(x_i) \ot \lambda(y_i), \hs \mbox{if} \;\; r = \sum_i x_i \ot y_i.
\end{equation}
The special case of the following \leref{le-sigma-r-Poisson} when $r$ is the $r$-matrix on a quadratic Lie algebra defined by a Lagrangian
splitting (see \exref{ex-double-0}) is proved in \cite[Theorem 2.3]{Lu-Yakimov:DQ}.

\ble{le-sigma-r-Poisson} \cite{Lu-Mou:mixed} Let $r$ be a quasitriangular $r$-matrix on $\g$.
If $\lambda(r)$ is skew-symmetric, i.e., if $\lambda(r)$ is a bivector field on $Y$, then it is Poisson, and $\lambda$ is a left Poisson action of the
Lie bialgebra $(\g, \delta_r)$ on the Poisson manifold $(Y, -\lambda(r))$.
\ele

Note that $\lambda(r)$ is skew-symmetric if and only if $\lambda(r + r^{21}) =0$,  or, equivalently, 
\begin{equation}\lb{eq-s-qy}
(r_+ -r_-)(\q_y^0) \subset \q_y, \hs y \in Y,
\end{equation}
where for $y \in Y$, $\q_y$ is the stabilizer  of $\lambda$ at $y$, i.e.,
$q_y =\ker (\lambda_y) \subset \g$, 
and $\q_y^0$ the annihilator of $\q_y$ in $\g^*$. In particular, \eqref{eq-s-qy} is independent on the skew-symmetric part of $r$.
Note also that if $r \in \gotg$ is factorizable defining the symmetric bilinear form $\lara_\g$ on $\g$, then 
\eqref{eq-s-qy} is equivalent to $\q_y \subset \g$ being coisotropic with respect to $\lara_\g$ for each $y \in Y$, i.e.,
\begin{equation}\lb{eq-s-qy-3}
q_y^\perp \subset \q_y, \hs \forall \; y \in Y.
\end{equation}

\bde{de-admi} By an {\it admissible quadruple} we mean a quadruple $(\g, r, Y, \lambda)$, where
$\g$ is a Lie algebra, $r \in \gotg$ is a quasitriangular $r$-matrix on $\g$, $Y$ is a manifold, and $\lambda$ is a 
left Lie algebra action of $\g$ on $Y$ such that $\lambda(r)$ is skew-symmetric, i.e., \eqref{eq-s-qy} holds for every $y \in Y$.
Given an admissible quadruple $(\g, r, Y, \lambda)$, we refer to $-\lambda(r)$ (and sometimes $\lam(r)$) as  {\it the Poisson structure
on $Y$ defined by $(r, \lambda)$}. 
\ede

Let $G$ be a connected Lie group with Lie algebra $\g$, and let $r \in \gotg$ be a quasitriangular $r$-matrix on
$\g$. Let $r^L$ (resp. $r^R$) be the
left (resp. right) invariant tensor field on $G$ with value $r$ at the identity element of $G$. Then the bivector field 
on $G$ given by 
\[
\piG = r^L-r^R
\]
is Poisson, making $(G, \piG)$ into a Poisson Lie group \cite{Etingof-Schiffmann}. If $\lambda: G \times Y \to Y$
is a left Lie group action of $G$ on a manifold $Y$ such that $\lambda(r)$ is skew-symmetric,
where $\lambda: \g \to \V^1(Y)$ also denotes the induced left Lie algebra action of $\g$ on $Y$ (see
\eqref{eq-sigma-induced}), we also call $(G, r, Y, \lambda)$ an
admissible triple. In this case,  $\lambda$ is 
a left Poisson action of the Poisson Lie group $(G, \piG)$ on $(Y, -\lambda(r))$.
When the action $\lam$ is transitive, we also say that $(G, r, Y, \lam)$
is a {\it homogeneous admissible quadruple}.

\bre{re-GQ}
It is clear from the definition that a quadruple $(G, r, Y, \lambda)$ is admissible  if and only if $(G, r, \O, \lambda)$
is admissible for every 
$G$-orbit in $\O$ in $Y$. In studying admissible quadruples, we may therefore restrict ourselves to homogeneous 
ones. Let $Q$ be any closed subgroup of $G$ with Lie algebra $\q$, and let 
$\lam_{\sGQ}$ be the left action of $G$ on $G/Q$ by left translation. Then the quadruple
$(G, r, G/Q, \lam_{\sGQ})$ is admissible if and only if $(r_+-r_-)(\q^0) \subset \q$, or equivalently, $\q^\perp \subset \q$ when $r$ is factorizable. 
As a special case, assume that the Lie algebra $\q$ of $Q \subset G$ satisfies $\q \supset \f_+$. Then
\[
(r_+-r_-)(\q^0) \subset (r_+-r_-)(\f_+^0) \subset r_+(\f_+^0) \subset \f_+ \subset \q,
\]
so $(G, r, G/Q, \lam_{\sGQ})$ is admissible. On the other hand, by \leref{le-m-sub},
$Q$ is a Poisson Lie subgroup of $(G, \piG)$, so  $\piG$ projects to a well-defined Poisson structure, denoted by
$\pi_{\scriptscriptstyle{G/Q}}$, on $G/Q$. It is shown in \cite[$\S$7]{Lu-Mou:mixed} that $\pi_{\scriptscriptstyle{G/Q}} = -\lambda_{\sGQ}(r)$. \hfill$\diamond$
\ere

Let $(G, r, Y, \lambda)$ be an
admissible quadruple and consider the Poisson structure $\pi = -\lambda(r)$ on $Y$. 
By the definition of $\pi$ and by \eqref{eq-s-qy}, one has
\begin{equation}\lb{eq-pi-y-sharp}
\pi_y^\#(\alpha_y) = \lambda_y\left(r_+ (\lambda_y^*(\alpha_y))\right) = \lambda_y\left(r_- (\lambda_y^*(\alpha_y))\right),
\hs \alpha_y \in T_y^*Y.
\end{equation}
It follows that for any Lie subalgebra $\m$ of $\g$ containing $\f_+$ or $\f_-$ and for any $y \in Y$, one has
\begin{equation}\lb{eq-m-contain}
{\rm Im} (\pi^\#_y) \subset \lambda_y(\m) \subset T_yY.
\end{equation}
The following \leref{le-M-O} follows immediately from \eqref{eq-m-contain} and \leref{le-sigma-r-Poisson}.

\ble{le-M-O} Let $(G, r, Y, \lambda)$ be an admissible quadruple, and let $M$ be a connected Lie subgroup of $G$ 
such that the Lie algebra $\m$ of $M$ contains $\f_+$ or $\f_-$. Then
every orbit $\O_{\sM}$ of $M$ in $Y$ is a Poisson submanifold of $Y$ with respect to the Poisson structure 
$\pi =-\lambda(r)$, and $\lambda$ restricts to 
a left Poisson action of the Poisson Lie group $(M, \piG|_\sM)$ on $(\O_{\sM}, \pi|_{{\scriptscriptstyle \O_{\sM}}})$.
\ele

In particular,  by taking $M = G$ in \leref{le-M-O}, 
every $G$-orbit in $Y$, equipped with the Poisson structure $\pi$, is a Poisson homogeneous space
of the Poisson Lie group $(G, \piG)$.

Let $\O$ be any $G$-orbit in $Y$. For $y \in \O$, let $({\rm Im}(\pi_y^\#))^0 \subset T_y^*\O$
be the co-normal space in $\O$ at $y$ of the symplectic leaf of $\pi$ through $y$. Let 
$[\lambda_y]:  \g/\q_y \to T_y\O$
 be the vector space isomorphism induced by $\lambda_y: \g \to T_yY$. 
Identify $(\g/\q_y)^*$ with $\q_y^0 \subset \g^*$. Then one has the vector space isomorphism
$[\lambda_y]^*:  T_y^*\O \to \q_y^0$.

\ble{le-corank} For any $y \in \O$, one has
\[
[\lambda_y]^*(({\rm Im}(\pi_y^\#))^0) = \q_y^0 \cap r_+^{-1}(\q_y) = \q_y^0 \cap r_-^{-1}(\q_y).
\]
Consequently, $\dim \O - \dim ({\rm Im}(\pi_y^\#)) = \dim (\q_y^0 \cap r_+^{-1}(\q_y)) = \dim (\q_y^0 \cap r_-^{-1}(\q_y))$.
\ele

\begin{proof}
By \eqref{eq-pi-y-sharp}, one has
${\rm Im}(\pi_y^\#) = \lambda_y(r_+(\q_y^0)) = [\lambda_y]\left((\q_y + r_+(\q_y^0))/\q_y\right)$. 
Thus
\[
[\lambda_y]^*(({\rm Im}(\pi_y^\#))^0) 
=\q_y^0 \cap \left(r_+(\q_y^0)\right)^0 = \q_y^0 \cap r_-^{-1}(\q_y)
=\q_y^0 \cap r_+^{-1}(\q_y).
\]
\end{proof}

We recall a result of Drinfeld \cite{dr:homog} on Poisson homogeneous spaces: suppose that
$(\O, \pi, \lambda)$ is a Poisson homogeneous space of a Poisson Lie group with Lie bialgebra $(\g, \delta_\g)$.
Then the Lie algebroid $(\O \times \g)\bowtie T^*\O$ over $\O$ in \eqref{eq-A-sigma} is transitive, so
the kernel of its anchor map $-\lambda + \pi^\#$ is a bundle
of Lie algebras over $\O$. Consider the map 
\[
\Psi: \;\; (\O \times \g)\bowtie T^*\O \lrw \d, \;\;\; (x, \alpha_y) \longmapsto x + \lambda_y^*(\alpha_y), \hs 
x \in \g, \, y \in \O, \, \alpha_y \in T_y^*\O,
\]
where $(\d, \lara_\d)$ is the Drinfeld double Lie algebra of $(\g, \delta_\g)$. 
For $y \in \O$, let 
\[
\l_y = \Psi(\ker (-\lambda_y+\pi_y^\#)) =\{x + \xi: \; x \in \g, \; \xi \in \q_y^0, \;
\lambda_y(x) = \pi_y^\#(([\lambda_y]^*)^{-1}(\xi))\} \subset \d.
\]
By \cite{Lu:Duke-Poi}, for  $y \in Y$, $\Psi|_{\ker (-\lambda_y+\pi_y^\#)}: \ker (-\lambda_y+\pi_y^\#) \to \d$
is injective and $\l_y$ is a Lagrangian subalgebra of $(\d, \lara_\d)$, called the {\it Drinfeld Lagrangian
subalgebra at $y$} associated to the Poisson structure $\pi$. 

The following \leref{le-ly} follows directly from the definition of $\l_y$ and is basic to Drinfeld's theory on
Poisson homogeneous spaces \cite{dr:homog}.

\ble{le-ly}
(i). The conormal subspace $({\rm Im}(\pi_y^\#))^0 \subset T_y^*\O$ can be identified with $\g^* \cap \l_y$ under
the vector space 
isomorphism $[\lambda_y]^*: T_y^*\O \to \q_y^0$;

(ii). For $x \in \g$, $\lambda_y(x) \in {\rm Im}(\pi_y^\#)$ if and only if $x \in \g \cap (\g^* + \l_y) = {\rm pr}_\g(\l_y)$,
where ${\rm pr}_\g: \d \to \g$ is the projection with respect to the decomposition $\d = \g + \g^*$.
\ele

In the context of \leref{le-ly}, one checks from the definition that 
\begin{equation}\lb{eq-de-ly}
\l_y = \{x + \xi \in \d: \; x \in \g, \; \xi \in \q_y^0, \; x+r_+(\xi) \in \q_y\}.
\end{equation}
Thus $\g^* \cap \l_y = \q_y^0 \cap r_+^{-1}(\q_y)$. This gives another proof of \leref{le-corank}
using Drinfeld's general theory 
in \cite{dr:homog}.
We will return to the Drinfeld Lagrangian subalgebras in $\S$\ref{subsec-strongly-admi}.

\subsection{Strongly admissible quadruples}\lb{subsec-strongly-admi} 

\bde{de-admi-strong}
By a {\it strongly admissible quadruple} we mean a quadruple $(G, r, Y, \lambda)$, where
$G$ is a connected Lie group, $r$ is a factorizable quasitriangular $r$-matrix on the Lie algebra 
$\g$ of $G$,  $Y$ is a manifold, and $\lambda$ is a 
left action of $G$ on $Y$, such that 
the stabilizer subgroup $Q_y$ of $G$
at every $y \in Y$ is connected and its Lie algebra $\q_y$ satisfies 
\begin{equation}\lb{eq-qy}
[\q_y, \; \q_y] \subset \q_y^\perp \subset \q.
\end{equation}
\ede

\bre{re-admi-2}
A strongly admissible quadruple is thus an admissible quadruple $(G, r, Y, \lam)$ with the additional requirements
that $r$ be factorizable, the stabilizer subgroup $Q_y$ of $G$ at each $y \in Y$ be connected and its Lie algebra satisfy 
$[\q_y, \q_y] \subset \q_y^\perp$.
\hfill$\diamond$
\ere

\bex{ex-homog}
Let $G$ be a connected Lie group with Lie algebra $\g$ and $r$ a factorizable quasitriangular $r$-matrix on $\g$.
{\it Homogeneous} strongly admissible quadruples are then of the form $(G, r, G/Q, \lam)$,   
where $Q$ is a closed and connected Lie subgroup of $G$ whose Lie algebra $\q$ 
satisfies 
\begin{equation}\lb{eq-on-qq}
[\q, \q] \subset \q^\perp \subset \q,
\end{equation}
and $\lam_{\sGQ}$ is the left action of $G$ on $G/Q$ by left translation.
In the special case when $r
=r_{(\u, \u')}$ is the $r$-matrix on $\g$ defined by a Lagrangian splitting
$\g = \u + \u'$ of $(\g, \lara_\g)$, Condition \eqref{eq-on-qq}  was first introduced in \cite{Lu-Yakimov:DQ} and some properties of the Poisson 
structure $\lam_{\sGQ}(r)$ on $G/Q$ were also studied in \cite{Lu-Yakimov:DQ}.
Note also that \eqref{eq-on-qq} holds automatically
if $\q$ is Lagrangian with respect to $\lara_\g$.
\hfill$\diamond$
\eex

\bex{ex-BD-1}
Continuing with \exref{ex-BD-0}, 
let $G$ be a connected complex simple Lie group with Lie algebra $\g$, and let $\lara_\g$ be
a fixed non-zero scalar multiple of the Killing form of $\g$. Recall from \exref{ex-BD-0} that 
Lagrangian splittings of $(\gog, \lara_{\gog})$ of the form $\gog = \gdia + \l$ have been classified by 
Belavin-Drinfeld in \cite{BD1}. Let $\rBD \in \gotg$ be the factorizable quasitriangular 
$r$-matrix on $\g$ corresponding to such an $\l \subset \gog$. 
Let $P$ be any parabolic subgroup of $G$. As
the Lie algebra $\p$ of $P$ is coisotropic with respect to $\lara_\g$, the quadruple 
$(G, \rBD, G/P, \lambda_{\sGP})$ is admissible. On the other hand, the Lie algebra $\p$ of $P$ satisfies 
$[\p, \p] \subset \p^\perp$ if and only if $\p$ is Borel. Thus 
$(G, \rBD, G/P, \lambda_{\sGP})$ is strongly admissible if and only if $P$ is a Borel subgroup of $G$.
\hfill$\diamond$ \eex

We now prove some preliminary properties of strongly admissible quadruples.

\ble{le-O-yy}
Assume that  $(G, r, Y, \lambda)$ is strongly admissible. Then for any $y_1, y_2 \in Y$ in the same $G$-orbit in $Y$, the linear isomorphism
\begin{equation}\lb{eq-phi-yy}
I_{y_2, y_1}: \;\; \q_{y_1}/\q_{y_1}^\perp \lrw \q_{y_2}/\q_{y_2}^\perp, \;\;\; x + \q_{y_1}^\perp \longmapsto 
\Ad_g(x) + \q_{y_2}^\perp, \hs x \in \q_{y_1},
\end{equation}
is independent of the choice of $g \in G$ such that $gy_1 = y_2$.
\ele

\begin{proof}
For any  $y \in Y$, as 
the stabilizer subgroup $Q_{y}$ of $G$ at $y \in Y$ is connected, condition \eqref{eq-qy} implies that
the action of $Q_{y}$ on $\q_{y}/\q_{y}^\perp$ induced by the adjoint action of $Q_{y}$ on $\q_{y}$ is trivial. 
Consequently, the map $I_{y_2, y_1}$ is independent of the choice of $g \in G$ such that $gy_1 = y_2$.
\end{proof}
 
It is also clear that if $y_1, y_2, y_3$ are in the same $G$-orbit, then
\begin{align}\lb{eq-phi-yyy}
I_{y_2, y_1}^{-1} &= I_{y_1, y_2}: \;\;\; \q_{y_2}/\q_{y_2}^\perp \lrw \q_{y_1}/\q_{y_1}^\perp, \\
\lb{eq-phi-yyy-1}
I_{y_3, y_2} \circ I_{y_2, y_1} &= I_{y_3, y_1}: \;\;\; \q_{y_1}/\q_{y_1}^\perp \lrw \q_{y_3}/\q_{y_3}^\perp.
\end{align}

\bde{de-l-yy}
Let $(G, r, Y, \lambda)$ be a strongly admissible quadruple.  For $y_1, y_2 \in Y$ in the same $G$-orbit, define
\begin{equation}\lb{eq-l-yy}
\l_{y_1, y_2} = \{(x_1, x_2) \in \q_{y_1} \oplus \q_{y_2}: \; I_{y_2, y_1}(x_1 + \q_{y_1}^\perp) = x_2 + \q_{y_2}^\perp\} \subset \gog.
\end{equation}
\ede

Recall from $\S$\ref{subsec-factorizable} that, as $r$ is factorizable,  the Drinfeld double Lie algebra of the Lie bialgebra $(\g, \delta_r)$ can be identified with the
quadratic Lie algebra $(\gog, \lara_{\gog})$, where $\lara_{\gog}$ is defined in \eqref{eq-bra-gog}.

\ble{le-l-yy}
Let $(G, r, Y, \lambda)$ be a strongly admissible quadruple.  For any $y_1, y_2 \in Y$ in the same $G$-orbit,
$\l_{y_1, y_2}$ is a Lagrangian subalgebra of $(\gog, \lara_{\gog})$. Moreover, for any $y \in Y$, 
\begin{equation}\lb{eq-l-yy-2}
\l_{y,y} = \{(x_1, x_2) \in \q_y \oplus \q_y: \; x_1 - x_2 \in \q_y^\perp\}
\end{equation}
is the Drinfeld Lagrangian subalgebra $\l_y$ at $y$ associated to the Poisson structure $\pi = -\lambda(r)$ on 
the $G$-orbit $Gy \subset Y$.
\ele

\begin{proof}
Let $y \in Y$ be arbitrary and let $\O = Gy \subset Y$. By \eqref{eq-de-ly}, the Drinfeld 
Lagrangian subalgebra $\l_y$ at $y$ associated to the Poisson structure $\pi = -\lambda(r)$ on $\O$ is given by
\[
\l_y = \{(x, \,x) + (\rfp(x'), \,\rfm(x')): \; x \in \g, \; x' \in \q_y^\perp, \, x + \rfp(x') \in \q_y\},
\]
from which it follows that $\l_y = \l_{y, y}$. For any  $y_1, y_2 \in Gy$, let  
$y_1 = g_1 y$ and $y_2 = g_2 y$, where $g_1, g_2 \in G$. It follows from the definitions that
\begin{equation}\lb{eq-l-yy-yy}
\l_{y_1, y_2} = (\Ad_{g_1} \oplus \Ad_{g_2}) (\l_{y, y}).
\end{equation}
Thus $\l_{y_1, y_2}$ is a Lagrangian subalgebra of  
$(\gog, \lara_{\gog})$.
\end{proof}

\bre{re-ly-2}
By \leref{le-ly} and \leref{le-l-yy}, for a strongly admissible quadruple $(G, r, Y, \lambda)$,
 the rank of $\pi$ at $y \in Y$ is equal to 
$\dim (Gy) - \dim (\l_r \cap \l_{y, y})$. \hfill$\diamond$
\ere

\sectionnew{Orbit intersections for strongly admissible quadruples}\lb{sec-orbits}

\subsection{The set-up}\lb{subsec-setup}
Assume that $(G, r, Y, \lambda)$ is an admissible quadruple, in which $r$ is factorizable.
Let $M_+$ and $M_-$ be connected Lie subgroups of $G$ whose respective Lie algebras $\m_+$ and $\m_-$ satisfy
\begin{equation}\lb{eq-fmm}
\f_+ \subset \m_+, \hs \f_- \subset \m_-.
\end{equation}
 Let $\T$ be the connected component of $M_+ \cap M_-$
containing the identity element. The Lie algebra $\t$ of $\T$ is then given by
$\t = \m_+ \cap \m_-$.
By \leref{le-m-sub}, $\T$ is a Poisson Lie subgroup of $(G, \piG)$, where $\piG = r^L - r^R$, and
$(\t, \delta_r|_{\t})$ is a sub-Lie bialgebra of the Lie bialgebra $(\g, \delta_r)$.
Consider the decomposition
\[
Y = \bigsqcup_{\O_+, \O_-} \OO,
\]
where $\O_+$ and $\O_-$ are respectively $M_+$-orbits and $M_-$-orbits in $Y$.
As $r$ is factorizable, one has $\f_+ + \f_- =\g$, and thus
$\m_+ + \m_- = \g$. Consequently, each non-empty intersection $\O_+ \cap \O_-$ is a smooth submanifold of $Y$, and by \leref{le-M-O}, also
a Poisson submanifold with respect to the Poisson structure $\pi = -\lambda(r)$ on $Y$. Denote the restriction of $\pi$ to 
$\O_+ \cap \O_-$ also by $\pi$. As $\O_+ \cap \O_-$ is $\T$-invariant, the Poisson action $\lambda$ of $(G, \piG)$ on $(Y, \pi)$
restricts to a Poisson action on $(\O_+ \cap \O_-, \pi)$ by the Poisson Lie group $(\T, \piG|_{{\scriptscriptstyle{{\mathbb T}}}})$ and the Lie bialgebra $(\t, \delta_r|_\t)$.

\bde{de-admi-6}
The six-tuple $(G, r, Y, \lambda, M_+, M_-)$ is said to be {\it admissible} if the quadruple ($G, r, Y, \lambda)$ is
strongly admissible and if $\T$ is abelian, acting on $(Y, \lambda(r))$ by Poisson isomorphisms, and
the $\T$-leaves of $\lambda(r)$ in $Y$ are precisely all the connected components of non-empty intersections
$\OO$ where $\O_+$ is an $M_+$-orbit and $\O_-$ an $M_-$-orbit in $Y$.
\ede

\bre{re-MM-connected}
If $G$ is an affine algebraic group over $\Cset$ and if $M_+$ and  $M_-$  are algebraic subgroups
of $G$ such that $M_+ \cap M_-$ is connected, by \cite[Corollary 1.5]{R}, every non-empty intersection of an $M_+$-orbit and an $M_-$-orbit in $Y$ is irreducible and thus connected. 
\hfill$\diamond$
\ere

In  $\S$\ref{sec-orbits}, we develop a test for a six-tuple $(G, r, Y, \lambda, M_+, M_-)$
to be  admissible and we compute the $\T$-leaf stabilizers.
The test will be used in $\S$\ref{sec-mixed} to identify a class of admissible
six-tuples.

More precisely, for a strongly admissible quadruple $(G, r, Y, \lambda)$, and a pair $(M_+, M_-)$ of
connected Lie subgroups of $G$ satisfying \eqref{eq-fmm}, we show that $\S$\ref{subsec-k-regular} that the
Lie bialgebra action $\lam|_\t$ of $(\t, \delta_r|_\t)$ on any non-empty intersection
$\OO$ of $M_+$- and $M_-$-orbits is regular, and we give a formula for the integer $\delOO$
that measures how far the action 
$\lambda|_\t$ of $(\t, \delta_r|_\t)$ on $\OO$ is from being full
(see \deref{de-full}).
When $\delta_r|_\t = 0$, so $\T$ acts on $(Y, \pi)$ by Poisson isomorphisms, we show that
the Poisson structure $\pi$ is also regular on each 
non-empty intersection $\OO$. In $\S$\ref{subsec-leaf-stabi},
under a further assumption on the
pair $(M_+, M_-)$, called $r$-admissible, which implies that $\T$ is abelian and $\delta_r|_\t = 0$, we compute the leaf
stabilizer for each $\T$-leaf in $Y$. \thref{th-main}, a summary of the main results in $\S$\ref{sec-orbits}, is proved in
$\S$\ref{subsec-proof-th-main}.

\bre{re-k-lr}
The assumptions that $\f_+ \subset \m_+$ and $\f_- \subset \m_-$ are equivalent to
\begin{equation}\lb{eq-k-lr}
\t_{\rm diag} + \l_r = \m_+ \oplus \m_-,
\end{equation}
where $\l_r$ is the Lagrangian subalgebra of $\gog$ given in \eqref{eq-lr}, and $\t_{\rm diag} = \{(x, x): x \in \t\}$.
Indeed, clearly $\f_+ \subset \m_+$ and $\f_- \subset \m_-$ if and only if 
$\l_r \subset \m_+ \oplus \m_-$, which, by the direct sum decomposition $\gog = \g_{\rm diag} + \l_r$, is 
 equivalent to $\m_+\oplus \m_- = \gdia \cap (\m_+ \oplus \m_-) + \l_r = \t_{\rm diag} + \l_r$.
\hfill$\diamond$
\ere

%\subsection{Statements of main results of $\S$\ref{sec-orbits}}\lb{subsec-statements-1}
%The following \prref{pr-regular-sigma} will be proved in $\S$\ref{subsec-full}.

\subsection{The action of $(\t, \delta_r|_{\t})$ on $(\OO, \pi)$ is regular}\lb{subsec-k-regular}
Assume first that $(G, r, Y, \lambda)$ is admissible, in which $r$ is factorizable, and let 
$(M_+, M_-)$ be a pair of connected Lie subgroups of $G$ satisfying \eqref{eq-fmm}.
Let $\pi = -\lambda(r)$ on $Y$. For $y \in Y$, define 
\begin{equation}\lb{eq-delta-y-00}
\delta_y = \dim ((M_+y) \cap (M_-y)) - \dim \,( \lambda_y(\t) +{\rm Im}(\pi_y^\#)).
\end{equation}
Then for any pair $(\O_+, \O_-)$ of $M_+$ and $M_-$-orbits in $Y$ such that $\OO \neq \emptyset$, the Poisson action $\lambda|_\t$ of
the Lie bialgebra $(\t, \delta_r|_\t)$ on $(\OO, \pi)$ is regular (see \deref{de-full}) if and only of the integer-valued 
function $y \mapsto \delta_y$ on $\OO$ is a constant function.
For $y \in Y$, let $p_y: \q_y \to \q_y/\q_y^\perp$ be the natural projection, and define 
two subspaces $\a_y^\pm$ of $\q_y/\q_y^\perp$  by
\begin{equation}\lb{eq-ay-pm}
\a_y^+ = p_y\left(\m_+^\perp \cap \q_y\right) = \frac{\m_+^\perp \cap \q_y}{\m_+^\perp \cap \q_y^\perp},\hs\hs
\a_y^- = p_y\left(\m_-^\perp \cap \q_y\right) = \frac{\m_-^\perp \cap \q_y}{\m_-^\perp \cap \q_y^\perp}.
\end{equation}

\ble{le-delta-y-1} 
One has $\delta_y = \dim (\a_y^+ \cap \a_y^-)$ for every $y \in Y$.
\ele
 
\begin{proof}
Let $\O = G y\subset Y$ and 
and let $[\lambda_y]: \g/\q_y \cong T_y\O$ be the vector space isomorphism induced by the linear map 
$\lambda_y: \g \to T_yY$. Let $\O_+ = M_+ y \subset \O$ and $\O_- = M_-y \subset \O$. Then
\begin{align*}
T_y(\OO)& = \lambda_y(\m_+) \cap \lambda_y(\m_-) = \lambda_y(\m_+ \cap (\m_- + \q_y)) = 
[\lambda_y]\left((\q_y+\m_+ \cap (\m_- + \q_y))/\q_y\right),\\
\lambda_y(\t) +{\rm Im}(\pi_y^\#) & = \lambda_y(\t + r_+(\q_y^0)) = [\lambda_y]\left(
(\q_y + \m_+ \cap \m_- + r_+(\q_y^0))/\q_y\right).
\end{align*}
On the other hand, it follows from $(r_+-r_-)(\q_y^0) = \q_y^\perp$ that 
\begin{equation}\lb{eq-mm-perp}
\m_+ \cap \m_- + r_+(\q_y^0) = \m_+ \cap (\m_- + \q_y^\perp).
\end{equation}
Indeed, as $\f_+ \subset \m_+$, one has $\m_+ \cap \m_- + r_+(\q_y^0)\subset \m_+$, and as 
$(r_+-r_-)(\q_y^0) = \q_y^\perp$, one has 
\[
\m_+ \cap \m_- + r_+(\q_y^0) \subset \m_- + (r_+-r_-)(\q_y^0)+ r_-(\q_y^0) \subset \m_- + \q_y^\perp,
\]
so $\m_+ \cap \m_- + r_+(\q_y^0) \subset \m_+ \cap (\m_- + \q_y^\perp)$. Conversely, let $x \in 
\m_+ \cap (\m_- + \q_y^\perp)$ and write $x = x_- + (r_+-r_-)(\xi)$ for $x_- \in \m_-$ and $\xi \in \q_y^0$. Then 
$x_- - r_-(\xi) = x -r_+(\xi) \in \m_+ \cap \m_-$ and $x =  x_- - r_-(\xi) +r_+(\xi) \in 
\m_+ \cap \m_- + r_+(\q_y^0)$. Thus $\m_+ \cap (\m_- + \q_y^\perp) \subset \m_+ \cap \m_- + r_+(\q_y^0)$, and this 
proves \eqref{eq-mm-perp}. Consequently, under the vector space isomorphism $[\lambda_y]: \g/\q_y \cong T_y\O$, one has
\[
T_y(\OO) \cong (\q_y+\m_+ \cap (\m_- + \q_y))/\q_y, \hs
\lambda_y(\t) +{\rm Im}(\pi_y^\#) \cong (\q_y+\m_+ \cap (\m_- + \q_y^\perp))/\q_y.
\]
Using $\q_y^\perp \subset \q_y$, one has
\begin{align*}
(\q_y+\m_+ \cap (\m_- + \q_y^\perp))^\perp & = \q_y^\perp \cap (\m_+^\perp + \m_-^\perp \cap \q_y) =
\q_y^\perp \cap (\m_+^\perp \cap \q_y + \m_-^\perp \cap \q_y),\\
(\q_y+\m_+ \cap (\m_- + \q_y))^\perp & = \q_y^\perp \cap (\m_+^\perp + \m_-^\perp \cap \q_y^\perp) =
\q_y^\perp \cap (\m_+^\perp \cap \q_y^\perp + \m_-^\perp \cap \q_y^\perp).
\end{align*}
It follows that
\[
\delta_y = 
\dim \left(\frac{\q_y^\perp \cap (\m_+^\perp \cap \q_y + \m_-^\perp \cap \q_y)}
{\q_y^\perp \cap (\m_+^\perp \cap \q_y^\perp + \m_-^\perp \cap \q_y^\perp)}\right).
\]
Note that $\m_+^\perp \cap \m_-^\perp = (\m_+ + \m_-)^0 = \g^0 =0$. Writing an element 
$x \in \q_y^\perp \cap (\m_+^\perp \cap \q_y + \m_-^\perp \cap \q_y)$ uniquely as $x = x_+ + x_-$, 
where $x_+ \in \m_+^\perp \cap \q_y$ and $x_-\in \m_-^\perp \cap \q_y$, the map
\[
\q_y^\perp \cap (\m_+^\perp \cap \q_y + \m_-^\perp \cap \q_y) \lrw \q_y/\q_y^\perp, \;\;\; x \longmapsto 
x_+ + \q_y^\perp = -x_- + \q_y^\perp
\]
induces a well-defined vector space isomorphism
\[
\frac{\q_y^\perp \cap (\m_+^\perp \cap \q_y + \m_-^\perp \cap \q_y)}
{\q_y^\perp \cap (\m_+^\perp \cap \q_y^\perp + \m_-^\perp \cap \q_y^\perp)} \lrw
\a_y^+ \cap \a_y^-.
\]
It follows that $\delta_y = \dim (\a_y^+ \cap \a_y^-)$.
\end{proof}

For the remainder of $\S$\ref{subsec-k-regular}, assume that $(G, r, Y, \lambda)$ is strongly admissible, and let
$(M_+, M_-)$ be a pair of connected Lie subgroups of $G$ satisfying \eqref{eq-fmm}.

\ble{le-a-pm} For any $M_+$-orbit $\O_+$ and any $M_-$-orbit $\O_-$ in $Y$, one has

1) $I_{y_2, y_1}(\a_{y_1}^+) = \a_{y_2}^+$ for all $y_1, y_2 \in \O_+$;

2) $I_{y_2, y_1}(\a_{y_1}^-) = \a_{y_2}^-$ for all $y_1, y_2 \in \O_-$.
\ele

\begin{proof}
We only prove 1), the proof of 2) being similar. Assume thus $y_2 = m_+ y_1$, where $m_+ \in M_+$. As
$\Ad_{m_+} \m_+^\perp = \m_+^\perp$, one has 
\[
I_{y_2, y_1}(\a_{y_1}^+) = \frac{\Ad_{m_+}(\m_+^\perp \cap \q_{y_1})}{\Ad_{m_+}(\m_+^\perp \cap \q_{y_1}^\perp)}
= \frac{\m_+^\perp \cap \q_{y_2}}{\m_+^\perp \cap \q_{y_2}^\perp} = \a_{y_2}^+.
\]
\end{proof}
 
Let $(\O_+, \O_-)$ be any pair of $M_+$- and $M_-$-orbits contained in the same $G$-orbit $\O$ in $Y$, possibly $\OO 
=\emptyset$. 
Let $y_0 \in \O$, $y_+ \in \O_+$ and $y_- \in \O_-$ be arbitrary, and let $g_+, g_- \in G$ be such that
$y_+ = g_+ y_0$ and $y_- = g_- y_0$. Define 
\begin{align}\lb{eq-delOO-2}
\delOO &= \dim \left((I_{y_0, y_+} (\a_{y_+}^+)) \cap (I_{y_0, y_-}(\a_{y_-}^{-}))\right) \\
\nonumber
&=\dim \left( \left(\frac{\left(\Ad_{g_+^{-1}} \m_+^\perp\right) \cap \q_{y_0}}{\left(\Ad_{g_+^{-1}} \m_+^\perp\right) \cap \q_{y_0}^\perp}
\right)
\; \cap \; \left(\frac{\left(\Ad_{g_-^{-1}} \m_-^\perp\right) \cap \q_{y_0}}{\left(\Ad_{g_-^{-1}} \m_-^\perp\right) \cap \q_{y_0}^\perp}\right)\right),
\end{align}
where the intersection on the right hand side of \eqref{eq-delOO-2} is in the vector space $\q_{y_0}/\q_{y_0}^\perp$.

\bpr{pr-delOO} For any pair 
$(\O_+, \O_-)$ of $M_+$- and $M_-$-orbits contained in the same $G$-orbit in $Y$, the integer $\delOO$ in 
\eqref{eq-delOO-2} is independent of the choices of  $y_0 \in \O$, $y_+ \in \O_+$ and $y_- \in \O_-$. 
Moreover, when $\OO \neq \emptyset$, one has $\delOO = \delta_y$ for any $y \in \OO$.
\epr

\begin{proof}
For $y_0 \in \O$, $y_+ \in \O_+$, and $y_- \in \O_-$, let
\begin{equation}\lb{eq-a-O-pm}
\a_{\sO_+, \,y_0} = I_{y_0, y_+}(\a_{y_+}^+) \subset \q_{y_0}/\q_{y_0}^\perp, \hs
\a_{\sO_-, \,y_0} = I_{y_0, y_-}(\a_{y_-}^-) \subset \q_{y_0}/\q_{y_0}^\perp.
\end{equation}
By \leref{le-O-yy} and \eqref{eq-phi-yyy-1},
$\a_{\sO_+, \,y_0}$ and $\a_{\sO_-, \, y_0}$ are independent of the choices of $y_+ \in \O_+$ and $y_- \in \O_-$, and the integer
$\delOO = \dim(\a_{\sO_+, \, y_0} \cap \a_{\sO_-, \, y_0})$ is independent of the choice of $y_0 \in \O$. 

Assume that $\OO \neq \emptyset$, and let $y \in \OO$. Then 
\begin{align*}
I_{y_0, y} (\a_y^+) & = I_{y_0, y_+}( I_{y_+, y}(\a_y^+)) = I_{y_0, y_+}(\a_{y_+}^+) = \a_{\sO_+, y_0}, \\
I_{y_0, y} (\a_y^-) & = I_{y_0, y_-}( I_{y_-, y}(\a_y^-)) = I_{y_0, y_-}(\a_{y_-}^-) = \a_{\sO_-, y_0}.
\end{align*}
It follows that $\delOO = \dim(\a_{\sO_+, \, y_0} \cap \a_{\sO_-, \, y_0}) = \dim (\a_y^+ \cap \a_y^-) = \delta_y.$
\end{proof}

\bco{co-regular-sigma}
Assume that $(G, r, Y, \lambda)$ is strongly admissible
and let $(M_+, M_-)$ be a pair of connected Lie subalgebras of $G$ satisfying \eqref{eq-fmm}. Then for every  pair $(\O_+, \O_-)$ of $M_+$- and $M_-$-orbits in
$Y$ such that $\OO \neq \emptyset$, the Poisson action $\lambda|_\t$ of the Lie bialgebra $(\t, \delta_r|_\t)$ on
$(\OO, -\lam(r))$ is regular. 
\eco

\bre{re-algo} For a pair $(\O_+, \O_-)$ of $M_+$ and $M_-$-orbits in $Y$ such that $\OO \neq \emptyset$,
the Poisson action $\lambda|_\t$ of the Lie bialgebra $(\t, \delta_r|_\t)$ on $(\OO, -\lam(r))$ is full if and only if 
$\delOO = 0$. Thus,
the integer $\delOO$ measures how far
the Poisson action $\lambda|_\t$ of $(\t, \delta_r|_\t)$ on $(\OO, -\lam(r))$ is from being full.  
In examples, one can compute $\delOO$ using \eqref{eq-delOO-2} by choosing  $y_0 \in \O$ and $y_\pm \in \O_\pm$ that are convenient for
the computation.  
This is the case for the Poisson structures to be considered in $\S$\ref{sec-mixed}.
\hfill$\diamond$
\ere

\bco{co-lag}
Assume that $(G, r, Y, \lambda)$ is strongly admissible
and let $(M_+, M_-)$ be a pair of connected Lie subalgebras of $G$ satisfying \eqref{eq-fmm}. If $\q_y \subset \g$ is Lagrangian for every $y \in Y$, then
for every  pair $(\O_+, \O_-)$ of $M_+$- and $M_-$-orbits in
$Y$ such that $\OO \neq \emptyset$, the Poisson action $\lambda|_\t$ of the Lie bialgebra $(\t, \delta_r|_\t)$ on
$(\OO, -\lam(r))$ is full.
\eco

\begin{proof}
As $\q_{y_0} = \q_{y_0}^\perp$ for any $y_0 \in Y$, it follows trivially from \prref{pr-delOO} that $\delOO = 0$.
\end{proof}

Still assuming that $(G, r, Y, \lambda)$ is strongly admissible and that $(M_+, M_-)$ satisfies \eqref{eq-fmm}, we now give another formula for the integer
$\delOO$. Let $(\O_+, \O_-)$ be any pair of $M_+$- and $M_-$-orbits in the same $G$-orbit in 
$Y$, and let $y_+ \in \O_+$ and $y_- \in \O_-$. Recall the Lagrangian Lie subalgebra
$\l_{y_+, y_-}$ of $(\gog, \lara_{\gog})$ given by 
\[
\l_{y_+, y_-} = \{(x_+, \, x_-): \; x_+ \in \q_{y_+}, \, x_- \in \q_{y_-}, \, I_{y_-, y_+} 
(x_+ + \q_{y_+}^\perp) = x_- + \q_{y_-}^\perp\}.
\]
Let $p_+: \l_{y_+, y_-} \to \q_{y_+}/\q_{y_+}^\perp$ be the given by
\[
p_+(x_+, x_-)= x_+ + \q_{y_+}^\perp = I_{y_+, y_-}(x_- + \q_{y_-}^\perp), \hs (x_+, x_-) \in \l_{y_+, y_-}.
\]

\ble{le-delOO-3}
Let
$(\O_+, \O_-)$ be any pair of $M_+$- and $M_-$-orbits in
$Y$ contained in the same $G$-orbit, and let $y_+ \in \O_+$ an $y_- \in \O_-$ be arbitrary. Then
\begin{equation}\lb{eq-delOO-3}
\delOO = \dim \left(p_+ \left((\m_+^\perp \oplus \m_-^\perp)\cap \l_{y_+, y_-}\right)\right).
\end{equation}
\ele

\begin{proof}
Let $y_0$ be any point in the unique $G$-orbit $\O$ containing both $\O_+$ and $\O_-$. Then the vector space isomorphism
$I_{y_0, y_+}: \q_{y_+}/\q_{y_+}^\perp \to \q_{y_0}/\q_{y_0}^\perp$ induces an isomorphism 
\[
p_+ \left((\m_+^\perp \oplus \m_-^\perp)\cap \l_{y_+, y_-}\right) \lrw 
I_{y_0, y_+}(\a_{y_+}^+) \cap I_{y_0, y_-}(\a_{y_-}^-).
\]
Consequently, $\delOO = \dim (I_{y_0, y_+}(\a_{y_+}^+) \cap I_{y_0, y_-}(\a_{y_-}^-))$ is also given by \eqref{eq-delOO-3}.
\end{proof}

We now turn to the ranks of the Poisson structure $\pi = -\lambda(r)$ in $Y$. 
Recall that 
$\pi$ is said to be regular on a Poisson submanifold $X$ of $(Y, \pi)$ if 
it has constant rank on $X$.

Recall our assumptions that $\f_+ \subset \m_+$ and $\f_- \subset \m_-$ in \eqref{eq-fmm}, which are equivalent to 
$\l_r \subset \m_+ \oplus \m_-$, which is in turn equivalent to 
$\t_{\rm diag} + \l_r  = \m_+ \oplus \m_-$, where $\t = \m_+ \cap \m_-$ and 
$\l_r$ is the Lagrangian subalgebra of $\gog$ given by \eqref{eq-lr} (see \reref{re-k-lr}).
Let $\n_{\gog}(\l_r)\subset \gog$ be the normalizer of $\l_r$ in $\gog$. Recall the co-bracket
$\delta_r: \g \to \wedge^2 \g$
on $\g$ defined in \eqref{eq-delta-r}. 
%The following \prref{pr-regular-pi} will be proved in $\S$\ref{subsec-regular}.

\ble{le-on-m-equi}
Under the assumption that  $\l_r  \subset \m_+ \oplus \m_-$, the following are equivalent:
\begin{align}\lb{eq-on-m-00}
&\m_+ \oplus \m_- \subset \n_{\gog} (\l_r),\\
\lb{eq-on-m-01}
&[\m_+, \; \f_+]\subset \f_+^\perp, \hs\hs [\m_-, \; \f_-] \subset \f_-^\perp,\\
\lb{eq-on-m-02}
&\delta_r(x) = 0, \hs \forall \; x \in \t.
\end{align}
\ele

\begin{proof}
It is clear that \eqref{eq-on-m-00} is equivalent to 
\[
\m_+ \oplus 0 \subset \n_{\gog} (\l_r) \hs \mbox{and} \hs 0 \oplus \m_- \subset \n_{\gog} (\l_r).
\]
By definition, $\m_+ \oplus 0\subset \n_{\gog} (\l_r)$ if and only if 
for any $x_+ \in \m_+$ and $x \in \g$, there exists $x' \in \g$ such that
$[(x_+, \, 0), \; (r^\flat_+(x), \; r^\flat_-(x))] = (\rfp(x'), \; \rfm(x'))$, 
which is equivalent to $[x_+, \, \rfp(x)] = \rfp(x')$ and $\rfm(x') = 0$, which, in turn, are equivalent to 
$[x_+, \, \rfp(x)] = x'$ and $x' \in \f_+^\perp$. Thus $\m_+ \oplus 0 \subset \n_{\gog} (\l_r)$ if and only if 
$[\m_+, \f_+]\subset \f_+^\perp$. Similarly, 
$\m_- \oplus 0 \subset \n_{\gog} (\l_r)$ if and only if 
$[\m_-, \f_-]\subset \f_-^\perp$. This shows that \eqref{eq-on-m-00} is equivalent to \eqref{eq-on-m-01}.

As $\t_{\rm diag} + \l_r  = \m_+ \oplus \m_-$, it
is also clear that \eqref{eq-on-m-01} is equivalent to $\t_{\rm diag} \subset \n_{\gog} (\l_r)$. 
Using the identification of the quadratic Lie algebra $(\gog, \lara_{\gog})$ with the
Drinfeld double of the Lie bialgebra $(\g, \delta_r)$, one sees that 
$\t_{\rm diag} \subset \n_{\gog} (\l_r)$ if and only if $\ad_\xi^* x = 0$ for all $\xi \in \g^*$ and $x \in \t$,
where $\ad_\xi^*$ is the co-adjoint action of $\xi \in \g^*$ on $\g$, and $\g^*$ has the Lie bracket defined by the dual map of $\delta_r$. It follows that 
$\t_{\rm diag} \subset \n_{\gog} (\l_r)$ if and only if $\delta_r(x) = 0$ for all $x \in \t$.
Thus \eqref{eq-on-m-00} is equivalent to \eqref{eq-on-m-02}.
\end{proof}

\bpr{pr-regular-pi}
Let $(G, r, Y, \lambda)$ be strongly admissible and assume furthermore that
\begin{equation}\lb{eq-on-m-0}
\l_r \subset \m_+ \oplus \m_- \subset \n_{\gog} (\l_r).
\end{equation}
Then $\T$ acts on $(Y, \lambda(r))$ by Poisson isomorphisms, and the Poisson structure $\lambda(r)$ is regular on every 
non-empty intersection $\OO$ of an $M_+$-orbit $\O_+$ and an $M_-$-orbit $\O_-$ in $Y$. 
\epr

\begin{proof}
Fix a pair $(\O_+, \O_-)$ of $M_+$- and $M_-$-orbits in $Y$ such that $\OO \neq \emptyset$, and let $\O$ be the unique
$G$-orbit in $Y$ containing $\OO$.
Let $y_1, y_2 \in \OO$, and let $m_+ \in M_+$ and $m_- \in M_-$ be such that $y_2 = m_+ y_1 = m_- y_1$.
By \reref{re-ly-2} and \eqref{eq-l-yy-yy} and the fact that 
$\left(\Ad_{m_+} \oplus \Ad_{m_-}\right) (\l_r) = \l_r$, 
\begin{align*}
\dim \O - \dim ({\rm Im} (\pi_{y_2}^\#))& = \dim (\l_r \cap \l_{y_2, y_2}) =
\dim \left(\l_r \cap \left(\Ad_{m_+} \oplus \Ad_{m_-}\right)(\l_{y_1, y_1}) \right) \\
&= \dim (\l_r \cap \l_{y_1, y_1})  = \dim \O -\dim ({\rm Im} (\pi_{y_1}^\#)).
\end{align*}
This proves \prref{pr-regular-pi}.
\end{proof}

\bex{ex-double-regular}
When 
$r =r_{(\u, \u')} \in \gotg$ is the $r$-matrix on a quadratic Lie algebra $(\g, \lara_\g)$
defined by a Lagrangian splitting $\g = \u + \u'$ (see \exref{ex-double-0}),
one can take $\m_- = \u$ and $\m_+ = \u'$, so that $\l_r = \m_+ \oplus \m_-$. \prref{pr-regular-pi} in this case was proved in 
\cite[Theorem 2.7]{Lu-Yakimov:DQ}. In the case of 
$\m_- = \n_\g(\u)$ and $\m_+ = \n_\g(\u')$, the normalizer subalgebras in $\g$ of $\u$ and $\u'$ respectively, 
one has $\m_+ \oplus \m_- = \n_{\gog} (\l_r)$, and
\prref{pr-regular-pi} in this case was proved in 
\cite[Proposition 2.13]{Lu-Yakimov:DQ}.
\hfill$\diamond$ \eex

\subsection{$\T$-leaf stabilizers}\lb{subsec-leaf-stabi}
Let again $G$ be a connected Lie group with Lie algebra $\g$, and let $r$ be a
factorizable quasitriangular $r$-matrix on $\g$. Recall the symmetric bilinear form $\lara_\g$ on $\g$ associated to $r$ and
the two Lie subalgebras $\f_+$ and $\f_-$ of $\g$. 

\bde{de-r-mm-admi}
A pair $(M_+, M_-)$ of Lie subgroups  of $G$ is said to be {\it $r$-admissible} if they are connected and 
their respective Lie algebras $\m_+$ and $\m_-$ satisfy
\begin{equation}\lb{eq-on-m-1}
\f_+ \subset \m_+, \hs \f_- \subset \m_-, \hs [\m_+, \, \m_+] \subset \m_+^\perp, \hs [\m_-, \, \m_-] \subset \m_-^\perp.
\end{equation}
\ede

Let $(M_+, M_-)$ be a pair of $r$-admissible Lie subgroups of $G$ with respective Lie algebras $\m_+$ and $\m_-$.
Recall that $\T$ is the connected component of $M_+ \cap M_-$ containing the
identity element and $\t = \m_+ \cap \m_-$ is the Lie algebra of $\T$.
Since
$\m_+^\perp \cap \m_-^\perp = (\m_+ + \m_-)^\perp = \g^\perp = 0$, 
\eqref{eq-on-m-1} implies that $\t = \m_+ \cap \m_-$ is abelian, and thus $\T$ is also
abelian. 
By \leref{le-on-m-equi},   \eqref{eq-on-m-1} implies  \eqref{eq-on-m-0}.
If $(G, r, Y, \lam)$ is a strongly admissible quadruple, then \prref{pr-regular-pi} applies and $\T$ acts on $(Y, \lambda(r))$ by Poisson isomorphisms.
Thus $(Y, \pi = -\lambda(r), \lambda)$ is a $\T$-Poisson manifold. 
Recall from \deref{de-leaf-stabilizer} that for each $y \in Y$, one has the leaf stabilizer
\[
\t_y = \{x \in \t: \, \lambda_y(x) \in {\rm Im}(\pi_y^\#)\} \subset \t.
\]
By \leref{le-leaf-stabilizer}, $\t_{y_1} = \t_{y_2}$ for any $y_1, y_2$ on the same $\T$-leaf of $\pi$ in $Y$. 
To compute $\t_y$ for $y \in Y$, consider the projection 
\begin{equation}\lb{eq-pmt}
p_\t:  \;\; \m_+ \oplus \m_- =\t_{\rm diag} + \l_r \lrw \t,  \;\;\; (x, x) + (\rfp(x'), \, \rfm(x')) \longmapsto x, 
\hs \;\;\; x \in \t, \,x' \in \g.
\end{equation}
Note that $p_\t$ is the restriction to $\m_+ \oplus \m_-$ of the projection
\[
p_\g: \;\; \g \oplus \g  = \g_{\rm diag} + \l_r \lrw \g, \;\;\; (x, x) + (\rfp(x'), \, \rfm(x')) \longmapsto x, 
\hs \;\;\; x \in \g, \,x' \in \g.
\]
Using the fact that $\rfp - \rfm = {\rm Id}_\g$, one sees that $p_\g: \g \oplus \g \to \g$
is also given by 
\begin{equation}\lb{eq-gg-g}
p_\g(x_1, \,x_2) =-\rfm(x_1) + \rfp(x_2), \hs x_1, x _2 \in \g.
\end{equation} 

Let $(\O_+, \O_-)$ be any pair of $M_+$- and $M_-$-orbits contained in the same $G$-orbit in $Y$, possibly $\OO = \emptyset$.  
Let $y_+ \in \O_+$ and $y_- \in \O_-$ be arbitrary, and define
\begin{equation}\lb{eq-ty-0}
\tOO= p_\t((\m_+ \oplus \m_-) \cap \l_{y_+, y_-}),
\end{equation}
where $\l_{y_+, y_-}$ is the Lagrangian subalgebra of $(\gog, \lara_{\gog})$ defined in \eqref{eq-l-yy}.

\bpr{pr-ty-0}
Assume that $(G, r, Y, \lambda)$ is strongly admissible and that $(M_+, \!M_-)$ is an $r$-admissible pair of Lie subgroups of $G$.
Then for any pair $(\O_+, \O_-)$ of $M_+$- and $M_-$-orbits contained in the same $G$-orbit in $Y$, $\tOO \subset \t$ is independent of the choices
of  $y_+ \in \O_+$ and $y_- \in \O_-$. Moreover, when $\OO \neq \emptyset$, one has $\tOO = \t_y$ for
any $y \in \OO$.
\epr

\begin{proof}
Let $y_+^\prime =m_+ y_+ \in \O_+$ and $y_-^\prime  = m_- y_- \in \O_-$, where
$(m_+, m_-) \in M_+ \times M_-$. Then
\[
(\m_+ \oplus \m_-) \cap \l_{y_+^\prime, y_-^\prime} = 
(\m_+ \oplus \m_-) \cap \Ad_{(m_+, m_-)} \l_{y_+, y_-} = \Ad_{(m_+, m_-)} \left((\m_+ \oplus \m_-) \cap \l_{y_+, y_-}\right).
\]
Let $x \in \t$. Then $x \in p_\t((\m_+ \oplus \m_-) \cap \l_{y_+^\prime, y_-^\prime})$ if and only if
$\Ad_{(m_+^{-1}, m_-^{-1})}(x, x) \in \l_r + \l_{y_+, y_-}$.
Writing 
\[
\Ad_{(m_+^{-1}, m_-^{-1})}(x, x) = (x, \; x) + \left(\Ad_{m_+^{-1}} (x) -x, \;\;  \Ad_{m_-^{-1}}(x) -x\right)
\]
and noting that \eqref{eq-on-m-1} implies 
$\displaystyle \left(\Ad_{m_+^{-1}} (x) -x, \;  \Ad_{m_-^{-1}} (x) -x\right) \in \m_+^\perp \oplus 
\m_-^\perp \subset \l_r$, one has
\[
\Ad_{(m_+^{-1}, m_-^{-1})}(x, x) \in \l_r + \l_{y_+, y_-}\hs \mbox{iff} \hs (x, x) \in \l_r + \l_{y_+, y_-},
\]
i.e., $x \in p_\t((\m_+ \oplus \m_-) \cap \l_{y_+, y_-})$. Thus $\tOO$ given in \eqref{eq-ty-0} is independent of the
choice of $y_+ \in \O_+$ and $y_- \in \O_-$.

Assume now that $\OO \neq \emptyset$ and let $y \in \OO$. 
Let $x \in \t$.  
By (ii) of \leref{le-ly}, $x \in \t_y$ if and only if $(x, x) \in \t_{\rm diag} \cap (\l_r + \l_{y, y})$, which is equivalent to
$x \in p_\t((\m_+ \oplus \m_-) \cap \l_{y, y}) \in \tOO$. Thus $\t_y = \tOO$.
\end{proof}

\subsection{Proof of \thref{th-main}}\lb{subsec-proof-th-main}
\thref{th-main}, which is a summary of the main results of $\S$\ref{sec-orbits}, now follows directly from
\prref{pr-delOO}, \prref{pr-ty-0} and
\reref{re-T-Pfaffian}.

For the convenience of the reader, we 
repeat the statements of \thref{th-main} here: Let $(G, r, Y, \lambda)$ be 
a strongly admissible quadruple and let $(M_+, M_-)$ of be an $r$-admissible pair
of Lie subgroups of $G$ with respective Lie algebras $\m_+$ and $\m_-$.
For each pair $(\O_+, \O_-)$ of $M_+$- and $M_-$-orbits contained in the same $G$-orbit in $Y$, 
let $\delOO$ be the integer defined in \eqref{eq-delOO-3} and let
the subspace $\tOO$ of $\t = \m_+ \cap \m_-$ be defined in \eqref{eq-ty-0}.
Let $\T$ be the connected component of $M_+ \cap M_-$ containing the identity element.
If $\delOO =0$ for every pair $(\O_+, \O_-)$ of $M_+$- and $M_-$-orbits contained in the same $G$-orbit in $Y$, then the
six-tuple $(G, r, Y, \lambda, M_+, M_-)$  is 
admissible, i.e., the $\T$-leaves of $\lambda(r)$ in $Y$ are
precisely the connected components of the non-empty intersections $\OO$ of $M_+$ and $M_-$-orbits in $Y$; Moreover,
the leaf stabilizer of each $\T$-leaf in $\OO$ is $\tOO$, and the co-rank of the Poisson structure $\pi = -\lambda(r)$ in $\OO$ is equal to the co-dimension of $\tOO$ in $\t$.

\subsection{Homogeneous admissible six-tuples}\lb{subsec-homog}
Consider  a homogeneous strongly admissible quadruple $(G, r, G/Q, \lam_\sGQ)$ (see 
\exref{ex-homog}), where
$Q$ is a closed and connected Lie subgroup of $G$ whose Lie algebra $\q$ satisfies 
\[
[\q, \, \q] \subset \q^\perp \subset \q.
\]
Here recall that $\q^\perp$ is defined in \eqref{eq-v-perp} using the symmetric non-degenerate bilinear form
$\lara_\g$ on $\g$ associated to  $r$. 
We first state a consequence of the assumption that $\q^\perp \subset \q$.

\bln{ln-q-decom}
For any subspace $\c$ of $\q$ such that $\q = \c + \q^\perp$ is a direct sum decomposition,
the restriction
of $\lara_\g$ on $\c$ is non-degenerate, and  $\g = \c + \c^\perp$
is a direct sum decomposition. Let $p_{\c}: \g \to \c$ be the projection with respect to the decomposition
$\g = \c + \c^\perp$. Then $p_\c \circ \Ad_g = \Ad_g\circ  p_\c: \g \to \c$, where $g$ is any element in the 
normalizer $N_G(\c)$ of $\c$ in $G$;
\eln

\begin{proof}
If $x \in \c \cap \c^\perp$, then $x \in \c \cap \q^\perp$, so $x = 0$. The second statement follows from the
fact that $\Ad_g \c^\perp = \c^\perp$ for every $g \in N_G(\c)$.
\end{proof}

\ble{le-as-homog}
Let $\q = \c + \q^\perp$ be a direct sum decomposition. If $y_+, y_- \in G/Q$ are of the form 
$y_+ = g_+Q$ and $y_- = g_- Q$, where $g_+, g_- \in N_G(\c)$, 
then the Lagrangian subalgebra  $\l_{y_+, y_-}$ of $(\gog, \lara_{\gog})$
defined in \eqref{eq-l-yy} is given by 
\begin{equation}\lb{eq-lyy-homog}
\l_{y_+, y_-} = \left\{(x_+, \, x_-) \in \Ad_{g_+} \q  \oplus \Ad_{g_-} \q: \;\; p_\c(x_+) = \Ad_{g_+g_-^{-1}} p_\c(x_-)\right\}.
\end{equation}
\ele

\begin{proof}
\leref{le-as-homog} follows directly from the the definition of $\l_{y_+, y_-}$ and the decompositions  
\[
\q_{y_+} =\Ad_{g_+} \q = \c + \Ad_{g_+} \q^\perp \hs \mbox{and} \hs \q_{y_-} =\Ad_{g_-} \q = \c + \Ad_{g_-} \q^\perp.
\]
\end{proof}

Let again $(M_+, M_-)$ be a pair of $r$-admissible Lie subgroups of $G$ with respective Lie algebras $\m_+$ and $\m_-$.
If every $M_+$-orbit $\O_+$ and every $M_-$-orbit $\O_-$ in $G/Q$ contain elements of the form $gQ$ with
$g \in N_G(\c)$, one can use \eqref{eq-lyy-homog} to compute the 
integers $\delOO$ and the subspaces $\tOO$ of $\t$, as in the 
following \prref{pr-as-homog}. 

\bpr{pr-as-homog}
Suppose that $\q = \c + \q^\perp$ is a direct decomposition, and assume that 
all $(M_+, Q)$- and $(M_-, Q)$-double cosets in $G$ contain elements in $N_G(\c)$.
Let $(\O_+, \O_-)$ be any pair of $M_+$ and $M_-$-orbits in $G/Q$, and choose any
$g_+, g_- \in N_G(\c)$ such that $g_+Q \in \O_+$ and $g_-Q \in \O_-$. Then
\begin{align*}
\delOO &= \dim \left(p_\c(\m_+^\perp \cap \Ad_{g_+}\q) \cap 
\Ad_{g_+g_-^{-1}}\left(p_\c(\m_-^\perp \cap \Ad_{g_-} \q)\right)\right),\\
\tOO & = p_\t (V_{g_+, g_-}),
\end{align*}
where  $V_{g_+, g_-} =\displaystyle \left\{(x_+, x_-) \in (\m_+ \cap  \Ad_{g_+} \q) \oplus (\m_-\cap  \Ad_{g_-} \q): 
p_\c(x_+) = \Ad_{g_+g_-^{-1}} p_\c(x_-)\right\}$. In particular, the 
six-tuple $(G, r, G/Q, \lam_{\sGQ}, M_+, M_-)$ is  admissible if
\[
p_\c(\m_+^\perp \cap \Ad_{g}\q) \cap 
\Ad_h\left(p_\c(\m_-^\perp \cap \Ad_k \q)\right) = 0, \hs \hs \forall \; g, h, k \in N_G(\c).
\]
\epr

\begin{proof} \prref{pr-as-homog} follows directly from 
\eqref{eq-delOO-3} for $\delOO$ and \eqref{eq-ty-0} for $\tOO$.
\end{proof}

\bre{re-lag-c}
Note that if $\q$ is Lagrangian with respect to $\lara_{\g}$, 
by taking $\c = 0$, the six-tuple $(G, r, G/Q, \lam_{\sGQ}, M_+, M_-)$ is automatically 
admissible, and the leaf stabilizer $\tOO$ is given by
\[
\tOO = p_\t\left((\m_+ \cap  \Ad_{g_+} \q) \oplus (\m_-\cap  \Ad_{g_-} \q)\right)
\]
for any $g_+Q \in \O_+$ and $g_-Q \in \O_-$.
\hfill$\diamond$ 
\ere

\bex{ex-G-MM}
Let $(G, \piG)$ be a connected Poisson Lie group, where $\piG = r^L - r^R$ for a
factorizable quasitriangular $r$-matrix $r$ on the Lie algebra $\g$ of $G$. Equip the direct product Lie algebra $\gog$ with the
direct product factorizable quasitriangular $r$-matrix $(r, -r)$, and let $\lam$ be the left action of $G \times G$ on $G$ given by
\[
(g_1, \, g_2) \cdot g = g_1gg_2^{-1}, \hs g_1, \, g_2, \, g \in G.
\]
As the stabilizer subgroup of $\lam$ at $g \in G$ is $\{(g_1, g^{-1}g_1g): g_1 \in G\}$, which is connected and its Lie algebra
Lagrangian with respect to the the symmetric bilinear form
on $\gog$ associated to $(r, -r)$, 
the quadruple $(G \times G, \;(r, -r), \; G, \; \lam)$ is strongly admissible. It is easy to see that
$-\lam(r, -r) = \piG$. If $(M_+, M_-)$ is a pair of $r$-admissible
Lie subgroups of $G$, then $(M_+ \times M_+,\; M_- \times M_-)$ is
a pair of $(r, -r)$-admissible Lie subgroups of $G \times G$. By \reref{re-lag-c}, the 
six-tuple $(G \times G, \,(r, -r), \, G, \, \lam, \, M_+ \times M_+, \, M_- \times M_-)$
is admissible, and hence the $(\T, \T)$-leaves of $(G, \piG)$ are precisely the connected components of 
non-empty intersections of $(M_+, M_+)$-double cosets and $(M_-, M_-)$-double cosets in $G$.
As a special case, if $r$ is defined by a Lagrangian splitting $\g = \u + \u'$,
by taking $M_+ = U'$ and $M_- = U$, where $U$ and $U$ are the connected Lie subgroups of $G$ with Lie algebras 
$\u$ and $\u'$ respectively, the symplectic leaves of 
$(G, \piG)$ are precisely the connected components of $(U, U)$-double cosets and $(U', U')$-double cosets
in $G$ (see, for example, \cite{STS2}). 
\eex

\sectionnew{Mixed product Poisson structures associated to admissible quadruples}\lb{sec-mixed}
\subsection{The construction}\lb{subsec-mixed}
Let $r$ be a quasitriangular $r$-matrix on a Lie algebra $\g$, and let
$n \geq 1$ be any integer.
Writing $r =\sum_i x_i \ot y_i \in \gotg$, we introduced in \cite{Lu-Mou:mixed} a 
quasitriangular $r$-matrix
$r^{(n)}$ on the direct product Lie algebra $\g^n = \g \oplus \cdots \oplus \g$ ($n$-copies), which is given by
\begin{equation}\lb{eq-r-(n)}
r^{(n)} = \sum_{1 \leq j \leq n,  \; j\, {\rm is \;odd}} (r)_j \;\;+ \sum_{ 1 \leq j \leq n,  \;j \;{\rm is \;even}} (-r^{21})_j
\;\; - \sum_{1 \leq j < k \leq n} \sum_i (y_i)_j \wedge (x_i)_k, 
\end{equation}
where for $X \in \g^{\otimes l}$, $l = 1, 2$, and $1 \leq j \leq n$, $(X)_j \in (\g^n)^{\ot l}$ is the image of $X$ under the embedding of $\g$ into $\g^n$ as the $j$'th summand.  
Since the symmetric part of $r^{(n)}$ is $(s, -s, s, -s, \ldots)$, where $s$ is the symmetric part of $r$, if $r$ is factorizable, so is $r^{(n)}$.

Assume that $(G, r, Y_i, \lambda_i)$ is an admissible quadruple for each $1 \leq i \leq n$, and consider the product manifold 
$Y = Y_1 \times \cdots \times Y_n$ and the 
direct product action $\lambda$ of $\g^n$ on $Y$, i.e.,  
\[
\lambda: \g^n \lrw \V^1(Y), \;\; \lambda(x_1, x_2, \ldots, x_n) = (\lambda_1(x_1), \; \lambda_2(x_2),
\ldots, \; \lambda_n(x_n)), \hs x_j \in \g.
\]
Then the quadruple $(G^n, r^{(n)}, Y, \lambda)$ is admissible. Moreover, when each $(G, r, Y_i, \lambda_i)$ is strongly admissible, so
is $(G^n, r^{(n)}, Y, \lambda)$. Being the sum of 
the direct product Poisson structure $(-\lam_1(r), \ldots, -\lam_n(r))$ on $Y$ and some ``mixed terms", 
the Poisson structure 
$-\lambda(r^{(n)})$
on $Y$ is an example of a {\it mixed product Poisson structure} on the product manifold $Y$ (see \cite{Lu-Mou:mixed}).

In $\S$\ref{sec-mixed}, we apply the general theory in $\S$\ref{sec-orbits} to the admissible quadruples
$(G^n, r^{(n)}, Y, \lambda)$.
We first establish a property of the quasitriangular $r$-matrix $r^{(n)}$ on $\g^n$. For notational simplicity, set
\begin{equation}\lb{eq-tilde-l}
 \tilde{r} = r^{(n)} \in \g^n \otimes \g^n, \;\;\;\;\;
\tilf_\pm = {\rm Im} (\tilde{r}_{\pm}) \subset \g^n.
\end{equation}
Recall that when $r$ is factorizable, one has the Lie subalgebra $\l_r$ of $\gog$ given by
\[
\l_r = \{(r_+(\xi), \, r_-(\xi)): \xi \in \g^*\} =\{(\rfp(x), \, \rfm(x)): \; x \in \g\}\subset \gog,
\]
and the Lie subalgebras $\f_\pm = {\rm Im}(r^\flat_\pm)$ of $\g$. 
Let $\tau$ be the automorphism of $\g^n$ defined by
\[
\tau(x_1, \, x_2, \, x_3, \, \ldots, \, x_{n-1}, \, x_n) = (x_1, \, x_3, \, \ldots, \, x_{n-1}, \, x_n, \, x_2),
\hs x_j \in \g.
\]

\ble{le-tilde-all} Assume that $r$ is factorizable. Then
\begin{align*} 
\tilf_+ & =  
\tau (\l_r \oplus \overbrace{\g_{\rm diag} \oplus \cdots \oplus \g_{\rm diag}}^{m-1}),\hs
\tilf_-   = \overbrace{\g_{\rm diag} \oplus \cdots \oplus \g_{\rm diag}}^m,   \hs \mbox{if} \;\; n = 2m \; \mbox{is even},\\
\tilf_+ & = \f_+ \oplus \overbrace{\g_{\rm diag} \oplus \cdots \oplus \g_{\rm diag}}^m, \hs
\tilf_- = \overbrace{\g_{\rm diag} \oplus \cdots \oplus \g_{\rm diag}}^m  \oplus \;\f_-, \hs
\mbox{if} \;\; n = 2m+1 \; \mbox{is odd}.
\end{align*}
\ele

\begin{proof}
Assume first that $n = 2m$ is even, and, for notational simplicity, 
set 
\[
\tilf_+^\prime =  
\tau (\l_r \oplus \overbrace{\g_{\rm diag} \oplus \cdots \oplus \g_{\rm diag}}^{m-1}) \hs \mbox{and} \hs
\tilf_-^\prime   = \overbrace{\g_{\rm diag} \oplus \cdots \oplus \g_{\rm diag}}^m.
\]
Let
$\tilde{\xi} = (\xi_1, \xi_2, \ldots, \xi_n) \in 
(\g^*)^n \cong (\g^n)^*$, and write 
\begin{equation}\lb{eq-tilr-xi}
\tilr_+(\tilde{\xi}) = (x_1, x_2, \ldots, x_n), \hs \tilr_-(\tilde{\xi}) = (y_1, y_2, \ldots, y_n).
\end{equation}
By the definition of $\tilr$, one has, for $1 \leq j \leq n$,
\begin{align}\lb{eq-xy-j}
x_j &= \begin{cases} r_-(\xi_1 + \cdots + \xi_{j-1}) + r_+(\xi_j + \xi_{j+1} + \cdots + \xi_n), & \;\;
j \; \;\mbox{odd}, \\
r_-(\xi_1 + \cdots + \xi_{j-1} + \xi_j) + r_+(\xi_{j+1} + \cdots + \xi_n), & \;\;
j \; \;\mbox{even}\end{cases}, \\
\nonumber
y_j &= \begin{cases} r_-(\xi_1 + \cdots + \xi_{j-1}+ \xi_j) + r_+(\xi_{j+1} + \cdots + \xi_n), & \;\;
j \; \;\mbox{odd}, \\
r_-(\xi_1 + \cdots + \xi_{j-1}) + r_+(\xi_j +\xi_{j+1} + \cdots + \xi_n), & \;\;
j \; \;\mbox{even}\end{cases}.
\end{align}
 As $x_1 = r_+(\xi_1 + \cdots + \xi_n)$, $x_n = r_-(\xi_1 + \cdots + \xi_n)$, and for $1 \leq k \leq m-1$,
\begin{align*}
x_{2k} = x_{2k+1}& = r_-(\xi_1 + \cdots + \xi_{2k}) + r_+(\xi_{2k+1} + \cdots + \xi_{n})\\ 
& = r_+(\xi_1 + \cdots + \xi_{n}) + (r_- - r_+) (\xi_1 + \cdots + \xi_{2k}),
\end{align*}
one has $\tilr_+(\tilde{\xi}) \in \tilf_+^\prime$. 
Moreover, since $r_- - r_+: \g^* \to \g$ is invertible, $\tilr_+(\tilde{\xi})=0$ if
and only if  $\xi_{2k-1} + \xi_{2k} = 0$ for every $1 \leq k \leq m$. Thus $\dim \ker \tilr_+ = m (\dim \g)$. It follows that
$\dim ({\rm Im (\tilr_+})) = \dim (\tilf_+^\prime)$ and hence $\tilf_+ = \tilf_+^\prime$. 
Similarly, since
\[
y_{2k-1} = y_{2k} = r_-(\xi_1 + \cdots + \xi_{2k-1}) + r_+(\xi_{2k} + \cdots + \xi_n) = 
(r_- - r_+)(\xi_1 + \cdots \xi_{2k-1}) + r_+(\xi_1 + \cdots + \xi_n)
\]
for every $1 \leq k \leq m$, 
 $\tilr_-(\tilde{\xi}) \in \tilf_-^\prime$. Moreover, 
$\tilr_-(\tilde{\xi})=0$ if and only if $y_1 = 0$ and $y_j-y_{j+1} = 0$ for every $1 \leq j \leq n-1$, which is equivalent to
$r_-(\xi_1) + r_+(\xi_n) = 0$ and $\xi_{2k} + \xi_{2k+1} = 0$ for every $1 \leq k \leq m-1$. Since
the map $(\g^*)^2 \to \g, (\xi, \eta) \mapsto r_-(\xi) + r_+(\eta)$ is surjective, one sees that 
$\dim ({\rm Im (\tilr_-})) = \dim (\tilf_-^\prime)$ and hence $\tilf_- = \tilf_-^\prime$. 

Assume that $n = 2m +1$ is odd, and set
$\displaystyle \tilf_+^\prime = \f_+ \oplus \overbrace{\g_{\rm diag} \oplus \cdots \oplus \g_{\rm diag}}^m$
and $\displaystyle \tilf_-^\prime = \overbrace{\g_{\rm diag} \oplus \cdots \oplus \g_{\rm diag}}^m  \oplus \;\f_-$. Let 
$\tilde{\xi} = (\xi_1, \ldots, \xi_n)\in (\g^*)^n$ with $\tilr_+(\tilde{\xi})$ and
$\tilr_-(\tilde{\xi})$ as given in \eqref{eq-tilr-xi} and \eqref{eq-xy-j}. Again it is clear that 
$\tilr_+(\tilde{\xi}) \in \tilf_+^\prime$. Moreover, $\tilr_+(\tilde{\xi}) =0$ if and only if $r_+(\xi_n) = 0$ and 
$\xi_{2k-1} + \xi_{2k} = 0$ for every $1 \leq k \leq m$. It follows by dimension counting that 
$\tilf_+ = \tilf_+^\prime$. Similar arguments show that $\tilf_- = \tilf_-^\prime$. 
\end{proof}

Assume again that $r$ is factorizable, and let $(M_+, M_-)$ be an $r$-admissible pair of Lie subgroups of $G$
with respective Lie algebras $\m_+$ and $\m_-$.
Let $M_+^{(n)}, M^{(n)}_- \subset G^n$ be given by
\begin{align*}
M^{(n)}_+ & = M_+ \times \overbrace{\Gdia \times \cdots \times \Gdia}^{m-1} \times M_-, \hs
M^{(n)}_- = \overbrace{\Gdia \times \cdots \times \Gdia}^m, \hs \mbox{if} \;\; n = 2m \; \mbox{is even},\\
M^{(n)}_+ & = M_+ \times \overbrace{\Gdia \times \cdots \times \Gdia}^{m}, \hs
M^{(n)}_- = \overbrace{\Gdia \times \cdots \times \Gdia}^m \times M_-, \hs \mbox{if} \;\; n = 2m+1 \; \mbox{is odd}.
\end{align*}
Then their Lie algebras $\m^{(n)}_+$ and $\m^{(n)}_-$  are respectively given by
\begin{align*} 
\m^{(n)}_+ & =  
\m_+ \oplus \overbrace{\g_{\rm diag} \oplus \cdots \oplus \g_{\rm diag}}^{m-1} \oplus \,\m_-,\hs
\m^{(n)}_-   = \overbrace{\g_{\rm diag} \oplus \cdots \oplus \g_{\rm diag}}^m,   \hs \mbox{if} \;\; n = 2m \; \mbox{is even},\\
\m^{(n)}_+ & = \m_+ \oplus \overbrace{\g_{\rm diag} \oplus \cdots \oplus \g_{\rm diag}}^m, \hs
\m^{(n)}_- = \overbrace{\g_{\rm diag} \oplus \cdots \oplus \g_{\rm diag}}^m  \oplus \;\m_-, \hs
\mbox{if} \;\; n = 2m+1 \; \mbox{is odd}.
\end{align*}
It is clear from \leref{le-tilde-all} that the pair $(M^{(n)}_+, M^{(n)}_-)$ of subgroups of $G^n$ is 
$r^{(n)}$-admissible. Let again $\T$ be the connected component of $M_+ \cap M_-$ containing the identity element of $G$. Then
$\T^{(n)} := \{(g, \ldots, g): g \in \T\}$ is the connected component of $M^{(n)}_+ \cap M^{(n)}_-$ containing the identity element of $G^n$.

\subsection{A homogeneous case}\lb{subsec-homog-n} As in $\S$\ref{subsec-homog}, consider 
a six-tuple $(G,  r,  G/Q,  \lam_\sGQ,  M_+,  M_-)$, 
where  $G$ is a connected Lie group with Lie algebra $\g$, $r$ a factorizable quasitriangular
$r$-matrix on $\g$, $Q$ a closed and connected Lie subgroup of $G$ whose Lie algebra $\q$ satisfies 
\begin{equation}\lb{eq-on-q}
[\q, \, \q] \subset \q^\perp \subset \q,
\end{equation} 
and $(M_+, M_-)$  a pair of $r$-admissible Lie subgroups of $G$ with respective Lie algebras $\m_+$ and $\m_-$.
In $\S$\ref{subsec-homog-n}, we consider, for each integer $n \geq 1$, the six-tuple
\[
(G^n, \; r^{(n)}, \; (G/Q)^n, \; \lambda, \; M_+^{(n)}, \; M_-^{(n)}),
\]
where $\lam$ is the direct product of the action $\lam_\sGQ$ of $G$ on each factor of $(G/Q)^n$.
Identify
\[
\T \; \stackrel{\sim}{\lrw} \; \T^{(n)}, \;\;\; t \longmapsto (t, t, \ldots, t), \hs t \in \T.
\]
Then $\T$ acts on $(G/Q)^n$ diagonally.
As the
case of $n = 1$ is covered in \prref{pr-as-homog}, we assume that $n \geq 2$.

\bas{as-homog-n}
There exists a direct sum decomposition $\q = \c + \q^\perp$ such that 

1) every $(K, Q)$-double cosets in $G$, where $K \in \{M_+, M_-, Q\}$, contains elements in $N_G(\c)$;

2) $p_{\c}(\m_+^\perp \cap \Ad_g\q) \cap \Ad_h p_{\c}(\m_-^\perp\cap \Ad_k\q)=0$ for 
all $g, h , k\in N_G(\c)$,

\noindent
where recall that $p_\c: \g \to \c$ is the projection with respect to the decomposition $\g = \c + \c^\perp$.
\eas

\bre{re-as-homog-n}
Note that by taking $\c = 0$, \asref{as-homog-n} holds automatically if $\q = \q^\perp$.
\hfill$\diamond$
\ere

\bnota{nota-gh} For  
$g= (g_1, \ldots, g_n) \in G^n$, let 
$\underline{g} = (g_1Q, \, g_2Q, \, \ldots, \, g_nQ) \in (G/Q)^n$. 
Let
$e \in G$ be the identity element $G$. Let $(\widetilde{\O}_+, \widetilde{\O}_-)$ be
any pair of $M^{(n)}_+$ and $M^{(n)}_-$-orbits in $(G/Q)^n$.
\asref{as-homog-n} implies that 
 there exist
$g, h \in (N_G(\c))^n$ of the form 
\begin{align}\lb{eq-gh}
g &= \begin{cases} (g_1, \, e,\, g_3, \, e, \,g_5, \,\ldots, \,e, \,g_{2m-1}, \,g_{2m}), & \;\mbox{if} \; n = 2m \; \mbox{is even}\\
(g_1, \, e, \,g_3,\, e,\, g_5,\, \ldots, \, \,e, \,g_{2m-1},\,e, \,g_{2m+1}), & \;\mbox{if} \; n = 2m+1 \; \mbox{is even}\end{cases}, \\
\lb{eq-gh-1}
h &= \begin{cases} (e, \,h_2, \,e, \,h_4, \,\ldots, \, e, \,h_{2m-2}, \,e, \,h_{2m}), & \;\mbox{if} \; n = 2m \; \mbox{is even}\\
(e, \,h_2,\, e, \,h_4, \,\ldots,\, e, \,h_{2m-2}, \, e, \,h_{2m},\, h_{2m+1}), & \;\mbox{if} \; n = 2m+1
\; \mbox{is odd}\end{cases},
\end{align} 
such that $\underline{g} \in \widetilde{\O}_+$ and $\underline{h} \in \widetilde{\O}_-$. 
With $g, h \in (N_G(\c)^n$ so chosen, let
\begin{equation}\lb{eq-gh-2}
g \bowtie h = \begin{cases} (g_1, \, h_2,\, g_3, \, h_4, \,\ldots,  \,g_{2m-1}, \,h_{2m}), & \;\mbox{if} \; n = 2m \; \mbox{is even}\\
(g_1, \, h_2, \,g_3,\, h_4,\,  \ldots, \,g_{2m-1},\,h_{2m}, \,g_{2m+1}), & \;\mbox{if} \; n = 2m+1 \; \mbox{is odd}\end{cases},
\end{equation}
and let $(g \bowtie h)_{n+1} = g_n$ if $n$ is even and $(g \bowtie h)_{n+1} = h_n$ if $n$ is odd.
Let $c = (c_1, c_2, \ldots, c_{n+1})$, where $(c_1, c_2, \ldots, c_n) = g \bowtie h$ and $c_{n+1}=
(g \bowtie h)_{n+1}$, and let 
\begin{equation}\lb{eq-V-gh}
V_c = \left\{(x_+, x_-) \in \!\left(\m_+ \cap \Ad_{c_1}\q\right) \oplus \left(\m_- \cap
\Ad_{c_{n+1}}\q\right): p_\c(x_+) = \Ad_{c_1c_2 \cdots c_n c_{n+1}^{-1}} (p_\c(x_-))\right\}.
\end{equation}
Recall again that $\t = \m_+ \cap \m_-$, and one has the projection $p_\t: \m_+ \oplus \m_- \to \t \cong \t_{\rm diag}$
with respect to the decomposition $\m_- \oplus \m_- = \t_{\rm diag} + \l_r$. 
\enota

\bpr{pr-D} Under \asref{as-homog-n}, 

(1) the six-tuple $(G^n, \, r^{(n)}, \, (G/Q)^n, \, \lambda, \, M_+^{(n)}, \, M_-^{(n)})$ is  admissible
for every $n \geq 2$;

(2) for any pair $(\widetilde{\O}_+, \widetilde{\O}_-)$ of $M^{(n)}_+$ and $M^{(n)}_-$-orbits in $(G/Q)^n$, the leaf
stabilizer of every $\T$-leaf in $\widetilde{\O}_+ \cap \widetilde{\O}_-$ is given by
$\displaystyle \t_{\widetilde{\sO}_+, \widetilde{\sO}_-} = p_\t \left(V_c\right) \subset \t$ with $V_c \subset \m_+ \oplus \m_-$ given in \eqref{eq-V-gh}.
\epr

\begin{proof} Let $(\widetilde{\O}_+, \widetilde{\O}_-)$ be an arbitrary pair of $M^{(n)}_+$ and $M^{(n)}_-$-orbits in $(G/Q)^n$,
and let $g, h \in (N_G(\c))^n$ be as in \eqref{eq-gh} and \eqref{eq-gh-1} such that $y_+=\underline{g} \in \widetilde{\O}_+$
and $y_- = \underline{h} \in \widetilde{\O}_-$. We use the description in  \leref{le-as-homog} of
the Lagrangian subalgebra $\l_{y_+, y_-}$ of $\g^n \oplus \g^n$ to show that
$\delta_{\widetilde{\sO}_+, \widetilde{\sO}_-} = 0$
and to compute the subspace $\t_{\widetilde{\sO}_+, \widetilde{\sO}_-}$ of $\t = \m_+ \cap \m_-$. 
For notational simplicity, set $\widetilde{M}_\pm :=M_\pm^{(n)}$ and let $\tilde{\m}_\pm$ be their respective
Lie algebras.

Assume first that
$n = 2m$ is even. By \leref{le-as-homog},  
$\displaystyle \left(\tilde{\m}_+^\perp \oplus \tilde{\m}_-^\perp\right)\cap \l_{y_+, y_-}$ consists of  
elements $a \in \g^n \oplus \g^n$ of the form
\begin{equation}\lb{eq-a}
a = (x_1 + x_1^\prime, \; x_2 + x_2^\prime, \; \ldots, \; x_n + x_n^\prime, \;
z_1 + z_1^\prime, \; z_2 + z_2^\prime, \; \ldots, \, z_n + z_n^\prime),
\end{equation}
where, by writing $g = (g_1, g_2, \ldots, g_n)$ and $h = (h_1, h_2, \ldots, h_n)$,
\[
x_j, \; z_j \in \c, \hs x_j^\prime \in \Ad_{g_j}\q^\perp, \; z_j^\prime \in \Ad_{h_j}\q^\perp, \hs x_j = \Ad_{g_jh_j^{-1}} z_j, \hs j = 1, \ldots, n,
\]
$x_1 + x_1^\prime \in \m_+^\perp \cap \Ad_{g_1}\q$,
$\;x_{2m} + x_{2m}^\prime \in \m_-^\perp \cap \Ad_{g_{2m}}\q$, and
\begin{align*}
&x_{2j} = x_{2j+1}, \;\;\; x_{2j}^\prime = x_{2j+1}^\prime, \hs j = 1, \ldots, m-1,\\
&z_{2j-1} = z_{2j}, \; \;\; z_{2j-1}^\prime = z_{2j}^\prime, \hs j = 1, \ldots, m.
\end{align*}
It follows that
\[
x_1 = \Ad_{g_1 h_2g_3h_4\cdots g_{2m-1}h_{2m}g_{2m}^{-1}}(x_{2m}) \in p_\c(\m_+^\perp\cap \Ad_{c_1}\q) 
\cap \Ad_{c_1c_2 \cdots c_nc_{n+1}^{-1}} p_\c(\m_-^\perp \cap \Ad_{c_{n+1}}\q).
\]
By 2) of \asref{as-homog-n}, $x_1=0$, and it follows that $x_j = 0$ for every $j = 1, \ldots, n$. By \eqref{eq-delOO-3}, 
$\delta_{\widetilde{\sO}_+, \widetilde{\sO}_-} = 0$. The case of $n = 2m+1$ is odd is proved similarly. This proves (1).

To prove (2), note that by \leref{le-tilde-all}, 
\[
\t_{\widetilde{\sO}_+, \widetilde{\sO}_-} = p_\t \left((p\left(\tilde{\m}_+ \oplus \tilde{\m}_-\right))
\cap \l_{y_+, y_-}\right),
\]
where 
$p: \g^n \oplus \g^n \to \gog$ is given by
\[
p(a_1, \, a_2, \, \ldots, \, a_n,\, b_1, \, b_2, \, \ldots, \, b_n) = \begin{cases}
(a_1, \, a_n), & \;\; \mbox{if} \;\; n \; \mbox{is even},\\
(a_1, \, b_n), & \;\; \mbox{if} \;\; n \; \mbox{is odd}\end{cases}.
\]
Replacing $\tilde{\m}_+^\perp \oplus \tilde{\m}_-^\perp$ in the proof of (1) by 
$\tilde{\m}_+ \oplus \tilde{\m}_-$, one sees that 
\[
p\left((\tilde{\m}_+ \oplus \tilde{\m}_-)
\cap \l_{y_+, y_-}\right) \subset V_c.
\]
Assume again that $n = 2m$ is even and let $(x_+, x_-) \in V_c$.
Then there is a unique element  $a \in \left(\tilde{\m}_+ \oplus \tilde{\m}_-\right)
\cap \l_{y_+, y_-}$ of the form \eqref{eq-a} with 
$x_1+x_1^\prime = x_+$,
$x_n + x_n^\prime = x_-$, $x_j^\prime =0$ for $j = 2, \ldots, n-1$, and $z_j^\prime = 0$ for all
$j = 1, \ldots, n$. Moreover, $p(a) =(x_+, x_-)$.  This shows that 
$p\left((\tilde{\m}_+ \oplus \tilde{\m}_-)
\cap \l_{y_+, y_-}\right) = V_c$ when $n$ is even, and by
\prref{pr-ty-0}, one has $\displaystyle \t_{\widetilde{\sO}_+, \widetilde{\sO}_-} = p_\t (V_c) \subset \t$. The case when $n$ is odd is proved similarly. This proves (2). 
\end{proof}

\subsection{Main examples}\lb{subsec-two-main}
Consider again a homogeneous strongly admissible quadruple $(G, r, G/Q, \lam_{\sGQ})$, 
but where we assume that
the Lie algebra $\q$ of $Q$ satisfies
\begin{equation}\lb{eq-on-q-1}
\f_+\subset \q \hs \mbox{and} \hs [\q, \q] \subset \q^\perp.
\end{equation}
Since $\q^\perp \subset \f_+^\perp \subset \f_+ \subset \q$, the quadruple
$(G, r, G/Q, \lam_{\sGQ})$ is indeed strongly admissible.
Moreover, by \reref{re-GQ}, $Q$ is a Poisson Lie subgroup of the Poisson Lie group
$(G, \piG)$, where $\piG = r^L - r^R$,
and  the Poisson structure 
$-\lambda_\sGQ(r)$ on $G/Q$ coincides with the projection $\pi_{\sGQ}$ of $\piG$ to $G/Q$.

Let $n \geq 1$ be an integer, and let $G^n$ acts on itself from the right by \eqref{eq-Gn-Gn-1}. 
Then we have the two quotient manifolds
\[
Y_n = G\times_{Q} \cdots \times_{Q} G/Q \hs \mbox{and} \hs X_n = G \times_{Q} \cdots \times_{Q} G
\]
of $G^n$, each with the quotient Poisson structure, respectively denoted by $\pi_{\sY_n}$ and $\pi_{\sX_n}$,
which, by definition, are the projections of the direct product Poisson structure $\pi_\sG^n$ on $G^n$.
Let $(M_+, M_-)$ be a pair of $r$-admissible Lie subgroups of $G$, and let again $\T$ be the 
connected component of $M_+ \cap M_-$ containing the identity element. Then 
$\T$ acts on $(Y_n, \pi_{\sY_n})$ and $\T^2 = \T \times \T$ acts on $(X_n, \pi_{\sX_n})$ by Poisson isomorphisms
via
\begin{align*}
t \cdot [g_1,\, g_2, \,\ldots, \,g_n]_{\sY_n} &= [tg_1, \,g_2, \,\ldots,\, g_n]_{\sY_n}, \hs t \in \T, \; g_j \in G,\\
(t_1, t_2) \cdot [g_1,\, g_2, \,\ldots,\, g_n]_{\sX_n}& = [t_1g_1,\, g_2, \,\ldots,\, g_nt_2^{-1}]_{\sX_n}, 
\hs t_1, t_2 \in \T, \; g_j \in G.
\end{align*}
In this section, we study the $\T$-leaves of $(Y_n, \pi_{\sY_n})$ and the $\T^2$-leaves of $(X_n, \pi_{\sX_n})$.
As the case for $(Y_1, \pi_{{\scriptscriptstyle Y}_1})$ is covered in \prref{pr-as-homog} and the
case of $(X_1 = G, \pi_{\sX_1} = \piG)$ is covered in \exref{ex-G-MM}, we will assume that $n \geq 2$.

We first look at the case of $(Y_n, \pi_{\sY_n})$.  Consider the diffeomorphism
\begin{equation}\lb{eq-IZ}
J_{\sY_n}: \; Y_n \lrw (G/Q)^n, \;\;\; [g_1, g_2, \ldots, g_n]_{\sY_n} \longmapsto (g_1Q, \; g_1g_2Q, \; \cdots, \; 
g_1g_2 \cdots g_nQ).
\end{equation}
Let $\lam: \g^n \to \V^1((G/Q)^n)$ be again the direct product of the action $\lam_\sGQ$ on each factor.
We have the following crucial \prref{pr-GGG-QQQ} from \cite[$\S$8]{Lu-Mou:mixed}.

\bpr{pr-GGG-QQQ}
As Poisson structures on $(G/Q)^n$, one has $J_{\sY_n}(\pi_{\sY_n}) = -\lam\left(r^{(n)}\right)$.
\epr

We can thus apply \prref{pr-D} to the Poisson structure $\lam\left(r^{(n)}\right)$ on $(G/Q)^n$ and
use 
\[
J_{\sY_n}^{-1}: \;\; (G/Q)^n \lrw Y_n,\;\;\; (k_1Q, \, k_2Q, \, \ldots, \, k_nQ) \longmapsto 
[k_1, \, k_1^{-1}k_2, \, k_2^{-1}k_3, \ldots, \, k_{n-1}^{-1}k_n]_{\sY_n},
\]
to translate the results to the Poisson structure $\pi_{\sY_n}$ on $Y_n$. For
$a = (a_1, a_2, \ldots, a_n) \in G^n$, let
\[
C(a) = (M_+a_1Q) \times_Q (Qa_2Q) \times_Q \cdots \times_Q (Q a_nQ)/Q \subset Y_n.
\]
Let $\mu_{\sY_n}: Y_n \to G/Q$ be given by $\mu_{\sY_n}([a_1, a_2, \ldots, a_n]_{\sY_n}) = a_1a_2 \cdots a_n Q$.
 
\ble{le-transition}
For any $g, h \in G^n$ of the form \eqref{eq-gh} and \eqref{eq-gh-1}, one has
\[
J_{\sY_n}^{-1}\left(\left(M^{(n)}_+ \underline{g}\right)\cap \left(M^{(n)}_- \underline{h}\right)\right) =
C(g \bowtie h) \cap
\mu_{\sY_n}^{-1}\left(M_-(g \bowtie h)_{n+1}Q/Q\right).
\]
\ele

\begin{proof}
Write $g = (g_1, g_2, \ldots, g_n)$ and
$h = (h_1, h_2, \ldots h_n)$, and consider first the case when $n = 2m$ is even. Let $k = (k_1, k_2, \ldots, k_n) \in G^n$. Then 

1) $\underline{k} \in M^{(n)}_+ \underline{g}$ if and only if $k_1 \in M_+ g_1Q$ and
$k_{2i}^{-1}k_{2i+1} \in Q g_{2i}^{-1} g_{2i+1} Q$
for $i = 1, \ldots, m-1$, and $k_{2m} \in M_- g_{2m} Q$, and

2)  $\underline{k}\in M^{(n)}_- \underline{g}$ if and only if $k_{2j-1}^{-1}k_{2j} \in Q g_{2j-1}^{-1} g_{2j} Q$
for $j = 1, \ldots, m$.

\noindent
It follows that $\displaystyle \underline{k} \in 
\left(M^{(n)}_+ \underline{g}\right)\cap \left(M^{(n)}_- \underline{h}\right)$ if and only if
$J_{\sY_n}(\underline{k}) \in C(g \bowtie h) \cap \mu_{\sY_n}^{-1}(M_- g_nQ/Q)$. The case when $n =2m+1$ is odd is proved similarly. 
\end{proof}

\bth{th-Zn} Under \asref{as-homog-n}, the following holds for every $n \geq 2$:

(a) Every non-empty intersection 
\begin{equation}\lb{eq-DQQD}
Y_n(c) = \left((M_+c_1Q)\times_Q (Qc_2Q) \times_Q \cdots \times_Q (Qc_nQ)/Q\right) \cap \mu_{\sY_n}^{-1}(M_-c_{n+1}Q/Q) \subset Y_n,
\end{equation}
where $c = (c_1, \ldots, c_{n+1}) \in (N_G(\c))^{n+1}$, 
is transversal, and their connected components are precisely all the $\T$-leaves of $\pi_{\sY_n}$ in $Y_n$;

(b) The leaf stabilizer for 
every $\T$-leaf in $Y_n(c)$ in \eqref{eq-DQQD} is $p_\t(V_c)$, 
where
\[
V_c = \{(x_+, x_-) \in (\m_+\cap \Ad_{c_1} \q) \oplus (\m_-\cap \Ad_{c_{n+1}} \q):\;\;
p_\c(x_+) = \Ad_{c_1c_2 \cdots c_nc_{n+1}^{-1}} (p_\c(x_-))\}. 
\]

\noindent
Moreover, if $G$ is an affine algebraic group over $\Cset$ and $M_+, M_-$ and $Q$ are algebraic subgroups such that $M_+ \cap M_-$ is connected, then  every non-empty intersection $Y_n(c)$ in \eqref{eq-DQQD} is irreducible and is thus a single $\T$-leaf of $\pi_{\sY_n}$ in $Y_n$.
\eth

\begin{proof}
Parts (a) and (b) follow immediately from \prref{pr-D} and \leref{le-transition} by 
relabeling
$g \bowtie h$ by $(c_1, \ldots, c_n)$ and $(g\bowtie h)_{n+1}$ by $c_{n+1}$. The last part of \thref{th-Zn} follows from 
\reref{re-MM-connected}.
\end{proof}

We now turn to the case of $(X_n, \pi_{\sX_n})$, where $n \geq 2$.
%Note that, as $Q$ is a Poisson Lie subgroup of $(G, \piG)$, for each $n \geq 1$, one has the Poisson embedding
%\begin{equation}\lb{eq-D-n-n}
%(X_n, \, \pi_{\sX_n}) \lrw (D_{n+1}, \, \pi_{\sD_{n+1}}), \;\;
%[g_1, \, g_2,\, \ldots g_n]_{\sX_n} \longmapsto [e, \, g_1,\, g_2, \,\ldots,\, g_n]_{\sD_{n+1}}.
%\end{equation}
Let 
\[
\mu_{\sX_n}: X_n \lrw G, \;\; [g_1, g_2, \ldots, g_n]_{\sX_n} \longmapsto g_1g_2 \cdots g_n, \hs g_j \in G.
\]

\bas{as-GQQG}
There exists a direct sum decomposition $\q = \c + \q^\perp$ such that 

1) every $(K, Q)$-double cosets in $G$, where $K \in \{M_+, Q\}$, contains elements in $N_G(\c)$;

2) $p_{\c}\left(\m_+^\perp \cap \Ad_g \q\right) \cap \Ad_h p_{\c}\left(\m_+^\perp\cap \Ad_k \q\right)=0$ for 
all $g, h, k\in N_G(\c)$,

\noindent
where again $p_\c: \g \to \c$ is the projection with respect to the decomposition $\g = \c + \c^\perp$.

\eas

\bth{th-Dn} Under \asref{as-GQQG}, the following holds for every $n \geq 2$:

(a) Every non-empty intersection 
\begin{equation}\lb{eq-MQQM}
X_n(c) = \left((M_+c_1Q)\times_Q (Qc_2Q) \times_Q \cdots \times_Q (Qc_nM_+)\right) \cap \mu_{\sX_n}^{-1}(M_-c_{n+1}M_-) \subset X_n,
\end{equation}
where $c = (c_1, \ldots, c_{n+1})$ with $(c_1, \ldots, c_n) \in (N_G(\c))^{n}$ and $c_{n+1} \in G$, 
is transversal, and their connected components are precisely all the $\T^2$-leaves of $\pi_{\sX_n}$ in $X_n$;

(b) The leaf stabilizer for 
every $\T^2$-leaf in $X_n(c)$ in \eqref{eq-MQQM} is $(p_\t \oplus p_\t)(V_c)$, 
where
\begin{align*}
V_c &= \left\{\left(x_+, \, x_-, \, z_+, \; \Ad_{c_{n+1}^{-1}}(x_-)\right): \;\; x_+\in \m_+\cap \Ad_{c_1} \q, \;\; 
x_- \in \m_- \cap \Ad_{c_{n+1}} \m_-, \right.\\
& \hspace{4cm}\; \left. z_+  \in \m_+\cap \Ad_{c_{n}^{-1}} \q, \;\;
p_\c(x_+) = \Ad_{c_1c_2 \cdots c_n} (p_\c(z_+))\right\}. 
\end{align*}

\noindent
Moreover, if $G$ is an affine algebraic group over $\Cset$ and $M_+, M_-$ and $Q$ are algebraic subgroups such that $M_+ \cap M_-$ is connected, then  every non-empty intersection $X_n(c)$ in \eqref{eq-MQQM} is irreducible and is thus a single $\T^2$-leaf of $\pi_{\sX_n}$ in $X_n$.
\eth

\begin{proof}
Consider the diffeomorphism
$J_{\sX_n}: X_n \to (G/Q)^{n-1} \times G$  given by
\[
J_{\sX_n}([g_1, g_2, \cdots, g_n]_{\sX_n}) = \left(g_1Q, \; g_1g_2Q, \, \ldots, g_1g_2 \cdots g_{n-1}Q, \; g_1g_2 \cdots g_n\right),
\hs g_j \in G.
\]
We will again study $\pi_{\sX_n}$ by studying the Poisson structure $J_{\sX_n}(\pi_{\sX_n})$ on $(G/Q)^{n-1} \times G$.
Let $\lam$ be the left  action of $G^{n+1}$ on  the product manifold $(G/Q)^{n-1} \times G$ by 
\[
(g_1, \,\ldots, \,g_n, \,g_{n+1}) \cdot (h_1Q, \; \ldots, \; h_{n-1}Q, \; h_n) = 
(g_1h_1Q, \; \ldots, \; g_{n-1}h_{n-1}Q, \; g_nh_ng_{n+1}^{-1}),
\]
where $g_j, h_k \in G$,
and define
the direct sum quasitriangular $r$-matrix $r^{\la n+1\ra}$ on $\g^{n+1}$ by
\begin{equation}\lb{eq-r-modify}
r^{\la n+1\ra} = \begin{cases} (r^{(n)}, \, 0) + (0, \, -r), & \;\; \mbox{if} \;\; n = 2m+1 \;  \mbox{is odd},\\
(r^{(n)}, \, 0) + (0, \, r^{21}), & \;\; \mbox{if} \;\; n = 2m \; \mbox{is even}.\end{cases}
\end{equation}
Then the homogeneous quadruple 
$\displaystyle \left(G^{n+1}, \, r^{\la n+1 \ra}, \; (G/Q)^{n-1} \times G, \; \lam\right)$ is strongly admissible.
It is proved in \cite[$\S$8]{Lu-Mou:mixed} that, as Poisson structures on $(G/Q)^{n-1} \times G$,
\begin{equation}\lb{eq-J-Dn-pi}
J_{\sX_n}(\pi_{\sX_n}) = -\lam\left(r^{\la n+1 \ra}\right).
\end{equation}
Define the Lie subgroups $M_+^{\la n+1\ra} \subset G^{n+1}$ and $M_-^{\la n+1\ra} \subset G^{n+1}$ respectively by
\begin{align*}
& M_+^{\la n+1\ra} = M_+^{(n)} \times M_-  \hs \mbox{and} \hs M_-^{\la n+1\ra} = M_-^{(n)} \times M_+ \hs
\mbox{if} \; n = 2m \; \mbox{is even},\\
& M_+^{\la n+1\ra} = M_+^{(n)} \times M_+  \hs \mbox{and} \hs M_-^{\la n+1\ra} = M_-^{(n)} \times M_- \hs
\mbox{if} \; n = 2m +1\; \mbox{is odd}.
\end{align*}
Then the pair $\left(M_+^{\la n+1\ra}, \; M_-^{\la n+1\ra}\right)$ is $r^{\la n+1 \ra}$-admissible, and one thus has the
six-tuple 
\begin{equation}\lb{eq-six-Dn}
\left(G^{n+1}, \; r^{\la n+1 \ra}, \; (G/Q)^{n-1} \times G, \; \lam, \; M_+^{\la n+1\ra}, \; M_-^{\la n+1\ra}\right).
\end{equation}
Same as for \thref{th-Zn}, we will apply \thref{th-main} to the six-tuple in \eqref{eq-six-Dn} and use the diffeomorphism 
$J_{\sX_n}: X_n \to (G/Q)^{n-1} \times G$ to prove \thref{th-Dn}. 
Note that by identifying 
\[
\T^2 \; \stackrel{\sim}{\lrw} \; M_+^{\la n+1\ra} \cap  M_-^{\la n+1\ra}, \; \; (t_1, t_2) \longmapsto
(t_1, \, t_1, \, \ldots, \, t_1, \, t_2), \hs t_1, t_2 \in \T
\]
the diffeomorphism $J_{\sX_n}: X_n \to (G/Q)^{n-1} \times G$ is $\T^2$-equivariant.
For notational simplicity, set again $\widetilde{M}_\pm =M_\pm^{\la n+1\ra}$ and let $\tilde{\m}_\pm$ be their respective
Lie algebras.

Assume first that $n = 2m$ is even with $m \geq 1$. 
Let $(\widetilde{\O}_+, 
\widetilde{\O}_-)$
be an arbitrary pair of $\widetilde{M}_+$ and $\widetilde{M}_-$-orbits in $(G/Q)^{n-1} \times G$. By
1) of \asref{as-GQQG}, there exist 
\begin{align*}
g& =(g_1, \; e, \;g_3,\; e,\; g_5,\; \ldots, \; \;e, \;g_{2m-1},\; g_{2m}, \; e) \in (N_G(\c))^{n-1} \times G \times G,\\
h &= (e, \;h_2,\; e, \;h_4,\; e, \;\ldots,\;h_{2m-2}, \; e, \;h_{2m}, \; e)\in (N_G(\c))^{n+1},
\end{align*} 
such that $y_+:=\underline{g} \in \widetilde{\O}_+$ and $y_-:=\underline{h} \in \widetilde{\O}_-$, where 
$\underline{a} = (a_1Q, \ldots, a_{n-1}Q, a_na_{n+1}^{-1}) \in (G/Q)^{n-1} \times G$ for
$a = (a_1, \ldots, a_{n+1}) \in G^{n+1}$. 
Let
\[
c=  (g_1, \; h_2,\; g_3, \; h_4, \; \ldots,  \;g_{2m-1}, \;h_{2m}, \; g_{2m}) \in (N_G(\c))^n \times G.
\]
By \leref{le-as-homog},  
$\displaystyle \left(\tilde{\m}_+^\perp \oplus \tilde{\m}_-^\perp\right)\cap \l_{y_+, y_-}$ consists of  
elements $a \in \g^{n+1} \oplus \g^{n+1}$ of the form
\begin{equation}\lb{eq-a-2}
a = (x_1 + x_1^\prime,  \; \ldots, \; x_{n-1} + x_{n-1}^\prime, \; \Ad_{g_{n}} (x_{n}), \; x_{n}, \;
z_1 + z_1^\prime, \; \ldots, \; z_{n-1} + z_{n-1}^\prime, \; \Ad_{h_{n}} (z_{n}), \; z_{n}),
\end{equation}
where, by writing $g = (g_1, g_2, \ldots, g_{n+1})$ and $h = (h_1, h_2, \ldots, h_{n+1})$,
\[
x_j, \; z_j \in \c, \hs x_j^\prime \in \Ad_{g_j}\q^\perp, \; z_j^\prime \in \Ad_{h_j}\q^\perp, \hs x_j = \Ad_{g_jh_j^{-1}} z_j, \hs j = 1, \ldots, n-1,
\]
$x_1 + x_1^\prime \in \m_+^\perp \cap \Ad_{g_1}\q$,
$\;x_{n} \in \m_-^\perp \cap \Ad_{g_{n}^{-1}}\m_-^\perp$, $\; z_{n} \in \m_+^\perp$, $\;z_{n-1} + z_{n-1}^\prime = \Ad_{h_{n}}(z_{n})$,
and
\[
x_{2j} = x_{2j+1}, \;\;\; x_{2j}^\prime = x_{2j+1}^\prime, \hs 
z_{2j-1} = z_{2j}, \; \;\; z_{2j-1}^\prime = z_{2j}^\prime,\hs j = 1, \ldots, m-1.
\]
It follows that
\[
x_1 = \Ad_{g_1 h_2g_3h_4\cdots g_{2m-1}h_{2m}}(z_{2m}) \in p_\c(\m_+^\perp\cap \Ad_{c_1}\q) 
\cap \Ad_{c_1c_2 \cdots c_n} p_\c(\m_+^\perp \cap \Ad_{c_n^{-1}}\q).
\]
By 2) of \asref{as-GQQG}, $x_1=0$, and it follows that $x_j = 0$ for every $j = 1, \ldots, n-1$. By \eqref{eq-delOO-3}, 
$\delta_{\widetilde{\sO}_+, \widetilde{\sO}_-} = 0$. 
Let $p_{\t \oplus \t}$ be the 
projection from $\tilde{\m}_+ \oplus \tilde{\m}_-$ to $\tilde{\m}_+ \cap \tilde{\m}_- \cong \t \oplus \t$.
An argument similar to that in the proof of \thref{th-Zn} also shows that 
\[
p_{\t \oplus \t} \left((\tilde{\m}_+ \oplus \tilde{\m}_-) \cap \l_{y_+, y_-}\right) = (p_\t \oplus p_\t) (V_c).
\]
Moreover, 
an argument similar to that in the proof of \leref{le-transition} shows that 
$J_{\sX_n}^{-1}\left(\widetilde{\O}_+ \cap \widetilde{\O}_-\right) = X_n(c)$. 
\thref{th-Dn} is now a consequence of \thref{th-main}.
The case of $n = 2m+1$ is odd is proved similarly. 

\end{proof}

\sectionnew{Applications to Poisson structures related to flag varieties}\lb{sec-flags}
After reviewing the standard complex semi-simple Lie groups,  we prove
\thref{th-Fn} - \thref{th-Dbbn} stated in $\S$\ref{subsec-GBGB} on the four series of Poisson manifolds 
given in \eqref{eq-series}.

\subsection{Standard complex semi-simple Poisson Lie groups}\lb{subsec-rst}
Let $\g$ be a complex semi-simple Lie algebra, and let $\lara_\g$ be a fixed symmetric non-degenerate invariant
bilinear form on $\g$. Fix also a choice $(\b, \b_-)$ of opposite Borel subalgebras of $\g$ and let
$\h = \b \cap \b_-$, a Cartan subalgebra of $\g$. Let $\Delta$ and $\Delta_+\subset\Delta$ be respectively the set of roots 
for the pairs $(\g, \h)$ and $(\b, \h)$, and let 
$\displaystyle \g = \h + \sum_{\alpha \in \Delta_+} \g_{\alpha} + \sum_{\alpha \in \Delta_+} \g_{-\alpha}$
be the corresponding root decomposition. Let $\displaystyle \n_+ = \sum_{\alpha \in \Delta_+} \g_{\alpha}$ and 
$\displaystyle \n_- =  \sum_{\alpha \in \Delta_+} \g_{-\alpha}$. Equip again the direct product Lie algebra $\gog$ with the
bilinear form $\lara_{\gog}$ as in \eqref{eq-bra-gog}.  Then 
\[
\l_{\rm st} = \{(x_+ + x_0, \; -x_0 + x_-): \;\;  x_\pm \in \n_\pm, \; x _0 \in \h\}
\]
is a Lagrangian subalgebra of $(\gog, \lara_{\gog})$ such that $\gog = \gdia + \l_{\rm st}$. 
The decomposition $\gog = \gdia + \l_{\rm st}$ is called a {\it standard} Lagrangian splitting  of
$(\gog, \lara_{\gog})$ and the element $r_{\rm st} \in \gotg$ such that $\l_{r_{\rm st}} = \l_{\rm st}$ a 
{\it standard (factorizable quasitriangular) $r$-matrix} on $\g$. 
Let $\{h_i\}_{i=1}^r$ be a basis of $\h$ such that $\la h_i, h_j\ra_\g = \delta_{ij}$
for $1 \leq i, j \leq r$, and let $E_\alpha \in \g_\alpha$ and $E_{-\alpha} \in \g_{-\alpha}$ be root vectors for 
$\alpha \in \Delta_+$ such that $\la E_\alpha, E_{-\alpha}\ra_\g = 1$. Then $r_{\rm st}$ is explicitly given by
\begin{equation}\lb{eq-r-st-0}
r_{\rm st} = \frac{1}{2} \sum_{i=1}^r h_i \ot h_i + \sum_{\alpha \in \Delta_+} E_{-\alpha} \otimes E_\alpha.
\end{equation}
It is clear from the definition in \eqref{eq-de-l-pm} that the Lie subalgebras $\f_-$
and $\f_+$ of $\g$ associated to $r_{\rm st}$ are respectively given by
$\f_- = \b_-$ and $\f_+ = \b_+$. It is also clear that $r_{\rm st}$ is factorizable and that
the symmetric bilinear form on $\g$ associated to $r_{\rm st}$
(see $\S$\ref{subsec-factorizable}) is precisely $\lara_\g$.

By \eqref{eq-r-(n)}, one also has the factorizable quasitriangular $r$-matrix  $r_{\rm st}^{(2)}$ on $\gog$. Explicitly,
\begin{align}\lb{eq-r-st-2}
r_{\rm st}^{(2)} &=\frac{1}{2} \sum_i\left((h_i, \, 0) \otimes (h_i, \, 0) -(0, h_i)\otimes (0, h_i) -(h_i, 0) \wedge (0, h_i)\right)\\
\nonumber
& \hs + \sum_{\alpha \in \Delta_+}\left((E_{-\alpha}\otimes E_\alpha, \, 0) - (0, \, E_\alpha \otimes E_{-\alpha}) -
(E_\alpha, 0) \wedge (0, E_{-\alpha})\right)\\
\nonumber
& = \frac{1}{2} \sum_i (h_i, h_i) \otimes (h_i, -h_i) + \sum_{\alpha \in \Delta_+}\left(
(E_\alpha, \, E_\alpha) \otimes (0, \, -E_{-\alpha}) + (E_{-\alpha}, \, E_{-\alpha}) \otimes (E_\alpha, \, 0)\right).
\end{align}
One checks directly (see also \cite[$\S$6]{Lu-Mou:mixed}) that  $r_{\rm st}^{(2)}$ coincides with the $r$-matrix on $\gog$ defined by the
Lagrangian splitting $\gog = \gdia + \l_{\rm st}$ (see the definition in \exref{ex-double-0}). The Lie subalgebras $\f_+$ and $\f_-$ of $\gog$ associated to $r^{(2)}_{\rm st}$ are then $\l_{\rm st}$ and $\gdia$, respectively.

Let $G$ be any connected complex Lie group with Lie algebra $\g$. The Poisson structure 
\[
\pist :=r_{\rm st}^L - r_{\rm st}^R
\]
on $G$ is called a
{\it standard multiplicative holomorphic Poisson structure} on $G$, and the pair $(G, \pist)$ a 
{\it standard complex semi-simple Poisson Lie group}.
One also has the {\it Drinfeld double} Poisson Lie group
$(G \times G, \, \Pi_{\rm st})$ of the Poisson Lie group $(G, \pist)$, where 
\[
\Pist = \left(r^{(2)}_{\rm st}\right)^L - \left(r^{(2)}_{\rm st}\right)^R.
\]

Let $B$ and $B_-$ be the Borel subgroups of $G$ with Lie algebras $\b$ and $\b_-$ respectively, and let
$T = B \cap B_-$, a maximal torus of $G$.
Let $Q = B$ for the Poisson Lie group $(G, \pist)$ and $Q = B \times B_-$ for the Poisson Lie group
$(G \times G, \, \Pi_{\rm st})$. Then Condition \eqref{eq-on-q-1} is satisfied, and one arrives at
the four quotient Poisson manifolds $(F_n, \, \pi_n), \, (\Fbb_n, \, \Pi_n), \, (\wF_n, \, \tilde{\pi}_n)$ and 
$(\wFbb_n, \, \wPi_n)$, as in \eqref{eq-series} in $\S$\ref{subsec-GBGB}.  

For  $(F_n, \, \pi_n)$ and $(\wF_n, \, \tilde{\pi}_n)$, take 
$M_+ = B$ and $M_- = B_-$. Then $M_+ \cap M_- = T$.  The projection $p_\t: \m_+ \oplus \m_-
\to \t = \h$ in \eqref{eq-pmt} is now given by
\begin{equation}\lb{eq-pmt-bb}
p_\t: \;\; \b \oplus \b_- \lrw \h, \;\;\; (x_0+x_+, \, y_0+y_-) \longmapsto  x_0 + y_0, \hs x_0, \, y_0 \in \h, \, x_+ \in \n, \, y_- \in \n_-.
\end{equation} 
For $(\Fbb_n, \Pi_n)$ and $(\wFbb_n, \wPi_n)$, take  
$M_+ = B \times B_-$ and $M_- = G_{\rm diag} = \{(g, ,g): g \in G\}$. Then $M_+ \cap M_- = T_{\rm diag} =
\{(t, t): t \in T\}$, 
The projection $p_\t: \m_+ \oplus \m_-
\to \t = \h_{\rm diag} \cong \h$ in \eqref{eq-pmt} in this case given by
\begin{equation}\lb{eq-pmt-bb-gdia}
\p_\t:\;\; (\b\oplus \b_-) \oplus \gdia \lrw \h, \;\;\; ((x_0+x_+, \, y_0+y_-), \, (x, \, x)) \longmapsto x_0+y_0,
\end{equation}
where $x_0, \, y_0 \in \h, \, x_+ \in \n, \, y_- \in \n_-$, and $x \in \g$.

We make some further preparation for the proofs  of \thref{th-Fn} - \thref{th-Dbbn}. 

\ble{le-wuv} For $u_1, \ldots, u_n, v_1, \ldots, v_n, w \in W$, the following are equivalent:

1) $(Bu_1Bu_2 \cdots Bu_nB) \cap (B_- v_1 B_- v_2 \cdots B_- v_nB_-wB) \neq \emptyset$;

2) $w \leq (v_1 \ast \cdots \ast v_n)^{-1} \ast u_1 \ast \cdots \ast u_n$.
\ele

\begin{proof} For $n = 1$, the equivalent between 1) and 2) is proved in 
\cite[Proposition 4.1]{Webster-Yakimov}. Assume that $n \geq 2$ and write $u = 
u_1 \ast \cdots \ast u_n$ and $v = v_1 \ast \cdots \ast v_n$.  
Suppose that 1) holds. As 
\[
Bu_1B \cdots Bu_nB = \bigsqcup_{x \in {\mathcal U}} BxB \hs \mbox{and} \hs
B_-v_1B_-v_2 \cdots B_-v_nB_- = \bigsqcup_{y \in {\mathcal V}} B_-yB_-,
\]
where ${\mathcal U} \subset \{x \in W: x \leq u\}$ and ${\mathcal V} \subset \{y \in W: y \leq v\}$, 
there exist $x \leq u$ and $y \leq v$ such that 
$(BxB) \cap (B_-yB_-wB) \neq \emptyset$, and hence $w \leq y^{-1} \ast x \leq v^{-1} \ast u$.
Conversely, if $w \leq v^{-1} \ast u$, then $(BuB) \cap (B_-vB_-wB) \neq \emptyset$. As
$BuB \subset Bu_1B \cdots Bu_nB$ and $B_-vB_- \subset B_-v_1B_- \cdots B_-v_nB_-$, 1) holds.
\end{proof}

\ble{le-C}
For any $u, v \in W$ and any conjugacy class $C$ in $G$, $(BuBB_- v B_-) \cap C \neq \emptyset$. 
\ele

\begin{proof} As $B_- \cap (BuB) \neq \emptyset$, one has $B_- \cap (BuBB_-) \neq \emptyset$, which implies that
$BB_- \subset BuBB_-$. Similarly, $BB_- \subset BB_-vB_-$. Thus
$B \subset B B_- \subset BuBB_- \subset BuBB_- v B_-$.
As $B \cap C \neq \emptyset$, one has $(BuBB_- v B_-) \cap C \neq \emptyset$. 
\end{proof}

For $n \geq 1$, consider 
$P_n:(G \times G)^n \to G^n$ given by  $P_n(g_1, k_1, \ldots, g_n, k_n) = (g_1, \ldots, g_n)$,
and the induced projections, both denoted as
$[P_n]$, from $\Fbb_n$ to $F_n$ and from $\wFbb_n$ to  $\wF_n$, i.e., 
\begin{equation}\lb{eq-Pn}
 P_n([g_1, k_1, \ldots, g_n, k_n]_{\sFbb_n})= [g_1, \ldots, g_n]_{\sF_n}, \hs 
P_n([g_1, k_1, \ldots, g_n, k_n]_{\swFbb_n}) =  [g_1, \ldots, g_n]_{\swF_n}.
\end{equation}
As $P_1: (G \times G, \, \Pist) \to (G, \pist)$ is Poisson (this follows, for example, from the fact that
$P_1(r_{\rm st}^{(2)}) = r_{\rm st}$), one knows that 
$[P_n]: (\Fbb_n, \, \Pi_n) \to (F_n, \, \pi_n)$ and $[P_n]: (\wFbb_n, \, \wPi_n) \to (\wF_n, \, \tilde{\pi}_n)$
are Poisson.  This observation will be used to show that \thref{th-Fn} is a special case of \thref{th-Fbbn}
(see $\S$\ref{subsec-proof-AB}) and that \thref{th-Dn-intro} a special case of \thref{th-Dbbn} (see $\S$\ref{subsec-proof-CD}).

\subsection{Proofs of \thref{th-Fn} and \thref{th-Fbbn}}\lb{subsec-proof-AB}
We first prove \thref{th-Fbbn} by applying \thref{th-Zn} to the Poisson Lie group  $(G \times G, \Pist)$ and by taking
$Q = M_+ = B \times B_-$
and $M_- = \Gdia$. As $M_+ \cap M_- = T_{\rm diag} :=\{(t, t): t \in \T\}$, the intersection $R^{\bfu, \bfv}_{w}$, whenever non-empty, is
smooth and connected and has dimension  $l(\bfu) +l(\bfv)- l(w)$. The fact that $R^{\bfu, \bfv}_w \neq \emptyset$
if and only if $w \leq (v_1 \ast \cdots \ast v_n)^{-1} \ast u_1 \ast \cdots \ast u_n$ follows directly from \leref{le-wuv}
(it is also proved in 
\cite[Proposition 3.32]{Victor:IMRN} using distinguished double subexpressions). This proves 1) of \thref{th-Fbbn}.
Letting $c = (u_1, v_1, \ldots, u_n, v_n)$ and $c_{n+1} = (w, e)$ in (b) of \thref{th-Zn}, one proves 
the first part of 2) of \thref{th-Fbbn}. Computing explicitly the subspace $V_c \subset \m_+ \oplus \m_-$ as
described in \thref{th-Zn} and using \eqref{eq-pmt-bb-gdia} for $p_\t: \m_+ \oplus \m_- \to \t = \h$, one sees that the 
leaf stabilizer of $\lambda_{\sFbb_n}$ in $R^{\bfu, \bfv}_w$ is  precisely $\h^{\bfu, \bfv}_w \subset \h$.
This proves 3) of \thref{th-Fbbn}.  
The second part of 2) follows from 3) by \thref{th-Zn}.
This finishes the proof of \thref{th-Fbbn}.

\thref{th-Dn-intro} is proved either similarly, by applying \thref{th-Zn} to the Poisson Lie group
$(G, \pist)$ and taking $Q = M_+ = B$ and $M_- = B_-$, or by the following observation: let 
\[
\widehat{F}_n :=(G \times B_-)\times_{(B \times B_-)} \cdots \times_{(B \times B_-)} (G \times B_-)/(B \times B_-)
\subset \Fbb_n.
\]
By \leref{le-m-sub}, $G \times B_-$ is a Poisson Lie subgroup of 
$(G \times G, \Pist)$. Thus $\widehat{F}_n$ is a Poisson submanifold of $(\Fbb_n, \Pi_n)$. 
It is then clear that the Poisson morphism
$[P_n]: (\Fbb_n, \, {\Pi}_n) \to (F_n, \, {\pi}_n)$ in \eqref{eq-Pn} restricts to a Poisson isomorphism
from $(\widehat{F}_n, \, {\Pi}_n)$ to $(F_n, \, {\pi}_n)$. \thref{th-Fn} now follows by
applying \thref{th-Fbbn} to the Poisson submanifold $(\widehat{F}_n, \, {\Pi}_n)$ of $(\Fbb_n, \, {\Pi}_n)$.

\subsection{Proofs of \thref{th-Dn-intro} and \thref{th-Dbbn}}\lb{subsec-proof-CD}
We first prove \thref{th-Dbbn} by applying \thref{th-Dn} to the Poisson Lie group  $(G \times G, \Pist)$ and taking $Q = M_+ = B \times B_-$
and $M_- = \Gdia$. As $M_+ \cap M_- = T_{\rm diag}$, the intersection $R^{\bfu, \bfv}_{\sC}$, whenever non-empty, is
smooth and connected, and  
\[
\dim R^{\bfu, \bfv}_{\sC} =\dim ((B \times B_-)(\bfu, \bfv)(B \times B_-)) + 
\dim \mu_{\swFbb_n}^{-1} (\Omega_{\sC}) - \dim \wFbb_n 
= l(\bfu) + l(\bfv) + \dim C + \dim T.
\]
Let $\bfu = (u_1, \ldots, u_n), \bfv = (v_1, \ldots, v_n) \in W^n$, and $C \in {\mathcal C}$ be arbitrary.
By definition,
\[
\mu_{\swFbb_n}^{-1}(\Omega_{\sC}) = \{[g_1, k_2, \ldots, g_n,k_n]_{\swFbb_n}: \, (g_1g_1\cdots g_n)(k_1k_2 \cdots k_n)^{-1} \in C\}.
\]
By \leref{le-C}, $(Bu_1Bu_2 \cdots Bu_nBB_-v_n^{-1}\cdots B_-v_2^{-1}B_-v_1^{-1}B_-) \cap C \neq \emptyset$, so
$R^{\bfu, \bfv}_{\sC} \neq \emptyset$. This proves 1) of \thref{th-Dbbn}.
Letting $c = (u_1, v_1, \ldots, u_n, v_n)$ and $c_{n+1} = (c, e)$ in (b) of \thref{th-Dn}, where $c \in C$ is arbitrary, one 
sees that the $R^{\bfu, \bfv}_{\sC}$'s are precisely the 
$T^2$-leaves of $\wPi_n$ in $\wFbb_n$, where $T^2$-acts on $\wFbb_n$ by
\begin{equation}\lb{eq-TT-Dbbn}  (t_1, t_2) \cdot [g_1, \,k_1, \,\ldots, \, g_n, \, k_n]_{\swFbb_n} = 
[tg_1, \,tk_1, \,g_2, \,k_2,  \,  \ldots, \, g_{n-1}, \, k_{n-1}, \,g_nt_2^{-1}, \,k_nt_2^{-1}]_{\swFbb_n}.
\end{equation}
Computing explicitly the subspace $V_c \subset \m_+ \oplus \m_-$ as
described in \thref{th-Dn} and using \eqref{eq-pmt-bb-gdia} for $p_\t: \m_+ \oplus \m_- \to \t = \h$, one sees that the 
leaf stabilizer for the $T^2$-action on $R^{\bfu, \bfv}_{\sC}$ is given by
\[
(\h \oplus \h)^{\bfu, \bfv} :=\{(u(x) + v(y), \, x + y): \;x, y \in \h\},
\]
where $u = u_1 \ast \cdots \ast u_n$ and $v = v_1 \ast \cdots \ast v_n$. Note that 
$(\h \oplus 0) + (\h \oplus \h)^{\bfu, \bfv} = \h \oplus \h$, the action
$\lambda_{\swFbb_n}$ of $T \cong T \times \{e\} \subset T \times T$ on $\wFbb_n$ given in
\eqref{eq-T-Dbbn} is also full, so each $R^{\bfu, \bfv}_{\sC}$ is a single $T$-leaf for the action 
$\lambda_{\swFbb_n}$. Moreover, the leaf-stabilizer of $\lambda_{\swFbb_n}$ is $R^{\bfu, \bfv}_{\sC}$ is
\[
(\h \oplus 0) \cap ((\h \oplus \h)^{\bfu, \bfv}) = \{u(x) - v(x): \; x \in \h\} = \h^{\bfu, \bfv}.
\]
This finishes the proof of \thref{th-Dbbn}.

\thref{th-Dn-intro} is proved either similarly, by applying \thref{th-Dn} to the Poisson Lie group
$(G, \pist)$ and taking $Q = M_+ = B$ and $M_- = B_-$, or by the following observation: let 
\[
F^\vee_n :=((G \times B_-)\times_{(B \times B_-)} \cdots \times_{(B \times B_-)} (G \times B_-)\times_{(B \times B_-)}
(G \times G))\cap  \mu_{\swFbb_n}^{-1}(\Gdia) \subset \wFbb_n.
\]
By \thref{th-Dbbn}, $F^\vee_n$ 
is a Poisson submanifold of $(\wFbb_n, \wPi_n)$. It is clear that the Poisson morphism
$[P_n]: (\wFbb_n, \, \wPi_n) \to (\wF_n, \, \tilde{\pi}_n)$ in \eqref{eq-Pn} restricts to a Poisson isomorphism
from $(F^\vee_n, \, \wPi_n)$ to $(\wF_n, \, \tilde{\pi}_n)$. \thref{th-Dn-intro} now follows by
applying \thref{th-Dbbn} to the Poisson submanifold $(F^\vee_n, \, \wPi_n)$ of $(\wFbb_n, \, \wPi_n)$.

%\subsection{Extended double Bruhat cells}

\subsection{Other examples}\lb{subsec-other}
The results in
$\S$\ref{sec-orbits} and $\S$\ref{sec-mixed} 
can be used to give a unified approach to many other examples, old or new, of $T$-Poisson manifolds related to real or complex
semi-simple Lie groups.
We give two such example here, leaving other examples to be treated elsewhere.

\bex{ex-theta-twisted} Let $(G, \pist)$ be again a standard complex semi-simple Lie group as in 
$\S$\ref{subsec-rst}, and let
$\theta \in {\rm Aut}(G)$ be such that $\theta(T) = T$ and $\theta(B) = B$, and denote by the same letter the induced
automorphism on $\g$. Let $\lam$ be the left action of 
$G \times G$ on $G$ given by
\begin{equation}\lb{eq-GG-G}
(g_1, g_2)\cdot_\theta g = g_1 g \theta(g_2)^{-1}, \;\;\; g_1, g_2, g \in G.
\end{equation}
Orbits of $\Gdia \subset G \times G$ on $G$ under $\lambda$ are called {\it $\theta$-twisted conjugacy classes} of $G$.
Let $r = r_{\rm st}^{(2)}$ be the $r$-matrix on $\gog$ defined by 
the Lagrangian splitting $\gog = \gdia + \l_{\rm st}$ (see \eqref{eq-r-st-2}). 
As the stabilizer subalgebras for $\lam$ are Lagrangian with respect to $\lara_{\gog}$, the quadruple
$(G\times G, \, r, \,G, \, \lambda)$ is strongly admissible. The Poisson structure $\pi_\theta = -\lambda(r)$ is
studied in \cite{Lu:twisted}.

Recall that the Lie subalgebras
$\f_-$ and $\f_+$ of $\gog$ associated to $r$ is $\f_- = \gdia$ and $\f_+ = \l_{\rm st}$. Let $M_- = \Gdia$ and let $M_+ = B \times B_-$, so that $M_- \cap M_+ =T_{\rm diag} = \{(t, t): t \in T\}$ which we
identify with $T$. 
Note that $M_-$-orbits in $G$ are precisely the $\theta$-twisted conjugacy classes in $G$, and each $M_+$-orbit
is of the form $BwB_-$ for a unique $w \in W$.
It is shown in \cite[$\S$2.4]{Chan-Lu-To} that for a $\theta$-twisted conjugacy class $C$, there is a unique
element $m_{\sC} \in W$ such that for $w \in W$, 
$C \cap (BwB_-) \neq \emptyset$ if and only if $w \leq m_{\sC}$. 
Again as the stabilizer subalgebras of $\lam$ are Lagrangian with respect to $\lara_{\gog}$, it follows trivially
from \prref{pr-delOO}
that $\delOO = 0$ for every $\O_- = C$ and $\O_+ = BwB_-$. Thus the $T$-leaf
decomposition of $(G, \pi_\theta)$ is given by
$G = \bigsqcup_{C, w} C \cap (BwB_-)$, 
where $C$ is a $\theta$-twisted conjugacy class in $G$ and $w \in W$ such that $w \leq m_{\sC}$. 
Denote by
$\t_{C, w}$ the leaf stabilizer of $C \cap (BwB_-)$ in $\h$. Pick any $g \in C$ and let $y_- = g$ and $y_+ = \dot{w}$. 
Then 
\[
\l_{y_+, y_-} = \q_{y_+} \oplus \q_{y_-} = \Ad_{(\dot{w}, e)} \g_\theta \oplus \Ad_{(g, e)} \g_\theta \subset 
(\g \oplus \g) \oplus (\g \oplus \g), 
\]
where $\g_\theta = \{(\theta(x), x): x \in \g\}$.
It follows from \eqref{eq-ty-0} that
$\t_{C, w} = {\rm Im} (w \theta + 1)$, a result obtained in \cite{Lu:twisted}.
\eex

\bex{ex-real-BS} Let $G$ be a connected complex semi-simple Lie group with Lie algebra $\g$, and let
${\rm Im} \lara_\g$ be the imaginary part of Killing form $\lara_\g$ of $\g$. Let
$G = KAN$ be an Iwasawa decomposition of $G$, and let $\k, \a$, and $\n$ be respectively the Lie algebras of $K$, $A$, and $N$.
Regarding $(\g, {\rm Im}\lara_\g)$ as a real quadratic Lie algebra, one has the Lagrangian splitting
$\g = \k + (\a + \n)$ of $(\g, {\rm Im}\lara_\g)$. Let $\piG$ be the real analytic Poisson structure on $G$
given by $\piG = (r_0)^L - (r_0)^R$, where $r_0 \in \g \otimes_{{\mathbb R}} \g$ is the quasitriangular $r$-matrix
on $\g$ defined by the Lagrangian splitting $\g = \k + (\a + \n)$. Then $K$ is a Poisson Lie subgroup of $(G, \piG)$, and 
for each integer $n \geq 1$, one has the quotient space
${\mathcal P}_n = G \times_K \cdots \times_K G/K$
with the quotient Poisson structure $\pi_{\scriptscriptstyle{\mathcal{P}}_n}$ on ${\mathcal P}_n$. 
One can regard  ${\mathcal P}_n$ as the space of $n$-gons in the Riemannian symmetric space $G/K$
(see \cite{amw:linear, KLM}). Taking $M_+ = K$ and $M_- = AN$ so that $M_+ \cap M_- = \{e\}$,
it follows from \thref{th-Zn} and the decomposition $G = KAK$ that the symplectic leaves of $\pi_{\scriptscriptstyle{\mathcal{P}}_n}$ in ${\mathcal P}_n$ are of the forms $(Ka_1K)\times_K \cdots \times_K (Ka_nK)/K$
with $a=(a_1, \ldots, a_n) \in A^n$, the space of $n$-gons with fixed {\it side lengths} $a$. 
\eex

\end{document}